\documentclass[11pt,twoside,leqno]{aomamlt2e} 
  \pageno{731}
\received{January 15, 2001}
\revised{July 22, 2004}
\newtheorem{theorem}{Theorem}[section]
\newcounter{list_count}
{\end{list}}

\numberwithin{equation}{Subsec}

\newtheorem{corollary}[theorem]{Corollary}
\newtheorem{lemma}[theorem]{Lemma}
\newtheorem{proposition}[theorem]{Proposition}
\theorembodyfont{\upshape}
\newtheorem{definition}{{\it Definition}}[section]
\newtheorem{remark}{{\it Remark}}[section]
\newtheorem{example}{{\it Example}}[section]
\newcommand\hfq{\hfill\qed}
\newcommand\bup{\vglue-18pt\phantom{up}}

\newcommand\epsil{\varepsilon} 
\newcommand\older{H{\hskip1.5pt\rm \"{\hskip-7pt\it o}}lder}
\begin{document}

\currannalsline{164}{2006} 

\title{Global hyperbolicity of renormalization \\
for $C^r$ unimodal mappings}

\shorttitle{Global Hyperbolicity of Renormalization}

 \acknowledgements{Financially supported by CNPq Grant 301970/2003-3.
 
\hglue2.5pt$^{\textstyle\ast\ast}$Financially supported by CNPq Grant 304912/2003-4 and
Faperj Grant E-26/152.189/ 2002.

 \hglue-4pt$^{\textstyle\ast\!\ast\!\ast}$Financially supported by Calouste Gulbenkian Foundation, 
PRODYN-ESF, POCTI\hfill\break and POSI by FCT and Minist\'erio da CTES, and
CMUP.
}
\twoauthors{Edson~de~Faria$^{\textstyle\ast}$, Welington de Melo$^{\textstyle\ast\ast}$}{Alberto
Pinto\raise2.75pt\hbox{${\textstyle\ast\ast}$}}

\institution{Instituto de Matem\'atica e Estat\'{\i}stica, Universidade de S\~ao
Paulo,\\ S\~ao Paulo SP,
Brazil
\\
\email{edson@ime.usp.br}\\
\vglue-9pt
I.M.P.A., Rio de Janeiro, Brazil
\\
\email{demelo@impa.br}\\
\vglue-9pt
Universidade do Porto, Porto, Portugal
\\
\email{aapinto@fc.up.pt}}

\shortname{Edson~de~Faria, Welington de Melo, and Alberto
Pinto}

\centerline{\bf Abstract}
\vglue6pt

In this paper we extend M.~Lyubich's recent results on the global
hyperbolicity of renormalization of quadratic-like germs to the space
of $C^r$ unimodal maps with quadratic critical point.
We show that in this space the bounded-type limit sets of the
renormalization operator have an invariant hyperbolic structure
provided $r \ge 2+\alpha$ with $\alpha$ close to one.
As an intermediate step between Lyubich's results and ours,
we prove that the renormalization operator is hyperbolic
in a Banach space of real analytic maps.
We  construct the local stable manifolds and prove that they form a
continuous lamination whose leaves are $C^1$ codimension one, Banach
submanifolds of the ambient space, and whose holonomy is $C^{1+\beta}$
for some $\beta>0$.  We also prove that the global stable sets are
$C^1$ immersed (codimension one) submanifolds as well, provided $r \ge
3+\alpha$ with $\alpha$ close to one.
As a corollary, we deduce that in generic,
one-parameter families of $C^r$ unimodal maps,
the set of parameters corresponding
to infinitely renormalizable maps of
bounded combinatorial type is a Cantor set with Hausdorff dimension  
less than one.\footnote{There is a list of symbols used in this paper, before the references, for the
convenience of the reader.}
 
\medbreak\centerline{\bf Table of Contents}
\def\sni#1{\smallbreak\noindent{#1}\hskip5pt}
\def\ssni#1{\vglue-1pt\noindent\hskip14pt{#1}}
\def\tni#1{\vglue-1pt\noindent\hskip38pt{#1}}
\def\oni#1{\smallbreak\noindent\hskip-5pt{#1}\hskip5pt}
\sni{1.} Introduction

\sni{2.}  Preliminaries and statements of results

\ssni{2.1.} Quadratic unimodal maps

\tni{2.1.1.}  The Banach spaces $\Bbb A^r$

\tni{2.1.2.}  The Banach spaces $\Bbb B^r$

\ssni{2.2.}  The renormalization operator

\ssni{2.3.}  The limit sets of \pagebreak renormalization

\ssni{2.4.}  Hyperbolic basic sets

\ssni{2.5.}  Hyperbolicity of renormalization

\sni{3.}  Hyperbolicity in a Banach space of real analytic maps

\ssni{3.1.}  Real analyticity of the renormalization operator

\ssni{3.2.}  Real analytic hybrid conjugacy classes

\ssni{3.3.}  Hyperbolic skew-products

\ssni{3.4.} Skew-product renormalization operator

\ssni{3.5.}  Hyperbolicity of the renormalization operator

\sni{4.}  Extending invariant splittings

\ssni{4.1.}  Compatibility

\sni{5.}  Extending the invariant splitting for renormalization

\ssni{5.1.}  H\"older norms and $L$-operators

\ssni{5.2.}  Bounded geometry

\ssni{5.3.}  Spectral estimates

\sni{6.}  The local stable manifold theorem

\ssni{6.1.}  Robust operators

\ssni{6.2.}  Stable manifolds for robust operators

\ssni{6.3.}  Uniform bounds

\ssni{6.4.}  Contraction towards the unstable manifolds

\ssni{6.5.}  Local stable sets

\ssni{6.6.}  Tangent spaces

\ssni{6.7.}  The main estimates

\ssni{6.8.}  The local stable sets are graphs

\ssni{6.9.}  Proof of the local stable manifold theorem

\sni{7.}  Smooth holonomies

\ssni{7.1.}  Small holonomies for robust operators

\sni{8.}  The renormalization operator is robust

\ssni{8.1.}  A closer look at composition

\ssni{8.2.}  Checking properties {\bf B2} and {\bf B3}

\ssni{8.3.}  Checking property {\bf B4}

\ssni{8.4.}  Checking properties {\bf B5} and {\bf B6}

\ssni{8.5.}  Proof of Theorem 8.1

\ssni{8.6.}  Proof of the hyperbolic picture

\tni{8.6.1.}  Proof of Theorem 2.5

\tni{8.6.2.}  Proof of Corollary 2.6

\sni{9.}  Global stable manifolds and one-parameter families

\ssni{9.1.}  The global stable manifolds of renormalization

\ssni{9.2.}  One-parameter families

\oni{10.}  A short list of symbols

\smallbreak\noindent References

\section{Introduction}
\label{sec:intro}

In 1978, M.~Feigenbaum \cite{feigen} and independently P.~Coullet and  
C.~Tresser
\cite{coullet} made a startling discovery concerning certain rigidity  
properties in
one-dimensional dynamics.
While analysing the transition between simple and
``chaotic'' dynamical behavior in ``typical'' one-parameter families of
unimodal maps -- such as the quadratic family $x\mapsto \lambda x(1-x)$  
--
they recorded the parameter values $\lambda_n$ at which successive
period-doubling bifurcations occurred in the family and found a
remarkable universal scaling law, namely
$$
\frac{\lambda_n-\lambda_{n-1}}{\lambda_{n+1}-\lambda_n}\;\to\;
4.669\dots
\ .
$$
They also found universal scalings within the geometry of the  
post-critical
set of the limiting map corresponding to the parameter
$\lambda_\infty=\lim{\lambda_n}$ (cf.\ the work of E.~Vul,
Ya.\ Sinai and K.~Khanin \cite{VKS}).
In an attempt to explain these phenomena, they introduced a certain  
nonlinear
operator acting on the space of unimodal maps -- the so-called {\it  
period
doubling operator}. They conjectured that the period-doubling operator  
has a
unique fixed point which is hyperbolic with a one-dimensional unstable
direction. They also conjectured that the universal constants they
found in their experiments are the eigenvalues of the derivative of the
operator at the fixed point.

A few years later (1982) this conjecture was confirmed by O.~Lanford
\cite{lanf} through a computer assisted proof. Working in a cleverly
defined Banach space of real analytic maps and using rigorous
numerical analysis on the computer,
Lanford  established at once the existence and
hyperbolicity of the fixed point of the period-doubling operator.
Subsequent work by M.~Campanino and\break H.~Epstein \cite{epstein} (also
Campanino et al. \cite{camp} and Epstein \cite{epst})
established the existence (but neither uniqueness nor hyperbolicity)
of the fixed point without essential help from the computer.

It was soon realized by Lanford and others that the period-doubling  
operator
was  just a restriction of another operator acting on the space of
unimodal maps -- the {\it renormalization operator} -- whose dynamical
behavior is much richer. The hopes were high that the iterates of this
operator would reveal the small scale geometric properties of the  
critical
orbits of many interesting one-dimensional systems. Hence, the initial
conjecture was generalized to the following.

\demo{\scshape Renormalization Conjecture}
{\it The limit set of the renormalization operator in the space
of maps of bounded combinatorial type is a hyperbolic Cantor set where  
the
operator acts as the full shift in a finite number of symbols.}
\Enddemo

(For a precise formulation of what is meant by bounded combinatorial  
type, see
\S \ref{sub:renorm} below.)

In the path towards a proof of this conjecture, several new ideas were\break
developed in the last 20 years by a number of mathematicians, especially\break
D.~Sullivan, C.~McMullen and M.~Lyubich. Among the deepest in
Dynamical Systems, these ideas have the complex dynamics of  
quadratic-like
maps (in the sense of Douady and Hubbard \cite{douady}) as a common  
thread.
Sullivan proved in \cite{sul} that
all limits of renormalization are quadratic-like maps with a definite
modulus. Then, constructing certain Teichm\"uller spaces
from quadratic-like maps and using a substitute of Schwarz's lemma in  
these
spaces, Sullivan established the existence of
horseshoe-like limit sets for renormalization.
Later, using a different approach based on Mostow rigidity,
McMullen \cite{mcmtwo} gave another proof of this result and went  
further by
showing that the convergence (in the $C^0$ sense) towards the limit set  
is
exponential.

The final breakthrough came with the work of Lyubich \cite{lyubich}.
He endowed the space of {\it germs} of quadratic-like maps (modulo  
affine
conjugacies) with a very subtle complex structure, showing that the
renormalization operator is complex-analytic with respect to such  
a structure.
In Lyubich's space, the stable sets of maps in the limit set of
renormalization coincide with the very {\it hybrid classes} of such  
maps, and
inherit a natural structure making them (complex codimension one)  
analytic
submanifolds. Combining McMullen's rigidity of towers with Schwarz's  
lemma in
Banach spaces, Lyubich proved exponential contraction along such stable
leaves. To obtain expansion in the transversal directions to such  
leaves at
points of the limit set, Lyubich argued by contradiction: if expansion  
fails,
then one can find a map in the limit set whose orbit under  
renormalization is
slowly shadowed by another orbit (the {\it small orbits theorem}, page  
323 of
\cite{lyubich}). This however contradicts another theorem of his,  
namely the
combinatorial rigidity theorem of \cite{lyubich2}. It follows that the
limit set is indeed hyperbolic in the space of germs.
Based on this result of Lyubich and using the real and complex bounds  
given by
Sullivan, we prove in   Theorem~\ref{lyubhyp} that
the attractor (for bounded combinatorics) is hyperbolic
in a Banach space of real analytic maps.

\pagegoal=50pc

In the present paper, we give the last step in the proof of the above
renormalization conjecture in the (much larger) space of $C^r$ smooth  
unimodal
maps with $r$ sufficiently large. The very formulation of the  
conjecture in
this setting requires some care, because the renormalization operator  
is {\it
not} differentiable in $C^r$. For the correct formulation, see Theorem
\ref{main} below.
To prove the conjecture, we combine Theorem \ref{lyubhyp}
  with some nonlinear functional analysis inspired by the
work of A.~Davie \cite{davie}. In that work, Davie constructs the
stable manifold of the fixed point of the period doubling operator in  
the
space of $C^{2+\epsil}$ maps ``by hand'', showing it to be a $C^1$
codimension-one submanifold of the ambient space, even though the  
operator is
not differentiable. To do this, he first extends the hyperbolic  
splitting of
the derivative at the fixed point from Lanford's \pagebreak Banach space of  
real-analytic
maps to the larger space of $C^{2+\epsil}$ maps (to which the  
derivative
extends as a bounded linear operator). This gives him an extended
codimension-one stable subspace in $C^{2+\epsil}$ to work with, and  
he views
the local stable set in $C^{2+\epsil}$ as the graph of a function  
over the
extended stable subspace. In attempting to prove that such function is  
$C^1$,
he goes around the inherent loss of differentiability of renormalization  
by
first noting that the local unstable manifold coming from Lanford's  
theorem is
still there (and is still smooth in $C^{2+\epsil}$) and then showing  
that
there is afterall a contraction in $C^{2+\epsil}$ towards that  
unstable
manifold, whose elements are {\it analytic} maps. Thus, the loss of
differentiability is somehow compensated by the contraction towards the
unstable manifold. Davie's crucial estimates show that the
renormalization operator in $C^{2+\epsil}$ is sufficiently
well-approximated by the extension of its derivative
in Lanford's space to a bounded linear operator in $C^{2+\epsil}$.
\vskip-3pt

Our approach is based on the idea that whatever Davie can do with\break
Lanford's Banach space relative to the fixed point, we can do with the\break  
Banach
space obtained in  Theorem \ref{lyubhyp}
relative to the whole limit set. There is one fundamental difference,
however. The linear and nonlinear estimates carried out by Davie rely
on the special fact that the period-doubling fixed point is {\it
   concave}. This allows him to prove his main theorems in
$C^{2+\epsil}$ for all $\epsil>0$. By contrast, we cannot -- and
do not -- rely on any such convexity assumptions. We derive our
estimates (in \S \ref{hfgrtyresd} and \S \ref{sec:robust}) directly  
from the
geometric properties of the postcritical set of maps in the limit set
(these properties -- proved in \S \ref{sec:bdedgeom} -- are a
consequence of the real {\it a~priori} bounds). As a result, our local
stable manifold theorem in $C^r$ requires
  $r \ge 2+\alpha$ with $\alpha$ close to one.

We go beyond the conjecture in at least three respects.
First, we show that the local stable manifolds form a $C^0$
lamination whose holonomy is $C^{1+\beta}$ for some $\beta >0$.
In particular, every smooth curve which is transversal
to such a lamination intersects it at a set of constant
Hausdorff dimension less than one.
Second, we prove that the global stable sets are
$C^1$ (immersed) codimension-one submanifolds in $C^r$
provided $r \ge 3+\alpha$ with $\alpha$ close to one
(we globalize the local stable manifolds
via the implicit function theorem,
hence the further loss of one degree of differentiability).
Third, we prove that in an open and dense set of $C^k$
one-parameter families of $C^r$ unimodal maps  (for any $k \ge 2$),
each family intersects the global  stable lamination transversally
at a Cantor set of parameters and the small-scale geometry of this
intersection is the same for all nearby families. In particular, its
Hausdorff dimension is strictly smaller than one.
\pagegoal=50pc

\vskip-2pt
In the path towards these results, we have made an attempt to abstract  
out
the more general features of the renormalization operator in the form  
of a few
properties or ``axioms'' -- the notion of a robust operator introduced  
in Section~\ref{sec:stable}.
We prove a general local stable manifold theorem for robust operators  
there. It is our hope that this might be useful in other
renormalization problems, for example in the case of critical circle  
maps (see
\cite{dFMone} and \cite{dFMtwo}).

\demo{Acknowledgement} We wish to thank M. Lyubich and A. Avila for
several useful discussions and A. Douady for his elegant proof of
Lemma \ref{fgtgtrhhh} (\S \ref{jhgyttre}). We are greatful to the
referee for his keen remarks and for pointing out several
corrections. We also thank FCUP, IMPA, IME-USP, KTH, SUNY Stony Brook
for their hospitality and support during the preparation of this paper.

\section{Preliminaries and statements of results}

In this section, we introduce the basic notions of the theory of
renormalization of unimodal maps. Then we state Sullivan's theorem on
the existence of topological limit sets for the renormalization
operator, the exponential convergence results of McMullen, and
Lyubich's theorem showing the full hyperbolicity of such limit sets in
the space of germs of quadratic-like maps. Finally, we state our main
results extending Lyubich's hyperbolicity theorem to the space of
$C^r$ unimodal maps with $r$ sufficiently large.

\vglue-4pt
\Subsec{Quadratic unimodal maps}\label{quadunimodal}
We describe here two types of ambient spaces of $C^r$ unimodal maps.
These will be determined by two families of Banach spaces, denoted  
$\mathbb{A}^r$
and $\mathbb{B}^r$.

\vglue-4pt\Subsubsec{The Banach spaces $\mathbb{A}^r$}
 Let $I=[-1,1]$ and for all $r \ge 0$ let $C^r(I)$ be the Banach space  
of $C^r$
real-valued functions on $I$.
Here $r$ can be either a nonnegative real number, say $r=k+\alpha$ with
$k \in \mathbb{N}$ and $0 \le \alpha < 1$,
in which case $C^r(I)$ is the space of $C^k$ functions whose
$k^{\rm th}$ derivative is $\alpha$-H\"older,
or else $r=k+\text{Lip}$, in which case $C^r(I)$ means
the space of $C^k$ functions whose $k^{\rm th}$ derivative is
Lipschitz (so whenever we say that $r$ is not an integer,
we include the Lipschitz cases).
Let us denote by $\mathbb{A}^r$ the space $C^r_e(I)$ consisting of all  
$C^r$
functions on $I$ which are {\it even\/} and vanish at the origin, in  
other words
$$
\mathbb{A}^r =
\left\{ v \in C^r(I): v \text{\ is even and } v(0)=0 \right\} \ .
$$
Then $\mathbb{A}^r$ is a closed linear subspace of $C^r(I)$ and
therefore also a Banach space   under the $C^r$ norm. Now, for each
$r\geq 2$, define
$$
\mathbb{U}^r \subset 1+ \mathbb{A}^r  \subset C^r(I)
$$
to be the set of all maps $f:I \to I$  of the form $f(x)=1+v(x)$,
where $v\in\mathbb{A}^r$ satisfies $v''(0)<0$, which are unimodal.
Then $\mathbb{U}^r$ is a Banach manifold; indeed it is an open
subset of the affine space $ 1+ \mathbb{A}^r$. Note that for all $f
\in \mathbb{U}^r$ the tangent space $T_f\mathbb{U}^r$ is naturally
identified with $\mathbb{A}^r$. The elements of $\mathbb{U}^r$ are
called {\it{$C^r$ unimodal maps with a quadratic critical point}}.

\vglue-4pt\Subsubsec{The Banach spaces $\mathbb{B}^r$} We define  
$\mathbb{B}^r$ to
be the space of functions $v:I\to \mathbb{R}$ of the form
$v=\varphi\circ q$ where $q(x)=x^2$ and $\varphi\in C^r([0,1])$
vanishes at the origin. The norm of $v$ in this space is given by
the $C^r$ norm of $\varphi$. This makes $\mathbb{B}^r$ into a Banach
space. Note that for each $s\leq r$ the inclusion map
$j:\mathbb{B}^r\to \mathbb{A}^s$ is linear and continuous (hence
$C^1$). Now, for each $r\geq 1$, let
$$
\mathbb{V}^r \subset 1+ \mathbb{B}^r
$$
be the open subset of the affine space $1+ \mathbb{B}^r$ consisting
of those $f=\phi\circ q$ such that $\phi([0,1]) \subseteq (-1,1]$,
$\phi(0)=1$ and $\phi'(x)<0$ for all $0\leq x\leq 1$. Just as
before, $\mathbb{V}^r$ is a Banach manifold. Note that each $f\in
\mathbb{V}^r$ is a unimodal map belonging to $\mathbb{U}^r$ when
$r\geq 2$. Moreover, the inclusion of $\mathbb{V}^r$ in
$\mathbb{U}^r$ is strict (for each $r\geq 2$).

\Subsec{The renormalization operator}
\label{sub:renorm}
A map $f\in \mathbb{U}^r$ is said to be {\it renormalizable} if there  
exist
$p=p(f)>1$ and $\lambda=\lambda(f)=f^p(0)$
such that $f^p|[-|\lambda|,|\lambda|]$ is unimodal and maps $[-|  
\lambda|,|\lambda|]$
into itself. In this case, with the smallest possible value of $p$,
the map $Rf:[-1,1]\to [-1,1]$ given by
\begin{equation}
Rf(x)\;=\;\frac{1}{\lambda}\,f^p(\lambda x)
\end{equation}
is called the {\it first renormalization} of $f$. We have $Rf\in
\mathbb{U}^r$. The intervals $f^j([-|\lambda|,|\lambda|])$, for $0\leq  
j\leq
p-1$, are pairwise disjoint and their relative order inside $[-1,1]$
determines a {\it unimodal} permutation $\theta$ of  
$\{0,1,\dots,p-1\}$. The
set of all unimodal permutations is denoted $\mathbf{P}$. The
set of $f\in\mathbb{U}^r$ that are renormalizable with the same unimodal
permutation $\theta\in\mathbf{P}$ is a connected subset of  
$\mathbb{U}^r$
denoted $\mathbb{U}^r_\theta$. Hence we have an operator
\begin{equation}
R:\bigcup_{\theta\in\mathbf{P}}\mathbb{U}^r_\theta\,
\to\,\mathbb{U}^r
\ ,
\end{equation}
the so-called {\it renormalization operator}.

Now let us fix a finite subset $\Theta\subseteq \mathbf{P}$. Given an  
infinite
sequence of unimodal permutations
$\theta_0,\theta_1,\dots,\theta_n,\dots\in\Theta$, write
$$
\mathbb{U}^r_{\theta_0,\theta_1,\cdots,\theta_n,\cdots}\;=\;
\mathbb{U}^r_{\theta_0}\cap R^{-1}\mathbb{U}^r_{\theta_1}
\cap\cdots\cap  R^{-n}\mathbb{U}^r_{\theta_n}\cap\cdots
\ ,
$$
and define
$$
\mathcal{D}^r_\Theta\;=\; 
\bigcup_{(\theta_0,\theta_1,\cdots,\theta_n,\cdots)\in
\Theta^{\mathbb{N}}}\mathbb{U}^r_{\theta_0,\theta_1,\cdots,\theta_n,\cdots}
\ .
$$
The maps in $\mathcal{D}^r_\Theta$ are {\it infinitely renormalizable}  
maps
with (bounded) combinatorics belonging to $\Theta$.
Note that $R(\mathcal{D}^r_\Theta)\subseteq
\mathcal{D}^r_{\Theta}$; in fact,
\begin{equation}
R(\mathbb{U}^r_{\theta_0,\theta_1,\cdots,\theta_n,\cdots})
\subseteq \mathbb{U}^r_{\theta_1,\theta_2,\cdots,\theta_{n+1},\cdots}
\ .
\end{equation}

We note that if $f$ is a renormalizable map in $\mathbb{V}^r$, then  
$R(f)$
belongs to $\mathbb{V}^r$ also. Hence, taking  
$\mathbb{V}^r_{\theta}=\mathbb{U}^r_{\theta}
\cap\mathbb{V}^r$, the restriction of the renormalization operator
\begin{equation}
R:\bigcup_{\theta\in\mathbf{P}}\mathbb{V}^r_\theta\,
\to\,\mathbb{V}^r
\end{equation}
is well-defined.\pagebreak

\Subsec{The limit sets of renormalization}
In~\cite{sul}, Sullivan established the existence of
horseshoe-like invariant sets for the renormalization operator,
showing that they all consist of real analytic maps of a special kind,
namely, restrictions to $[-1,1]$ of {\it quadratic-like maps}
in the sense of Douady-Hubbard.
We remind the reader that a quadratic-like map
$f:V \to W$ is a  holomorphic map with the
property that $V$ and $W$ are topological disks with
$V$ compactly contained in $W$,
and $f$ is a proper, degree two branched covering map with a continuous
extension to the boundary of $V$.
The {\it conformal modulus of $f$} is
the modulus of the annulus $W \setminus \overline{V}$.

We are interested only in quadratic-like maps that
commute with complex conjugation, for which $V$ is symmetric about the
real axis. Consider the real Banach space $\mathcal{H}_0(V)$ of
holomorphic functions which commute with complex conjugation and are
continuous up to the boundary of $V$, with the $C^0$ norm. Let
$\mathbb{A}_V\subset \mathcal{H}_0(V)$ be the closed linear subspace
of functions of the form $\varphi=\phi\circ q$, where $q(z)=z^2$ and
$\phi: q(V)\to \mathbb{C}$ is holomorphic with $\phi(0)=0$. Also, let
$\mathbb{U}_V$ be the set of functions of the form $f=1+\varphi$, where
$\varphi=\phi\circ q\in \mathbb{A}_V$ and $\phi$ is {\it univalent} on
some neighborhood of $[-1,1]$ contained in $V$, such that the  
restriction
of $f$ to $[-1,1]$ is unimodal. Then $\mathbb{U}_V$ is an open subset
of the affine space $1+\mathbb{A}_V$, which is linearly isomorphic to
$\mathbb{A}_V$ via the translation by~$1$, and we shall regard
$\mathbb{U}_V$ as an open subset of $\mathbb{A}_V$ itself via this
identification.
For each $a>0$, let us denote by $\Omega_a$ the set of points in the  
complex
plane whose distance from the interval $[-1,1]$ is smaller than $a$.
We may now state Sullivan's theorem as follows.

\begin{theorem} \label{sullii} Let $\Theta\subseteq\mathbf{P}$ be a
nonempty finite set. Then there exist $a>0${\rm ,} a  compact subset
$\mathbb{K}=\mathbb {K}_\Theta \subseteq \mathbb
{A}_{\Omega_a}\cap \mathcal{D}^{\omega}_\Theta$ and $\mu>0$ with
the following properties.
\begin{enumerate}
\item[{\rm (i)}] Each $f\in \mathbb{K}$ has a quadratic-like extension with  
conformal
modulus boun\-ded from below by $\mu$.
\item[{\rm (ii)}]   $R(\mathbb{K})\subseteq\mathbb{K}${\rm ,} and the
restriction of $R$ to $\mathbb{K}$ is a homeomorphism which is
topologically conjugate to the two-sided shift
$\sigma:\Theta^{\mathbb{Z}}\to\Theta^{\mathbb{Z}}$\/{\rm :}\/ in other words{\rm ,}  
there
exists a homeomorphism $H:\mathbb{K}\to\Theta^{\mathbb{Z}}$ such that  
the
diagram
$$
\begin{CD}
\mathbb{K}@>{R}>>\mathbb{K}\\
@V{H}VV             @VV{H}V\\
{\Theta^{\mathbb{Z}}}@>>{\sigma}>{\Theta^{\mathbb{Z}}}
\end{CD}
$$
\vglue-12pt
commutes.
\item[{\rm (iii)}] For all $g\in \mathcal{D}^r_\Theta\cap \mathbb{V}^r${\rm ,} with  
$r \ge 2${\rm ,}
there exists $f\in\mathbb{K}$ with the property that
$||R^n(g)-R^n(f)||_{C^0(I)} \to 0$ as $n\to\infty$.
\end{enumerate}

{\rm For a detailed exposition of this theorem, see Chapter VI of \cite{MS}.} 
\end{theorem}

Later, in \cite{mcmtwo}, C.~McMullen established the exponential  
convergence of
renormalization for bounded combinatorics (using rigidity of towers).  
His
theorem forms the basis for the contracting part of Lyubich's  
hyperbolicity
theorem in \cite{lyubich}.

\begin{theorem}
If $f$ and $g$ are infinitely renormalizable
quadratic-like maps with the same bounded combinatorial type in
$\Theta \subset P${\rm ,} and
with conformal moduli greater than or equal to $\mu${\rm ,}
then
$$
\|R^n f -  R^n g\|_{C^0(I)} \le C \lambda^n
$$
for all $n \ge 0$ where $C=C(\mu, \Theta) > 0$ and
$0 < \lambda= \lambda(\Theta) < 1$.
\end{theorem}

The above result was extended by Lyubich to all combinatorics. In
particular it follows, in the case of bounded combinatorics, that
the exponent $\lambda$ and the constant $C$ in Theorem 2 do not
depend on $\Theta$.
The conclusion of the above theorem can also be improved in  bounded
combinatorics: for $r\geq 3$; the exponential convergence holds in the  
$C^r$
topology if the maps are in  ${\mathbb V}^r$ (see \cite{MMM} and  
\cite{MP}).

In \cite{lyubich}, Lyubich considered the space of quadratic-like
germs modulo affine conjugacies in which the limit set $\mathbb {K}$
is naturally embedded. This space is a manifold modeled on a complex
topological vector space (arising as a direct limit of Banach spaces
of holomorphic maps). In this setting, Lyubich established in [8]
the full hyperbolicity of the renormalization operator. With the
help of Sullivan's real and complex bounds and  Lyubich's theorem we
prove the hyperbolicity of some iterate  of the renormalization
operator acting on a space $\mathbb {A}_{\Omega_a}$ for some $a>0$
(see Theorem \ref{lyubhyp} in \S \ref{jjhtrgfgf}). Then we extend
Davie's analysis for the Feigenbaum fixed  point to the context of
bounded combinatorics to conclude that the hyperbolic  picture also
holds true in the much larger space $\mathbb{U}^r$ (see Theorem
\ref{main} in \S \ref{jjhtrgfgf}).

\vglue-2pt

\Subsec{Hyperbolic basic sets}
\label{ehkytd}
We need to introduce the  well-known concept of hyperbolic basic set  
for nonlinear
operators acting on Banach spaces. Let us consider a Banach space
$\mathcal{A}$, and an open subset $\mathcal{O}\subseteq \mathcal{A}$.

\begin{definition}
\label{hyp}
Let $T: \mathcal{O} \to \mathcal {A}$ be a smooth
nonlinear operator.
A {\it hyperbolic basic set} of $T$ is a compact subset
$\mathbb{K} \subset \mathcal{O}$ with the following properties. \bup
\begin{enumerate}
\item
$\mathbb{K}$ is $T$-invariant and $T|\mathbb{K}$ is a topologically
transitive homeomorphism whose periodic points are dense.
\bup
\item
If $y\in\mathcal{O}$ and all $T$-iterates of $y$ are defined, then
$T^n(y)$ converges to~$\mathbb{K}$. \bup
\item
  There exist a continuous, $DT$-invariant splitting
$\mathcal{A}=E^s_x\bigoplus E^u_x$, for\break $x\in \mathbb{K}$, and
uniform constants $C>0$ and
$0<\theta<1$  such that
$$
\|DT^n(x)\,v\| \leq C \theta^n \|v\|
$$
for all $v\in E^s_x$, as well as
\end{enumerate}
\vglue-22pt
\begin{equation}
\label{rrrrrr}
\|DT^n(x)\,v\| \geq C \theta^{-n} \|v\|
\end{equation} \vglue-22pt\phantom{up}
\begin{enumerate}\item[]
for all $v\in E^u_x$.
\bup
\item[(iv)]
  The dimension of $E^u_x$ is finite and constant.
\end{enumerate}
\end{definition}

The following notions are also standard. Let $\mathcal{A}
(x,\epsil)$ be the ball in $\mathcal A$ with center $x$ and radius
$\epsil$. The  {\it local stable manifold $W^s_\varepsilon(x)$  of
$T$ at $x$} consists of all points $y \in  \mathcal{A} (x,
\varepsilon)$ such that,  for all $n >  0$, we have $T^n (y)  \in
\mathcal{A}  (T^{n}(x),\varepsilon)$ and
$$
\left \|T^n (y)- T^n (x) \right \| \to 0 ~{\rm when}~ n \to \infty \ .
$$
The {\it local unstable manifold $W^u_\varepsilon(x)$  of $T$ at
$x$} consists of all points $y \in  \mathcal{A} (x,\varepsilon)$
such that, when $y_0=y$,  for all $n \ge 1$  there exists $y_n
\in \mathcal{A} (T^{-n}(x),\varepsilon)$ such that $y_{n-1}=T(y_n)$
and
$$
\|T^{-n}(x)-y_n\| \to 0 ~{\rm when}~ n \to \infty  \ .
$$
Finally the {\it global stable set of $T$ at  $x$}  is defined as
$$
W^s(x) = \left\{y\in \mathcal{O}\ :\ \|T^n(y)-T^n(x)\|\to 0    ~{\rm  
when}~ n \to \infty     \right\}
\ .
$$
The question arises as to whether these sets have smooth manifold  
structures.
We have the following general result.

\begin{theorem}
\label {hpicture}
If $\mathbb{K}$ is a hyperbolic basic set of a $C^1$ operator
$T:\mathcal{O}\to \mathcal{A}$ then
\begin{enumerate}
\item  the local stable \/{\rm (}\/resp.\ unstable\/{\rm )}\/ set at
$x\in\mathbb{K}$ is a $C^1$ Banach submanifold of $\mathcal{A}$ which is
tangent to $E^s_x$ \/{\rm (}\/resp.\ $E^u_x${\rm )} at $x$.
\item   If $y\in W^s(x)$ then
$$
\left \|T^n(x)-T^n(y) \right \| \leq C  \theta^n \|x-y\|
\ .
$$
Moreover{\rm ,} $T(W^u_\epsil(x))\supseteq W^u_\epsil(T(x))${\rm ,} the  
restriction of
$T$ to $W^u_\epsil(x)$ is one-to-one and for all $y\in  
W^u_\epsil(x)${\rm ,}
$$
\left \|T^{-n}(x)-T^{-n}(y)\right \|\leq C   \theta^n \|x-y\|
\ .
$$
\item  If $y\in \mathcal {A}(x,\epsil)$
is such that $T^i(y) \in \mathcal {A}(T^i(x),\epsil)$ for $i\leq n$  
then
$$ {\rm dist} \left (T^n(y), W^u_\epsil (T^n(x)) \right )
\leq C \theta^n,
\  \text{as well as}\  \
{\rm dist} \left (y, W^s_\epsil(x)  \right ) \leq C    \theta^n \ .
$$
\item  The family of local stable manifolds \/{\rm (}\/and also the family of
local unstable manifolds\/{\rm )}\/ form a $C^0$ lamination\/{\rm :}\/ the tangent spaces  
to the
leaves vary continuously.
\end{enumerate}
\end{theorem}

We do not prove this theorem here since we will not use it, but
instead make the following comments.
Using the arguments of Hirsch-Pugh in \cite{HP}, we can prove that the  
local unstable
set is a smooth manifold. The local stable set is also a smooth  
manifold, but
a different proof is needed: one can use the ideas of Irwin in \cite  
{irwin}.
See also Theorem 2.1 on page 375 of \cite{Rob}.
In both cases the smoothness can be improved to $C^k$ if the operator  
$T$ is $C^k$.

For invertible operators the global stable set is also a smooth
submanifold.
In the noninvertible situation, this is not always true.
However, we will prove in \S \ref{sec:global} that this is the case for  
the
renormalization operator acting on ${\mathbb V}^r$, provided $r \ge 3  
+\alpha$
and $\alpha>0$ is close to one.

\Subsec{Hyperbolicity of renormalization}
\label{jjhtrgfgf}
In the present paper we prove three main theorems. The first main
theorem shows that there exists a real Banach space of analytic maps,
containing the topological limit set $\mathbb{K}$ of renormalization,
on which the renormalization operator $R$ acts as a real-analytic  
operator
and has $\mathbb{K}$ as a hyperbolic basic set. More precisely, we
have the following result.

\begin{theorem}[Hyperbolicity in a real Banach space]
\label{lyubhyp} \hskip-5pt
There exist $a\!>\!0${\rm ,} an open set
$\mathbb{O} \subset \mathbb {A}=\mathbb {A}_{\Omega_a}$
containing $\mathbb{K}=\mathbb{K}_\Theta$
and a positive integer $N$ with the following property.
There exists a real analytic operator
$T: \mathbb{O}  \to \mathbb{A}$
  having  $\mathbb{K}$
as a hyperbolic basic set with codimension-one stable manifolds at
each point{\rm ,} such that $T(f)|[-1,1]= R^N(f|[-1,1])${\rm ,} for all $f \in
\mathbb{O}${\rm ,} is the $N^{\rm th}$ iterate of the renormalization operator.
\end{theorem}

The proof of this theorem, presented in \S 3 (see Theorem \ref{B15}),  
combines
Lyubich's hyperbolicity results with Sullivan's real and complex bounds.

The second main theorem establishes the ``hyperbolicity'' of
renormalization in $\mathbb{U}^r$.
As we have mentioned before, the renormalization operator is not smooth
in $\mathbb{U}^r$, so the definition of hyperbolicity of an invariant  
set
does not even make sense. However, the hyperbolic picture holds in this
situation. More precisely, we have the following theorem.

\begin{theorem}[Hyperbolic Picture in $\mathbb{U}^r$]
\label{main}
If $r\geq 2+\alpha${\rm ,} where $\alpha >0$ is close to one{\rm ,} then
statements {\rm (i), (ii), (iii)} and {\rm (iv)} of Theorem {\rm \ref{hpicture}} hold  
true for
the renormalization operator acting on $\mathbb{U}^r$.
Furthermore{\rm ,}
\begin{enumerate}
\item
  the local
unstable manifolds are real analytic curves\/{\rm ;}\/
\item the local stable manifolds are
of class $C^1${\rm ,} and together they form a continuous lamination whose  
holonomy
is $C^{1+\beta}$ for some $\beta>0${\rm .}
\end{enumerate}
\end{theorem}

The main difficulty behind the proof of this theorem is the fact that  the
operator $T$ is not Fr\'echet differentiable in $C^r$ (in fact it is  
only
continuous in a dense subset of $\mathbb{U}^r$).
However, as we shall see in \S \ref{sec:axiom5}, it is a $C^1$ mapping
from its domain
in $\mathbb{U}^r$ into $\mathbb{U}^s$ if $s<r-1$ (even for $s=r-1$ if  
$r$ is
an integer). Hence its tangent map defines a continuous map
$L \colon \mathbb K\times \mathbb{A}^r \to \mathbb K \times  
\mathbb{A}^s$ by
$L(g, v)= (T(g), DT(g)(v) )= (T(g), L_g(v))$. The bounded
linear mappings $L_g \colon \mathbb{A}^r \to \mathbb{A}^s$ extend to  
bounded
linear operators $L_g \colon \mathbb{A}^t \to  \mathbb{A}^t$ for all
$0 \le t \le r$.
Although $L_g$ is not the derivative of $T$ at $g$ in $C^r$, it is  
nevertheless a
sufficiently good linear approximation to $T$ near $g$
(see the properties of Definition \ref{def4}, checked in \S  
\ref{sec:robust}).

\begin{corollary}[Hyperbolic Picture in $\mathbb{V}^r$]
\label{mainmain}
If $r\geq 2+\alpha${\rm ,} where $\alpha >0$ is close to one{\rm ,} then
statements {\rm (i), (ii), (iii)} and {\rm (iv)} of Theorem {\rm \ref{hpicture}} hold  
true for
the renormalization operator acting on $\mathbb{V}^r$.
Furthermore{\rm ,}
\begin{enumerate}
\item
  the local
unstable manifolds are real analytic curves\/{\rm ;}\/
\item the local stable manifolds are
of class $C^1${\rm ,} and together they form a continuous lamination whose
holonomy is $C^{1+\beta}$ for some $\beta>0$.
\end{enumerate}
\end{corollary}

For the proofs of Theorem \ref{main} and Corollary \ref{mainmain}, see
Section~\ref{sec:robust}.

By an argument using the implicit function theorem and the results in
\cite{MMM}, which {\it a priori\/} are valid only in
$\mathbb{V}^r$, we shall prove in Section~\ref{sec:mainglobal} our third  
main theorem,
which we state as follows.

\begin{theorem}
\label{corol} If  $r \ge 3  +  \alpha${\rm ,} where $\alpha>0$ is close to
one{\rm ,} then the following assertions hold true for the renormalization
operator acting in $\mathbb{V}^r$\/{\rm :}\/
\begin{enumerate}
\item The global stable sets are $C^1$ immersed
submanifolds.
\item  For each integer $2 \le k \le r${\rm ,} there exists an open dense set  
of
   $C^k$ one-parameter families of maps in $\mathbb{V}^r$
all of whose elements intersect the global stable lamination of
$(T,K_{\Theta})$ transversally.
\item In each such family{\rm ,} the set of parameters where the  
intersections occur
is a Cantor set which is locally $C^{1+\beta}$ diffeomorphic
to the corresponding Cantor set of the quadratic family. In
particular{\rm ,} its Hausdorff dimension is a universal number depending  
only on
$\Theta$ which lies strictly between zero and one if $\Theta$ has more  
than
one element.
\end{enumerate}
\end{theorem}

It is worth emphasizing that when a generic family (in the
sense of the above corollary) intersects the stable lamination at a  
point,
then any neighborhood of this point in parameter space contains a
renormalization window that is mapped under a suitable power of the
renormalization operator onto a {\it full} transversal family.

\section{Hyperbolicity in a Banach space of real analytic maps}
\label{grtsfd}

In this section we give a proof of   Theorem \ref{lyubhyp}.
Using the real and complex bounds given by
Sullivan in \cite{sul},
we  prove in \S \ref{Bfds3}  that there is an iterate of the  
renormalization operator
which extends as a real analytic map
  $T$ to an open set ${\mathbb O}_{\Omega_a}$
  of the Banach space $\mathbb {A}_{\Omega_a}$ consisting
  of real analytic maps whose domain is
an $a$-neighborhood of  the interval $[-1,1]$, for a suitable $a>0$.
The maps $g \in \mathbb K$  have unique extensions belonging to
   ${\mathbb O}_{\Omega_a}$. In \S \ref{hybriddd},
using Lemmas 4.16 and 4.17 in  Lyubich's paper \cite{lyubich},
we show that the  hybrid conjugacy classes
of  the maps  $g \in \mathbb K$
form a continuous lamination of codimension one real analytic
manifolds.
Then in \S \ref{skewww}  we construct a skew-product renormalization  
operator
that satisfies properties (W1) to (W4) on page 395 of
\cite{lyubich} in the real-analytic case (restated in \S  
\ref{ghjtyrert}).
By Theorems 8.2 and
  8.8 in \cite{lyubich}  the
  skew-product renormalization operator
will have  fiberwise stable and unstable  leaves (as defined in
\S \ref{fefefsds}).   The local stable leaf at $g \in \mathbb K$
is a relatively open set  of  the hybrid conjugacy class of $g$.
Then using the  skew-product renormalization operator,
  we prove in \S \ref{ghjtyrert} that $\mathbb K$
is a basic set for the real analytic renormalization operator
$T:{\mathbb O}_{\Omega_a} \to \mathbb {A}_{\Omega_a}$.

\Subsec{Real analyticity of the renormalization operator}
\label{Bfds3}
Using  Sullivan's real and complex bounds
  in \cite{sul}, we will show that there exists $a>0$ such that some  
iterate
$T:{\mathbb O}_{\Omega_a} \to \mathbb {A}_{\Omega_a}$ of the  
renormalization operator
is a (well-defined) real analytic operator with a compact derivative.

For each $f\in\mathbb{K}$, let $\mathcal{I}_f\subseteq [-1,1]$ be the
postcritical set of $f$ (the Cantor
attractor of $f$). For each $k\ge 0$, we can write
$$
R^k f (x) = \Lambda_k^{-1} \circ  f^{p_k} \circ \Lambda_k (x)
$$
where
\begin{eqnarray}
\label{scalings}
p_k&=&p(f,k)=\prod_{i=0}^{k-1} p(R^if)\ ,\\
\lambda_k&=&\lambda (f,k)=\prod_{i=0}^{k-1} \lambda{(R^if)}\  
,\nonumber\\
\Lambda_k(x)&=&\Lambda (f,k)(x)= \lambda_k \cdot x \ , \nonumber
\end{eqnarray}
with $p(\cdot)$ and $\lambda(\cdot)$ as defined in \S \ref{sub:renorm}.
Consider the renormalization intervals
$\Delta_{0,k}=\Delta_{0,k}(f)=[-|\lambda_k|,|\lambda_k|]
\subset [-1,1]$, and define
$\Delta_{i,k}=\Delta_{i,k}(f)=f^i(\Delta_{0,k})$ for
$i=0,1, \dots, p_k-1$.
The collection
${\mathbf C}_k=\{\Delta_{0,k}, \dots, \Delta_{p_k-1,k}\}$
consists of pairwise disjoint intervals at level $k$.
Moreover,
$\bigcup \{\Delta:\Delta\in {\mathbf C}_{k+1}\} \subseteq \bigcup\break
\{\Delta:\Delta\in {\mathbf C}_k\}$ for all $k \ge 0$ and we have
$$
{\mathcal I}_f=
\bigcap _{k=0}^{\infty}
\bigcup_{i=0}^{p_k-1}
\Delta_{i,k} \ .
$$
\begin{definition}
The set ${\mathcal I}_f$ has {\it geometry bounded by $0<\tau<1$} with  
respect
to $({\mathbf C}_k)_{k \in \mathbb N}$
if  the following conditions are met for $k\geq 1$.
\begin{enumerate}
\item If $\Delta_{j,k+1} \subset \Delta_{i,k}$
then $\tau <\left| \Delta_{j,k+1} \right|  /  \left| \Delta_{i,k}  
\right| <
1-\tau$.
\item If $I$ is a {\it connected component} of
$\Delta_{i,k} \setminus \bigcup_j \Delta_{j,k+1}$
then  $\tau < \left|I \right|  /  \left| \Delta_{i,k} \right| < 1-\tau$.
\end{enumerate}
\end{definition}

By  Sullivan's real bounds
(see \cite{sul} and Section VI.2 on page 453 of \cite{MS}),
there exists $\alpha > 0$, such that  for  every  $g \in  \mathbb{K}$
the set ${\mathcal I}_g$ has   geometry bounded by $\alpha$ with respect
to $({\mathbf C}_k)_{k \in \mathbb N}$.

The following
result is a consequence of Sullivan's complex bounds
(see  \cite{sul} and Section VI.5 on page 483 of \cite{MS}).

\begin{theorem}
\label{aaa3}
There exist $\mu > 0${\rm ,}  $N_0 > 0$ and a neighborhood $V$ of the dynamics
with the following properties. Every  $g \in  \mathbb{K}$
   extends to a holomorphic  map  $g:V \to \mathbb C$  and
  for every $N \ge N_0$
there exists a symmetric neighborhood $O_{g,N}$  of the interval
$\Delta_{0,N}(g)$ such that
\begin{enumerate}
\item  the diameter of the set $g^i(O_{g,N})\subset V$
is comparable to the length $|\Delta_{i,N}(g)|$
of the interval $\Delta_{i,N}(g)$  for every $0 \le i \le p=p(N,g)${\rm ;}
\item  the map $g^p:O_{g,N} \to  g^p(O_{g,N})$ is
a quadratic-like map with conformal modulus greater than $\mu> 0$.
\end{enumerate}
  \end{theorem}

Applying Theorem  \ref{aaa3} (ii) to $g \in \mathbb{K}$,
we see that $R^N(g)$ has a quadratic-like extension to
\begin{equation}
\label{crtyh}
U_{g,N} = \Lambda_{g}^{-1} (O_{g,N})
\end{equation}
(where $\Lambda_g=\Lambda(g,N)$) and such an extension has conformal  
modulus
greater than $\mu>0$.

Recall that the {\it  filled-in  Julia set ${\mathcal K}_f$
of a quadratic-like map $f:U \to U'$}
is the set $\{z:f^nz \in U, n=0,1,\dots\}$, and
its boundary is the Julia set ${\mathcal J}_f$ of $f$.
Since all maps in $\mathbb K$ have conformal modulus greater than
or equal to $\mu>0$, we deduce from
  Proposition 4.8 on page 83 of McMullen's book \cite{mcmtwo} that
there exists  $b>0$
such that for every $g \in \mathbb{K}$ we have
\begin{equation}
\label{frffr}
\Omega_{b}({\mathcal K}_{R^N(g)})
  \subset  U_{g,N}
\ .
\end{equation}
Here the notation $\Omega_{\varepsilon}(K)$ means the set of all points  
whose
distance from $K$ is less than $\varepsilon/2$ times the diameter of  
$K$.

For each neighborhood $U$ of $[-1,1]$ in $\mathbb{C}$, symmetric about  
the
real axis, we consider the real Banach space $\mathbb{A}_U$ of  
holomorphic
functions defined earlier. We denote by $\mathbb{A}_U(g,\delta)$ the  
open ball
of radius $\delta$ around $g$.
By \eqref{frffr}, the inclusion map $i_{g,N}:\mathbb{A}_{U_{g,N}}
\to \mathbb{A}_{\Omega_\alpha}$
is a well-defined compact linear operator for every $0<\alpha<b$.

\begin{lemma}
\label{B3}
Let $\mu>0$ and $N_0>0$  be as in Theorem  {\rm \ref{aaa3}} and $b>0$ as in
{\rm \eqref{frffr}.}
For every $0<\alpha<b$ there exist
$N>N_0$ and  $\delta_0 > 0$
such that
\begin{enumerate}
\item for  every $g \in \mathbb K${\rm ,}
  the operator $T_{g,N}:\mathbb{A}_{\Omega_\alpha}(g, \delta_0)
\to \mathbb{A}_{U_{g,N}}$
is well-defined if we set
$$
T_{g,N} (f) = {\Lambda_f}^{-1} \circ f^p \circ \Lambda_f
: U_{g,N}\to \mathbb{C}
\ ,
$$
where $p=p(f,N)=p(g,N)${\rm ,} $\Lambda_f=\Lambda(f,N)${\rm ,} and $T_{g,N} (f)$
is a quadratic-like map with conformal modulus greater than   $\mu/2${\rm ;}
\item
the  operator
$T:{\mathbb O}_{\Omega_\alpha}  \to \mathbb{A}_{\Omega_\alpha}$
  given by
$T =i_{g,N} \circ T_{g,N}$ is  real-analytic
with a compact derivative{\rm ,}
where
$${\mathbb O}_{\Omega_\alpha}=
\bigcup_{g\in \mathbb{K}}\mathbb{A}_{\Omega_\alpha}(g, \delta_0)
\ .
$$
\end{enumerate}
\end{lemma}

\Proof
By  Sullivan's real bounds,
there exist $C_1>1$ and $0<\nu_1<\nu_2<1$ such that
for all $g \in \mathbb K$, all $k \in \mathbb N$  and
all $0 \le  j  \le p(k,g)-1$, we have
$
C_1^{-1} \nu^k_1 < |\Delta_{j,k}(g)| < C_1 \nu^k_2 \ .
$
Thus, by property (i) in Theorem \ref{aaa3},
for every     $\alpha > 0$
there is $N>0$ so large that
the open sets $g^j( O_{g,N})$ have diameter smaller than
$\alpha/3$ for all $0 \le j \le p(N,g)$.  Recall that $O_{g,N}=
\Lambda_g( U_{g,N})$.
By a continuity argument,
there is $\delta_g > 0$ such that
  for every $f \in  \mathbb{A}_{\Omega_\alpha}(g, \delta_g)$,
the restriction  $f|[-1,1]$ is   $N$-times renormalizable,
$f^j (\Lambda_f( U_{g,N}))  \subset  \Omega_{\alpha/2}$
  for  every $0 \le j \le p=p(N,f)$, and moreover
$f^p:   \Lambda_f( U_{g,N}) \to {f}^p  (  \Lambda_f( U_{g,N}))$
is a quadratic-like map with   conformal modulus greater than $\mu/2$.
By compactness of $\mathbb K$ in $\mathbb{A}_{\Omega_\alpha}$,
there is a  finite set
$\left \{  g_i :i=1,\dots,l \right\}$ such that
$$
\mathbb K \subset \bigcup_{i=1}^l  \mathbb{A}_{\Omega_\alpha}(g_i,
\delta_{g_i}/2)
\ .
$$
Set $\delta_0 =\min_{i=1,\dots,l} \{\delta_{g_i}/2 \}$.
Then, for every $g \in {\mathbb K}$ there exists $i=i(g)$ such that
$ \mathbb{A}_{\Omega_\alpha}(g, \delta_0) \subset
  \mathbb{A}_{\Omega_\alpha}(g_i, \delta_{g_i})$.
Hence $T_{g,N} (f)$ is well-defined, and it is
a quadratic-like map with conformal modulus greater than $\mu/2$,
for every $f \in  \mathbb{A}_{\Omega_\alpha}(g, \delta_0)$ which proves  
(i).

Note that the real Banach space $\mathbb{A}_{\Omega_\alpha}$ is
naturally embedded in the complex Banach space
$\mathbb{A}_{\Omega_\alpha,\mathbb{C}}$ of maps $f:\Omega_\alpha \to
\mathbb{C}$ which are holomorphic and continuous up to the boundary
and that $T_{g,N}$ extends to an operator $T_{g,N}^\mathbb{C}$ in an
open set of $\mathbb{A}_{\Omega_\alpha,\mathbb{C}}$, given by the
same expression. Applying Cauchy's integral formula, we see that
$T_{g,N}^\mathbb{C}$ is complex-analytic, and so $T_{g,N}$ is real-analytic. Since by Montel's theorem
the inclusion $i_{g,N}$ is a compact linear operator, we deduce that $T:{\mathbb
O}_{\Omega_\alpha} \to \mathbb {A}_{\Omega_\alpha}$ is a
real-analytic operator with a compact derivative, which proves (ii).
\hfq

\Subsec{Real analytic hybrid conjugacy classes}
\label{hybriddd}
We will introduce later (in \S \ref{skewww}) a
skew-product renormalization operator.
The fiberwise local stable manifolds of such a skew-product
-- which will be used to determine the stable manifolds of the  
real-analytic operator
$T:\mathbb{O}_{\Omega_a} \to \mathbb{A}_{\Omega_a}$, for some suitable  
$a>0$ --
turn out to be openly contained in the hybrid conjugacy classes
of the maps in the limit set $\mathbb{K}$.
Here we analyze the manifold structure of hybrid classes in more detail.

A homeomorphism $h:U \to V$,
where $U$ and $V$ are contained in
${\mathbb C}$ or $\overline{\mathbb C}$,
  is {\it quasiconformal}
if it has locally square integrable distributional derivatives
$\partial h$, $\overline{\partial} h$, and
there exists $\epsil < 1$ with the property that
$\left|\overline{\partial} h/\partial h \right|\le \epsil$ almost  
everywhere.
The Beltrami differential $\mu_h$ of $h$ is
given by
$\mu_h =\overline{\partial} h/\partial h$.
A quasiconformal   map $h$ is $K$ quasiconformal   if
$K \ge (1+||\mu_h||_\infty)/(1-||\mu_h||_\infty)$.

Two quadratic-like maps $f$ and $g$ are
{\it hybrid conjugate} if there
is a quasiconformal conjugacy $h$ between $f$ and $g$
with the property that  $\overline{\partial} h (z)=0$
for almost every $z \in {\mathcal K}_f$.
Let us denote by ${\mathcal H}(f)$ the hybrid conjugacy class of $f$.

By a slight abuse of notation, we will denote by
${\mathbb K} \cap {\mathbb A}_V (g,\delta)$
the set of maps $f \in \mathbb{A}_{V} (g,\delta)$ with the property that
$f|[-1,1]$ belongs to ${\mathbb K}$.

In the proof of the following theorem, we will need to work with the  
complexification of
$\mathbb{A}_{V}$.
Let $\mathbb{A}_{V,\mathbb C}$ be the complex Banach space of all
holomorphic maps $f :V \to \mathbb C$ with a continuous extension to  
the boundary of $V$.
Let  $\mathbb{A}_{V,\mathbb C}(f,\delta)$ be the open ball in
$\mathbb{A}_{V,\mathbb C}$ centered in $f$ and with radius $\delta>0$.
Let $C:\mathbb{A}_{V,\mathbb C} \to \mathbb{A}_{V,\mathbb C}$
be the conjugation operator given by $C(f)=c \circ f \circ c$,
where $c(z)=\overline{z} \in \mathbb C$.
We note that  $f \in \mathbb{A}_{V}$ if and only if
  $f \in \mathbb{A}_{V,\mathbb C}$ and $C(f)=f$.

\begin{theorem}
\label{dfdfd}
For every  $g \in  \mathbb K${\rm ,}
there exists a symmetric  neighborhood
  $\hat V_g$ of the reals such that $g$ has a
quadratic-like extension to $\hat V_g$ \/{\rm (}\/also denoted by $g${\rm ),}  
$\hat V_g$
contains a definite neighborhood of ${\mathcal K}_g$ and
for every   neighborhood  $V \subset \hat V_g$
  symmetric with respect to $\mathbb R$
and with the property that $g|V$ is a quadratic-like map{\rm ,}
there is $\delta_{g,V} > 0$ such that
for all $f \in \mathbb K \cap \mathbb{A}_{V}(g, \delta_{g,V})${\rm ,}
$${\mathcal H}_V (f)={\mathcal H} (f) \cap \mathbb{A}_{V}(g,  
\delta_{g,V})$$  are
codimension-one{\rm ,} real-analytic leaves varying continuously with $f$.
\end{theorem}

\Proof
By   Lemmas 4.16 and 4.17  on page 354 of
    Lyubich's paper \cite{lyubich}, we obtain
that for all $f \in \mathbb {K} \cap \mathbb{A}_{V,\mathbb {C}}(g,  
\delta_{g,V})$,
$
{\mathcal H}_{V,\mathbb {C}} (f) ={\mathcal H} (f) \cap  
\mathbb{A}_{V,\mathbb {C}}(g, \delta_{g,V})
$ are
codimension-one complex analytic leaves varying continuously with $f$.
If $f \in  \mathbb{A}_{V}(g, \delta_{g,V})$ then the hybrid
conjugacy class of
$f$ in $\mathbb{A}_{V,\mathbb {C}}(g, \delta_{g,V})$
is invariant under the conjugation operator $C$.
Hence, the tangent space $T_f {\mathcal H}_{V,\mathbb {C}} (f)$ at $f$  
to
its hybrid conjugacy class  is invariant under the conjugation operator  
$C$,
and there is a one-dimensional transversal $E_f$
to $T_f {\mathcal H}_{V,\mathbb {C}} (f)$
which   is also
invariant under the conjugation operator $C$.
Locally ${\mathcal H}_{V,\mathbb {C}} (f)$ is a graph
of $\mathcal{G}:Z \subset T_f {\mathcal H}_{V,\mathbb {C}} (f) \to E_f$
with the property that if $h=v+\mathcal{G}(v)$ then  
$C(h)=C(v)+\mathcal{G}(C(v))$.
Thus, locally ${\mathcal H}_{V} (f)$ is also
the graph of $\mathcal{G}|Z \cap  \mathbb{A}_{V}(g, \delta_{g,V})$,
and so it is a codimension-one, real-analytic leaf.
Since the complex-analytic leaves ${\mathcal H}_{V,\mathbb {C}} (f)$
vary continuously with  $f$, we deduce that the real-analytic leaves ${\mathcal H}_{V} (f)$ also
vary continuously with $f$.
\hfq

\Subsec{Hyperbolic skew-products}
\label{fefefsds}
Before going further, we pause for a moment to introduce the elementary  
concept of
hyperbolic skew product in an abstract setting.
Let $\mathbb{K}$ be a compact metric space and
assume that $\mathbb{K}$ is totally disconnected.
Let $\mathcal{F}$ be a finite collection of (real) Banach spaces,
say $\mathcal{F}=\{\mathcal{A}_1,\mathcal{A}_2,\dots,\mathcal{A}_N\}$,
and assume we have a locally constant map $\varphi:\mathbb{K} \to  
\mathcal{F}$.
We write $\mathcal{A}_x=\varphi(x) \in \mathcal{F}$, for all $x \in  
\mathbb{K}$.
Let $E=\cup_{x \in \mathbb{K}} \{x\} \times \mathcal{A}_x$.
We endow $E$ with a topology as follows.
If $\mathbb{K}_i=\varphi^{-1}(\mathcal{A}_i)$,
then  $\mathbb{K}_i$ is an open and closed set in $\mathbb{K}$,
for each $i=1,2,\dots,N$. Note that $E$ is the disjoint union
of $\mathbb{K}_i \times \mathcal{A}_i$, $i=1,2,\dots,N$.
Hence endow each factor $\mathbb{K}_i \times \mathcal{A}_i$ with the
product topology and then $E$ with the union topology.
It is clear that $E$ is metrizable also.
The natural projection $E \to \mathbb{K}$ is open and continuous.
We shall assume that there exists a continuous injection
  $\mathbb{K}_i \to \mathcal{A}_i$ for each $i$,
and will accordingly identify each $x \in \mathbb{K}_i$
with its image in $\mathcal{A}_i$.

Now suppose $T:\mathbb{K} \to \mathbb{K}$ is a homeomorphism
(in the case we are interested, $T$ is transitive),
and also that for each $x \in \mathbb{K}$ we have a
real-analytic map $S_x:\mathcal{A}_x(x,\delta)\to \mathcal{A}_{T(x)}$,
where $\mathcal{A}_x(x,\delta)=\{x+v \in  
\mathcal{A}_x:\|v\|_{\mathcal{A}_x} < \delta\}$.
We define a skew-product operator $S:E(\delta) \to E$ over $T$,
where
$$
E(\delta) = \left\{(x,y):x \in \mathbb{K}, y \in \mathcal{A}_x,
\|y-x\|_{\mathcal{A}_x} < \delta\right\},
$$
by $S(x,y)=(T(x),S_x(y))$.

\begin{definition}
We say that $S$ is {\it fiberwise hyperbolic} if there exists a
continuous splitting $\mathcal{A}_x=E^s_x \bigoplus E^u_x$ with {\rm dim}  
$E^u_x=1$
which is invariant in the sense that $DS_x(E^s_x) \subseteq E^s_{T(x)}$  
and
$DS_x(E^u_x) \subseteq E^u_{T(x)}${\rm ,} satisfying for all $v^s \in E_x^s$  
and all
  $v^u \in E_x^u$ the inequalities
   \begin{eqnarray*}
\left \|D \left( S_{T^{n-1} (x)} \circ \dots S_x\right) (x)   v^s
\right\|_{\mathcal{A}_{T^n (x)}}
& \le & C \theta^n \|v^s \|_{\mathcal{A}_x} \nonumber \\
\left \|D \left( S_{T^{n-1}(x)} \circ \dots S_x\right) (x)   v^u
\right \|_{\mathcal{A}_{T^n (x)}}
& \ge & C^{-1} \theta^{-n} \| v^u \|_{\mathcal{A}_x} \ ,
\end{eqnarray*}
where  $C>1$ and $0<\theta < 1$ are uniform constants on $g$.
\end{definition}

\begin{definition}
The {\it fiberwise  local stable manifold} $W^s_\beta(x)$  of $S$
at $x$ consists of all points $y \in  \mathcal{A}_x (x, \beta)$
such that  for all $n \ge 1$, we have
$S_{T^{n-1}(x)} \circ \dots \circ S_x(y)  \in
\mathcal{A}_{T^n(x)}  (T^{n}(x), \beta)$
and
$$\left \|S_{T^{n-1}(x)} \circ \dots \circ S_x (y)-
      S_{T^{n-1}(x)} \circ \dots \circ S_x (x) \right
\|_{\mathcal{A}_{T^n(x)}} \le C \theta^n
$$
where $C>0$ and $0<\theta <1$ are uniform constants on $x \in \mathbb  
K$.
The {\it fiberwise  local unstable manifold $W^u_\beta(x)$  of $S$
at $x$} consists of all points $y \in  \mathcal{A}_x (x, \beta)$
such that when $y_0=y$,  for each $n \ge 1$  there exists
$y_n \in \mathcal{A}_{T^{-n}(x)}$ such that
$y_{n-1}=S_{T^{-n}(x)}(y_n)$  and
$
\|T^{-n}(x)-y_n\|_{\mathcal{A}_{T^{-n}(x)}} \le C \theta^n.
$
\end{definition}

\vglue-12pt

\Subsec{Skew-product renormalization operator}
\label{skewww}
Our goal in this section is to build a skew-product renormalization
operator that will
play a central role in the proof that $\mathbb K$
is a basic set for $T:{\mathbb O}_{\Omega_a} \to \mathbb  
{A}_{\Omega_a}$, for a suitable $a>0$.
Our  skew-product is constructed so as to satisfy properties (W1) to  
(W4) on page 395 of
\cite{lyubich} in the real analytic case -- restated in \S  
\ref{ghjtyrert} -- and
therefore will have fiberwise stable and unstable manifolds, as we will
explain in that section.

Using Theorem \ref{aaa3} and \eqref{frffr}, we know that for every $0  
<\alpha < b$,  $\mathbb{K}$
injects continuously into $\mathbb{A}_{\Omega_{\alpha}}$.
Hence for $f,g \in \mathbb{K}$
we define  
$\text{dist}_\mathbb{K}(f,g)=\|f-g\|_{\mathbb{A}_{\Omega_{\alpha}}}$.
We also denote by
$\mathbb{K} (g, \varepsilon)$ the ball of radius $\varepsilon$ centered  
at $g$
in this metric.
The metric is compatible with the natural topology of $\mathbb{K}$,
independently of which $\alpha$ we take.

\begin{lemma}
\label{B33}
The filled-in Julia set ${\mathcal K}_{g}$ varies continuously in
  the\break Hausdorff metric with respect to $g \in  \mathbb{K}$.
\end{lemma}

\Proof
We need to show that for every $\varepsilon>0$ there exists $\delta>0$
such that if $\text{dist}_\mathbb{K}(f,g) < \delta$ then
(a)
${\mathcal K}_{g} \subset \Omega_\varepsilon ({\mathcal K}_{f})$
and
(b)
${\mathcal K}_{f} \subset \Omega_\varepsilon ({\mathcal K}_{g})$.
Let $U = U_{R^{-N_0}(g),N_0} \subset \mathbb C$ be the symmetric
neighborhood of $[-1,1]$ given by  Lemma \ref{B3}.
Since the operator $T_{R^{-N_0}(g),N_0}$ is continuous,
every $f \in \mathbb{K}$ sufficiently close to $g$ in
$\mathbb{A}_{\Omega_\varepsilon}$ is quadratic-like on $U$
($f=T_{R^{-N_0}(g),N_0} (T^{-1}(f)):U \to \mathbb{C}$ )
and is also close to $g$ in $\mathbb{A}_U$.

To prove (a), cover $\mathcal{K}_{g}$ by finitely many disks
$D(z_i(g), \varepsilon/2)$, $i=1,2,\dots,m$, where each $z_i(g)$ is an
expanding periodic point of $g$. For $f$ sufficiently close to~$g$,
the corresponding periodic points $z_i(f) \in D(z_i(g), \varepsilon/2)$.
Hence each $z \in {\mathcal K}_{g}$ is at distance at most $\varepsilon$
from some $z_i(f)$, which proves (a).

To prove (b), let $n > 0$ be so large that
  $W=g^{-n}(U) \subset \Omega_\varepsilon ({\mathcal K}_{g})$.
Since   $f$ is   close to $g$   and $W \subseteq U$ is symmetric
   $f:W \to f(W)$ is quadratic-like also,
whence ${\mathcal K}_{f} \subset W \subset
\Omega_\varepsilon ({\mathcal K}_{g})$ and so (b) is proved.
\hfq

\begin{lemma}
\label{B55}
Let $g \in  \mathbb{K}$ and let $V \subset \mathbb C$ be a symmetric  
neighborhood of
$[-1,1]$  which is compactly contained in $\Omega_{ b/2}  ({\mathcal
K}_{g}) ${\rm ,}  where $b$ is given by  \eqref{frffr}. Then for all
$\varepsilon>0$ sufficiently small $\mathbb{K} \cap \mathbb{A}_{V}(g,
\varepsilon)$  is an open subset of $\mathbb{K}$.
\end{lemma}

\Proof
Take $0 < \alpha < b$ sufficiently small such that $\Omega_\alpha$ is
compactly contained in $V$.
By Theorem \ref{aaa3} and \eqref{frffr}, every $f \in \mathbb{K}$
is well-defined on $\Omega_b({\mathcal K}_{f})$. Since by Lemma  
\ref{B33}
the map $f \mapsto {\mathcal K}_{f}$ is continuous in the Hausdorff  
metric,
there exists $\varepsilon_0 > 0$ such that if $f \in \mathbb{K}$ is  
such that
$\text{dist}_\mathbb{K}(f,g) < \varepsilon_0$ then
$\Omega_{ b/2}  ({\mathcal K}_{g}) \subset \Omega_{ b}  ({\mathcal  
K}_{f})$.
Since $\overline{V} \subset \Omega_{ b/2}  ({\mathcal K}_{g})$,
it follows that $f$ is well-defined on $V$, that is, $f \in  
\mathbb{A}_V$.
Hence there is a well-defined injection
$\mathbb{K}(g, \varepsilon_0) \to \mathbb{A}_V$.
Such injection is continuous.
Indeed, for $f \in \mathbb{K}(g, \varepsilon_0)$, the $C^0$ norm of $f$  
in
$\Omega_{b/2}({\mathcal K}_{g})$ is uniformly bounded,
while $\|f\|_{{\mathbb{A}_{\Omega_\alpha}}}$ varies continuously with  
$f$.
Since $\overline{\Omega_\alpha} \subset V \subset \overline{V} \subset
\Omega_{ b/2}  ({\mathcal K}_{g})$,
we deduce from Hadamard's three circles theorem
(see Lemma 11.5 on page 415 of \cite{lyubich}) that
$\|f\|_{\mathbb{A}_V}$ varies continuously with $f$ also.
Therefore the map $\mathbb{K}(g, \varepsilon_0) \to \mathbb{A}_V$
is continuous as asserted.
Now let $f \in \mathbb{K} \cap \mathbb{A}_V(g, \varepsilon_0)$.
Since the inclusion $\mathbb{A}_V \to \mathbb{A}_{\Omega_\alpha}$
has Lipschitz constant one, we have that $f \in \mathbb{K}(g,  
\varepsilon_0)$.
Hence, by continuity of the map $\mathbb{K}(g, \varepsilon_0) \to  
\mathbb{A}_V$, there
exists $\varepsilon_1 > 0$ such that $\mathbb{K}(f,\varepsilon_1)  
\subseteq
\mathbb{K} \cap \mathbb{A}_V(g, \varepsilon_0)$, which shows that this  
last
set is open in $\mathbb{K}$. This completes the proof.
\hfq

\begin{lemma}
\label{B2}
Let $b > 0$ be as defined in \eqref{frffr} and $\delta_0>0$ as in Lemma  
{\rm \ref{B3}.}
There exist
$\nu>0$, $0<\delta < \delta_0${\rm ,} a finite set ${\mathcal V}$
   of symmetric neighborhoods of $[-1,1]$ and a locally constant map
$\mathbb{K}\ni g\mapsto V_g\in \mathcal{V}$ with the following  
properties\/{\rm :}\/
\begin{enumerate}
\item
The neighborhood $V_g$ is compactly contained in $\Omega_{ b/2}   
({\mathcal K}_{g})${\rm ;}
\item  Every  $f \in  \mathbb{A}_{V_g}(g, \delta)$
is a quadratic-like map
with conformal modulus larger than $\nu${\rm ;}
\item
If  $f \in \mathbb{K} \cap \mathbb{A}_{V_g}(g, \delta)$ then
${\mathcal H}(f) \cap \mathbb{A}_{V_g}(g, \delta)$ is a codimension-one{\rm ,}
real-analytic submanifold
varying continuously with $f$.
\end{enumerate}
\end{lemma}

\Proof
For every $g \in \mathbb K$, let $U_g \subset \mathbb C$ be a symmetric
neighborood of $[-1,1]$ where $g$ is quadratic-like,
and take $n_g > 0$ so large that
  $V_g'=g^{-n_g}(U_g) \subset \Omega_{b/3} ({\mathcal K}_{g})$
and $V_g'\subseteq \hat V_g$, where $\hat V_g$ is as given in Theorem  
\ref{dfdfd}.

Let $\delta_g>0$ be so small that each $f \in \mathbb{A}_{V_g'}(g,  
\delta_g)$
is quadratic-like in $V_g'$ with conformal modulus greater than $\nu_g  
 > 0$
and also so that Theorem  \ref{dfdfd} holds true (for $V_g'$ and  
$\delta_g$).
By Lemma \ref{B33}, making $\delta_g$ smaller if necessary, we see that
$V_g'=g^{-n_g}(U_g) \subset \Omega_{b/2} ({\mathcal K}_{f})$
for all $f \in \mathbb{K} \cap  \mathbb{A}_{V_g'}(g, \delta_g)$.

By Lemma \ref{B55}, each set $\mathbb{K} \cap  \mathbb{A}_{V_g'}(g,  
\delta_g/2)$
is open in $\mathbb{K}$. Since $\mathbb{K}$ is compact, there exists a  
finite set
$
\left \{  g_i :i=1,\dots,l \right \}
$ such that
$$\mathbb K \subset \bigcup_{i=1}^l   
\mathbb{A}_{V_{g_i}'}(g_i,\delta_{g_i}/2)
\ .
$$
Thus we can set
$${\mathcal V}= \left\{ V_{g_i}': i=1,\dots,l  \right  \},\
\delta=\min_{i=1,\dots,l}
\{\delta_{g_i}/2 \}~~~~~{\rm and }~~~~~ \nu=\min_{i=1,\dots,l}  
\{\nu_{g_i} \}
\ .
$$
Therefore, since ${\mathbb K}$ is totally disconnected,
there exists a locally constant map
$\mathbb{K} \ni  g \mapsto V_g \in {\mathcal V}$
so that properties (i), (ii) and (iii) are satisfied.
\Endproof

We are now in a position to define the skew-product renormalization  
operator. This is
accomplished in our next lemma. Let us define first its range and  
domain, respectively,
as follows
\begin{eqnarray*}
E & = & \left \{(g,f): g \in \mathbb{K}~{\rm and}~
f \in \mathbb{A}_{V_g} \right \}, \\
E(\delta)  & = &   \left \{(g,f) \in E:
f \in \mathbb{A}_{V_g}(g, \delta)  \right \} \ .
\end{eqnarray*}

Let us now fix once and for all $a>0$ so small that
$\overline{\Omega_a} \subset V_g$ for every $g \in \mathbb{K}$
(this is possible because ${\mathcal V}$ in Lemma \ref{B2}
is a finite set).
The inclusion
   $k_g:\mathbb{A}_{V_g}\to \mathbb{A}_{\Omega_a}$
is a well-defined compact linear operator.
By \eqref{frffr} and Lemma \ref{B2} (i)
we also have
$$V_{R^N (g)} \subset
\Omega_{ b/2 }  ({\mathcal K}_{R^N (g)})   \subset
\Omega_{ b}  ({\mathcal K}_{R^N (g)})
\subset U_{g,N} \ .
$$
Therefore the inclusion
  $j_{g,N}:\mathbb{A}_{U_{g,N}} \to \mathbb{A}_{V_{R^N (g)}}$
is also a  well-defined compact linear operator.

\begin{lemma}
\label{B4} \label{B11}
  Let $\delta\!>\!0$ and $V_g \in {\mathcal V}$ be as in
Lemma {\rm \ref{B2}}. Let $N\!=\!N(a)\!>\!0${\rm ,}  $T_{g,N}$ and
  $T:{\mathbb O}_{\Omega_a} \to \mathbb {A}_{\Omega_a}$
be as in Lemma {\rm \ref{B3}.}
\begin{enumerate}
\item For every $g \in \mathbb{K}${\rm ,}
the operator
$T_{g}:\mathbb{A}_{\Omega_a}(g, \delta) \to \mathbb{A}_{V_{R^N (g)}}$
  given by
$T_{g} =j_{g,N} \circ T_{g,N}$ is  real analytic
with a compact derivative.
\item The   skew-product renormalization operator  $S: E (\delta) \to E$
given by $S(g,f)\break= \left (T(g),  S_g(f) \right)${\rm ,}
where     $S_g=T_{g}  \circ  k_g:\mathbb{A}_{V_g}(g, \delta)  
\to\mathbb{A}_{V_{T(g)}}${\rm ,}
is well-defined. Furthermore{\rm ,}\end{enumerate} \vglue-28pt
\begin{equation}
\label{fsdfe}
k_{T(g)} \circ S_g =   T    \circ  k_g \ .
\end{equation}
\end{lemma}
\vglue8pt

\Proof The proof is similar to the proof of Lemma \ref{B3} (ii).
\hfq

\Subsec{Hyperbolicity of the renormalization operator}
\label{ghjtyrert}
The purpose of this section is to show that $\mathbb{K}$ is a  
hyperbolic basic set
for the operator $T:\mathbb{O}_{\Omega_a} \to \mathbb{A}_{\Omega_a}$.
This will follow from the fact (Lemma \ref{B12} below) that the  
skew-product
renormalization operator has fiberwise real analytic stable manifolds  
and fiberwise one dimensional real
analytic unstable manifolds.

We start by noting that our skew-product operator satisfies the  
conditions W1--W4
on page 395 of Lyubich \cite{lyubich} in the real analytic case.
Namely, we have
\begin{itemize}
\item[W1.] The conformal modulus of each $g \in \mathbb K$ is larger  
than
a uniform constant $\mu > 0$.
\item[W2.] There exists $\eta > 0$ such that if dist$_\mathbb{K}(f,g) <  
\eta$
for some $f,g \in \mathbb K$, then $\mathbb{A}_{V_f}=\mathbb{A}_{V_g}$.
\item[W3.] There exists $\delta>0$ such that
$S_g(\mathbb{A}_{V_g}(g,\delta)) \subseteq \mathbb{A}_{V_{T(g)}}$.
\item[W4.] The vertical fibers $\mathcal{Z}_g$
(consisting of those normalized symmetric\break quadratic-like
germs whose external class is the same as that of $g$)
sit locally in $\mathbb{A}_{V_g}$ for each $g \in \mathbb{K}$.
\end{itemize}

Condition W1 is satisfied because of the complex bounds (Theorem
\ref{aaa3}). Condition W2 follows from Lemma \ref{B2}.    Condition W3
holds by the construction of $S_g$ in  Lemma \ref{B4}.
Condition W4 is a consequence of Lemma \ref{B2} (iii).

Now we have the following result.

\begin{lemma}
\label{B12}
The skew-product renormalization operator $S:E(\delta) \to E$ defined  
in Lemma {\rm \ref{B4}}
is fiberwise hyperbolic. Moreover
\begin{enumerate}
\item The local stable set $W_\delta^s(g)$ of $S$ at $g$ is a
codimension-one submanifold of $\mathbb{A}_{V_g}$ which is relatively  
open in
$\mathcal{H}(g) \cap \mathbb{A}_{V_g}(g,\delta)${\rm ,} and $W_\delta^s(g)$  
is tangent to
$E^s_g$ at $g$.
\item The local unstable set $W_\delta^u(g) \subset \mathbb{A}_{V_g}$  
of $S$ at $g$ is a
one-dimensional{\rm ,} real-analytic manifold{\rm ,} and $\{g\} \times W_\delta^u(g)$
varies continuously with $g\in \mathbb{A}$ in~$E$.
\end{enumerate}
\end{lemma}

\Proof
Since the operator $S$ satisfies Lyubich's conditions W1--W4 stated  
above,
part (i) follows from Theorem 8.2 on page 392 of \cite{lyubich} and
Theorem \ref{dfdfd}, and
part (ii) follows from
Theorem 8.8
on page 398 of Lyubich's paper \cite{lyubich}.
\hfq

\begin{theorem}
\label{B15}
Let   $T: \mathbb{O}_{\Omega_a}
  \to \mathbb{A}_{\Omega_a}$  be the real analytic operator
    defined in Lemma {\rm \ref{B4}.}  Then   there  is
  a continuous, $DT$-invariant
    splitting $\mathbb{A}_{\Omega_a}=E^s_g \bigoplus E^u_g${\rm ,}
   for
$g \in \mathbb K${\rm ,} such that
if $v^u \in E^u_g$ and  $v^s \in E^s_g$ then
\begin{eqnarray}
\label{unstable2}
\left \|DT^n(g) v^u
\right \|_{\mathbb{A}_{\Omega_a}}
& \ge & C^{-1} \theta^{-n} \|v^u \|_{\mathbb{A}_{\Omega_a}}, \\[4pt]
\label{unstable33}
\left \|DT^n(g) v^s
\right \|_{\mathbb{A}_{\Omega_a}}
& \le & C  \theta^{n} \|v^s \|_{\mathbb{A}_{\Omega_a}} \ ,
\end{eqnarray}
where $C>1$ and $0<\theta < 1$ are uniform constants  on $g$.
\end{theorem}

\Proof
Since for every $g \in \mathbb K$ the map
$k_g: {\mathbb A}_{V_g} \to {\mathbb A}_{\Omega_a}$ is linear and  
injective,
  it follows from Lemma \ref{B12} (ii) that  $Z_g^u=k_g(W^u_\delta (g))$
is a real-analytic, one-dimensional manifold varying continuously with  
$g$.
Let $w_g$ be the unitary vector tangent to $W^u_\delta (g)$ at $g$.
Then $v_g=k_g(w_g)$ is a vector tangent to $Z_g^u$ at $g$
and also varies continuously with $g$.
Since $k_g$ and $k_{T(g)}$ are linear maps we see from \eqref{fsdfe}
that if $\lambda_g$ is such that
$DS_g(g) w_g = \lambda_g w_{T(g)}$ then
$D  T (g)   v_g = \lambda_g   v_{T(g)}$.
Thus a natural candidate for $E^u_g$ is the
one-dimensional linear subspace generated by $v_g$.
In particular, (\ref{unstable2}) is satisfied.

Let us find the natural candidate for $E^s_g$.
We have that $D  T_{g} (g) v_g = w_{T(g)}$ and by hypothesis
$w_{T(g)}$ is transversal to the tangent space of $W^s_{\delta}(T(g))$.
Thus, by the implicit function theorem
$Z^s_g=T_g^{-1}(W^s_{\delta}(Tg))$ is a codimension-one manifold
transversal to $Z^u_g$.
Taking $E^s_g$ equal to the tangent space of
$Z^s_g$, we obtain that
$E^s_g \bigoplus E^u_g = {\mathbb A}_{\Omega_a}$.
By \eqref{fsdfe}, we have that
a neighborhood of $T(g)$ intersected with $Z^u_{T(g)}$  is contained in  
$T(Z^u_g)$  and
a neighborhood of $T(g)$ intersected with $T(Z^s_g)$ is contained in $  
Z^s_{T(g)}$, which implies that
the splitting $E^s_g \bigoplus E^u_g$ is
invariant under $DT$.
 From  assertion (i) in Lemma \ref{B12},
we obtain that $E^s_g$ varies continuously with $g$
and so
the splitting $E^s_g \bigoplus E^u_g$ also
varies continuously with $g$.

Finally, let $M \!>\!0$ be such that $\|DT_{g} (g)\|_{{\mathbb  
A}_{\Omega_a}} \!  \le \! M$
and note that $\|k_g\|_{{\mathbb A}_{\Omega_a}}   \!\le  \!1$
for all $g \in \mathbb K$\/.
For all $v^s \in E^s_g$ with unit norm, let $u^s=D  T_{g} (g) v^s \in  
{\mathbb A}_{V_{T(g)}}$.
By Lemma \ref{B12} (i) and \eqref{fsdfe} there exists $C_1>1$ and  
$0<\theta<1$
such that
\begin{eqnarray*}
\|DT^n(g)v^s\|_{{\mathbb A}_{\Omega_a}}
& = & \| k_{T^n (g)} \circ
DS_{T^{n-1} (g)} (T^{n-1} (g))
\circ \dots \circ
DS_{T(g)} (T(g))  u^s \|_{{\mathbb A}_{\Omega_a}}  \nonumber \\
& \le & \|
DS_{T^{n-1} (g)} (T^{n-1} (g))
\circ \dots \circ
DS_{T(g)} (T(g))  u^s \|_{{\mathbb A}_{V_{T^n(g)}}}
\nonumber \\
& \le & C_1 M \theta^{n-1} \ ,
\end{eqnarray*}
which shows that \eqref{unstable33} is satisfied. This completes the  
proof.
\Endproof

With the above results, we have therefore established  
Theorem~\ref{lyubhyp},
to the effect that a suitable power of the renormalization operator is  
indeed
hyperbolic in a suitable (real) Banach space of real analytic mappings.  
From
now on, we shall concentrate on the problem of extending such  
hyperbolicity to
larger ambient spaces of smooth mappings. Our journey will take us far  
into
the wilderness of nonlinear functional analysis.

\section{Extending invariant splittings}
\label{hjerxc}
In this section we prove a certain result from functional analysis  
(Theorem
\ref{11ttt} below) that is absolutely crucial for the stable manifold  
theorem that
we shall prove later. This result deals with the notion of compatibility
presented below and is a strong generalization of a key idea of
Davie in \cite{davie}.
In Section~\ref{hfgrtyresd}, we shall use the results presented here to show  
that the invariant
splitting for the renormalization operator $T$ in  
$\mathbb{A}_{\Omega_a}$ of
Section~\ref{grtsfd}  extends to an invariant splitting for the action of $T$
in the larger spaces $\mathbb{A}^r$ of $C^r$ maps.

\Subsec{Compatibility}
\label{sub:compat}
We are interested in the answer to the following question. Given a  
smooth
operator $T:\mathcal{O}\to\mathcal{A}$ having a hyperbolic basic set
$\mathbb{K}$, and given a larger ambient space $\mathcal{B}\supseteq
\mathcal{A}$ to which $T$ extends (not necessarily smoothly),
under which conditions does $\mathbb{K}$ have a hyperbolic structure in
$\mathcal{B}$?  To give a precise meaning to this question (and then  
answer
it!) we introduce the following notion.

We have a natural continuous  map
$L:\mathbb{K}\to \mathcal{L}(\mathcal{A},\mathcal{A})$ given by
\begin{eqnarray*}
\mathbb{K}\ni x \ \mapsto && \, L_x :\mathcal{A}\to \mathcal{A}\\
&& L_x(v)=DT(x)\,v
\ .
\end{eqnarray*}
We will also assume that for every $x \in \mathbb{K}$,
$E_x^u$ is a one-dimensional subspace and  that we can chose
  a unit vector  $\mathbf{u}_x \in E_x^u$ varying continuously with $x$
   so that $L_x(\mathbf{u}_x)= \delta_x \cdot \mathbf{u}_{T(x)}$
with $\delta_x > 0$. In the case of the renormalization operator there  
is a
natural choice for the vectors $\mathbf{u}_x $: choose the unit vector  
pointing
in the direction of increasing topological entropy.

For every $x \in \mathbb{K}$, we denote
$DT^n(x)=L_{T^{n-1}(x)}\circ\cdots\circ L_x$ by
$L_x^{(n)}$ and $\delta_{T^{n-1}(x)}  \cdots \delta_x$ by
$\delta_x^{(n)}$.
By hyperbolicity of $\mathbb{K}$, there exist $C_0>0$ and
$\lambda>1$ such that for every $x \in \mathbb{K}$ and every $n \ge 1$  
we have
\begin{equation}
\label{fggtrs}
\delta_x^{(n)} >  C_0 \lambda^n \ .
\end{equation}

We denote by  $\mathcal {X}(r)$ the open ball  in the Banach space
$\mathcal {X}$  centered at the origin and with radius $r >0$.

\begin{definition}
\label{robust}
Let $\theta < \rho < \lambda$ where $\theta$ is the  contraction
  exponent of the hyperbolic basic set  $\mathbb{K}$  of the operator $T$
  and $\lambda$ is as in \eqref{fggtrs}.
  The pair $(\mathcal {B},
\mathcal {C})$ is {\it $\rho$-compatible} with $(T, \mathbb{K})$ if the  
following
conditions are satisfied.
\begin{enumerate}
\item[$\mathbf{A1.}$]
The inclusions $\mathcal {A} \to \mathcal {B}\to \mathcal {C}$ are
compact operators.
\item[$\mathbf{A2.}$]
  There exists $M>0$ such that each linear operator $L_x= DT(x)$ extends
to a linear operator $\hat{L}_x :\mathcal{C}\to \mathcal{C}$ with
\begin{eqnarray*}
\left\|\hat L_x \right\|_{\mathcal{C}}  &  <  & M   ,\\
\hat L_x(\mathcal {B}) & \subset & \mathcal {B} ,    \\
\left\|\hat L_x(v)\right\|_{\mathcal{B}} & < & M\|v\|_{\mathcal{B}}.
\end{eqnarray*}
\item[$\mathbf{A3.}$]
  The map $\tilde{L}:\mathbb{K} \to \mathcal{L}(\mathcal{B},  
\mathcal{C})$
given by $\tilde{L}_x=\hat{L}_x |\mathcal {B}$
is continuous.
\item[$\mathbf{A4.}$]
There exists $\Delta>1$ such that $\mathcal {B}(\Delta)\cap \mathcal
{A}$ is $\mathcal {C}$-dense in $\mathcal {B}(1)$.
\item[$\mathbf{A5.}$]
There exist $K>1$ and a positive integer $m$ such that
$$
\left\|\hat L_x^{(m)} (v) \right\|_\mathcal {B} \leq \max \left\{
\frac{\rho^m}{2(1+\Delta)}\|v\|_{\mathcal {B}}, K  
\|v\|_{\mathcal{C}}\right\}
\ .
$$
\end{enumerate}
\end{definition}

\begin{remark} Note that neither the map $\hat L:\mathbb{K}\times
\mathcal{C}\to \mathbb{K}\times \mathcal{C}$ given by $\hat L (x,v)=
(T(x),\hat L_x(v))$ nor its restriction
from $\mathbb{K}\times \mathcal{B}$ to $\mathbb{K}\times \mathcal{B}$  
are
necessarily continuous.
\end{remark}

\begin{example}
As we know from Theorem \ref{lyubhyp}, $\mathbb K$ is a
hyperbolic basic set of the
renormalization operator
$T=R^N:\mathbb{O} \to \mathbb{A}$.
In Section~\ref{hfgrtyresd} (see Theorem \ref{dfgdsfgg}),
we will show that the pair $(\mathbb{A}^r,\mathbb{A}^0)$ is
$\rho$-compatible for $r$   sufficiently close to $2$ and
$1$-compatible for $r>2$ noninteger.
\end{example}

Let $\pi_x^u:\mathcal{A} \to  E_x^u$ and
$\pi_x^s:\mathcal{A} \to  E_x^s$ be the canonical projections.
We define  $P_x=\pi_{T(x)}^u \circ L_x$ and $Q_x=\pi_{T(x)}^s \circ L_x$
which have the property  that
   $L_x=P_x + Q_x$ and that
$P_{T(x)}  Q_x = Q_{T(x)}  P_x = 0$.
We also define the  linear functional
  $\sigma_x:\mathcal{A} \to \mathbb{R}$
   by $\pi_x^u(v)=\sigma_x(v)\mathbf{u}_x$, and observe that $P_x(v) =
   \delta_x  \sigma_x(v) \mathbf{u}_{T(x)}$.
We note that the map $\sigma:\mathbb{K} \to \mathcal{L}(\mathcal{A},
\mathbb{R})$ which associates to each $x$ the linear functional
$\sigma _x$ is continuous.

\begin{theorem}
\label{11ttt} If $(\mathcal {B},\mathcal {C})$ is $\rho$-compatible with
$(T, \mathbb{K})$ then each stable functional $\sigma_x$ extends to a
unique linear functional
$\hat \sigma_x \in \mathcal {B}^*$ satisfying
\begin{equation}
\label{ineqsig}
\left \|\hat L_x^{(n)}(v) - \delta_x^{(n)} \hat \sigma_x(v)
\mathbf{u}_{T^n(x)} \right \|_\mathcal{B} \leq C \hat \theta^n
\|v\|_\mathcal{B}
\end{equation}
for some $C>0$ and $0<\hat \theta <\rho $.
Furthermore{\rm ,} the map $\hat{\sigma}:\mathbb{K} \to  
\mathcal{L}(\mathcal{B},
\mathbb{R})$ which associates to each $x$ the linear functional
$\hat{\sigma} _x$ is continuous.
\end{theorem}
\Proof
Let $m$ and $M$ be as given in Definition~\ref{robust}.
Since by property $\mathbf{A1}$ the $\mathcal{C}$-closure of
$\mathcal{B}(1)$ is compact, and since by property $\mathbf{A4}$ the
intersection
$\mathcal{A}\cap\mathcal{B}(\Delta)$ is $\mathcal{C}$-dense in
$\mathcal{B}(1)$ we can find a finite set
\begin{equation*}
\Phi\subseteq \mathcal{A}\cap\mathcal{B}(\Delta)
\end{equation*}
such that for each $w\in\mathcal{B}(1)$ there exists $w'\in \Phi$ such  
that
\begin{equation*}
\|w-w'\|_{\mathcal{C}}\;<\;
  \frac{\rho^m}{4K}
\ .
\end{equation*}
Now let $v\in \mathcal{B}(1)$, and let $v_{0}\in \Phi$ be such that
\begin{equation*}
\|v-v_{0}\|_{\mathcal{C}}\;<\;
  \frac{\rho^m}{2K}
\ .
\end{equation*}
Since $\|v-v_{0}\|_{\mathcal{B}}< 1+\Delta$,
applying the inequality of property $\mathbf{A5}$ to $v-v_{0}$ yields
\begin{eqnarray*}
\left \|\hat{L}_x^{(m)}(v-v_{0}) \right \|_{\mathcal{B}}
& \le &
\max \left  \{
  \frac{\rho^m}{2(1+\Delta)}  \|v-v_{0}\|_{\mathcal{B}}, K
  \|v-v_{0}\|_{\mathcal{C}}     \right \}  \\
& < &
\rho^m/2  \ .
\end{eqnarray*}
Therefore $\hat{L}_x^{(m)}(v)=\hat{L}_x^{(m)}(v_{0})+(\rho^m/2)w_1$ for
some $w_1\in \mathcal{B}(1)$.
Repeating the argument with $w_1$ replacing $v$ and proceeding
inductively in this fashion, we get after $k$ steps
\begin{equation*}
\hat{L}_{x}^{(km)}(v)\;=\;
\sum_{j=0}^{k-1}\frac{\rho^{jm}}{2^j}L_{T^{jm}(x)}^{((k-j)m)}(v_j)+
\frac{\rho^{km}}{2^k}w_k
\end{equation*}
for some $w_k\in \mathcal{B}(1)$ and $v_j\in \Phi$.
Now recall that
\begin{equation*}
\begin{split}
L_{T^{jm}(x)}^{((k-j)m)}(v_j)\;&=\;P_{T^{jm}(x)}^{((k-j)m)}(v_j)
+Q_{T^{jm}(x)}^{((k-j)m)}(v_j) \\
\;&=\; \delta_{T^{jm}(x)}^{((k-j)m)}  
\sigma_{T^{jm}(x)}(v_j)u_{T^{km}(x)}
+Q_{T^{jm}(x)}^{((k-j)m)}(v_j)
\ .
\end{split}
\end{equation*}
Hence we can write
\begin{eqnarray}\label{sigma}
\hat{L}_{x}^{(km)}(v)&=  
& \delta_{x}^{(km)}
\left(\sum_{j=0}^{k-1}\frac{1}{2^{j}}
\frac{\rho^{jm}}{\delta_x^{(jm)}}
  \sigma_{T^{jm}(x)}(v_j)\right)\, \mathbf{u}_{T^{km}(x)}\\
&&+  \sum_{j=0}^{k-1}\left(\frac{\rho^{m}}{2}\right)^{j}
Q_{T^{jm}(x)}^{((k-j)m)}(v_j) +
\frac{\rho^{km}}{2^k}w_k\ .
\nonumber
\end{eqnarray}
The first summation in parentheses converges to a limit because, by
\eqref{fggtrs}, $\left|\delta_x^{(jm)}\right |\ge  
C_0\lambda^{jm}>C_0\rho^{jm}$
and $\left \{\sigma_{T^{jm}(x)}(v_j) \right \}$ is bounded, as the  
$v_j$ run
through finitely many values and
$\left \|\sigma_{T^i(x)} \right \|\leq M$ for all $i$.
We therefore define
\begin{equation}
\label{qsqsqsq}
\hat\sigma_x(v)\;=\;\lim_{k\to\infty}\,
\sum_{j=0}^{k-1}\frac{1}{2^{j}}
\frac{\rho^{jm}}{\delta_x^{(jm)}}
\sigma_{T^{jm}(x)}(v_j)
\ .
\end{equation}
It will be clear in a moment that this extension of $\sigma_x$ is
  independent of the choices of approximants $v_j$ performed above,
linear,
continuous, and the unique extension satisfying
\eqref{ineqsig}. We know that
$$
\left \|Q_{T^{jm}(x)}^{((k-j)m)} \right \|_{\mathcal{A}}\le
C_1\theta^{m(k-j)}
$$
for all $j<k$. Thus the second summation plus the last term in
{\eqref{sigma}} add up to a vector with $\mathcal{B}$-norm
bounded by
\begin{equation*}
C_2\left[\sum_{j=0}^{k 
-1}\frac{1}{2^j}\left(\frac{\theta}{\rho}\right)^{m(k-j)}
+\frac{1}{2^k}\right]\rho^{km}
\ .
\end{equation*}
  This gives
\begin{equation}\label{theta}
\left\|\hat{L}_{x}^{(km)}(v)-\delta_x^{(km)}
\hat\sigma_x(v)\mathbf{u}_{T^{km}(x)}\right\|_{\mathcal{B}}\le\;
C_3(k+1)\beta^{k}\rho^{km}
\ ,
\end{equation}
where $\beta=\max \{1/2,\theta/\rho\}<1$.
Now choose $0<\hat\theta<\rho $ so that  
$(k+1)\beta^k<C_4(\frac{\hat\theta}{\rho})^{km}$
for all $k$.  Since by property $\mathbf{A2}$
for all $x \in \mathbb{K}$ and for all $v \in \mathcal{B}$ we have
$\|\hat{L}_x (v) \|_{\mathcal{B}} < M  \|v\|_{\mathcal{B}}$, writing  
$n=km+r$
and using the  above estimates we obtain the
desired inequality {\eqref{ineqsig}}.

Let us now verify that $\hat\sigma_x(v)$ is the unique value satisfying
{\eqref{ineqsig}}.
In particular, it does not depend
  on the choices of approximants $v_j$ taken in \eqref{qsqsqsq}.
To do this we represent by $\sigma_x^*(v)$ a value
satisfying {\eqref{ineqsig}}, for instance that
obtained in
\eqref{qsqsqsq} by taking another choice of   approximants. Therefore  
we have
\begin{eqnarray*}
\left \|\hat\sigma_x(v)\mathbf{u}_{T^{km}(x)} -
\sigma_x^*(v)\mathbf{u}_{T^{km}(x)}
\right\|_{\mathcal{B}}
& \le &
\left\|
\hat\sigma_x(v)\mathbf{u}_{T^{km}(x)}-\left(\delta_x^{(km)}\right)^{-1}
\hat{L}_{x}^{(km)}(v) \right \|_{\mathcal{B}}  \\
&&  \ +
\left \|\left(\delta_x^{(km)}\right)^{-1}
\hat{L}_{x}^{(km)}(v) - \sigma_x^*(v)\mathbf{u}_{T^{km}(x)}
\right\|_{\mathcal{B}}   \\
& \leq &
2C   \hat{\theta}^{km}\left(\delta_x^{(km)}\right)^{-1} \ .
\end{eqnarray*}
Letting $k\to\infty$ in this inequality we
deduce  that $\hat\sigma_x(v)=\sigma_x^*(v)$.
A similar argument shows that   $\hat\sigma_x$ is
linear. Using inequality \eqref{theta} with $k=1$, we obtain that
$\|\hat\sigma_x\|_{\mathcal{B}}$
is bounded. Finally, the fact that $\hat\sigma_x$ is continuous in $x$  
can be
deduced from \eqref{qsqsqsq} using property $\mathbf{A3}$.
\Endproof

\begin{corollary}
\label{corrob}
Let $(\mathcal {B},\mathcal {C})$ be $\rho$-compatible with
$(T, \mathcal {A})$.
Let the linear
functional  $\hat \sigma_x \in \mathcal {B}^*$
be the   extension  of the stable functional $\sigma_x$
satisfying inequality \eqref{ineqsig} for all $x \in \mathbb{K}$.
Then{\rm ,}   there exists a continuous splitting
$\mathcal {B}=\hat{E}_x^s \bigoplus \hat{E}_x^u$
   with the following properties\/{\rm :}\/
\begin{enumerate}
\item $\hat{E}_x^u$ is the inclusion
in $\mathcal {B}$ of
the unstable linear space $E_x^u \subset \mathcal {A}${\rm ;}
\item  $\hat{E}_x^s= {\rm Ker} (\hat{\sigma}_x )${\rm ;}
\item The  splitting is invariant by $\hat{L}_x${\rm ;}
\item There exists a constant  $C>0$
such that
$$
\left \|\hat{L}_x^m\,(v) \right \|_\mathcal{B} \geq C \lambda^m  
\|v\|_\mathcal{B}
$$
for all  $x \in \mathbb{K}${\rm ,}  for all
$v\in \hat{E}^u_x${\rm ,} and for all   $m \in \mathbb{N}$
  {\rm (}where $\lambda> 1$ is the same as in  \eqref{fggtrs}{\rm );}
\item There exist constants $C>0$ and $0 < \hat{\theta} < \rho <  
\lambda$
such that
$$
\left \|\hat{L}_x^m\,(v) \right \|_\mathcal{B}
  \leq C \hat{\theta}^m \|v\|_\mathcal{B}
$$
for all  $x \in \mathbb{K}${\rm ,} for all $v\in \hat{E}^s_x${\rm ,}
and for all   $m \in \mathbb{N}$.
  In particular{\rm ,} if $\rho \le 1$ then $\hat{\theta} <1$.
\item  Let  $\hat{\pi}^s_x:\mathcal{B} \to \hat{E}_x^s$
and $\hat{\pi}^u_x:\mathcal{B} \to \hat{E}_x^u$
be the natural projections such that
$$\hat{\pi}^s_x \circ \hat{\pi}^s_x  =\hat{\pi}^s_x, \
\hat{\pi}^u_x \circ \hat{\pi}^u_x=\hat{\pi}^u_x \ \text{and} \
   \hat{\pi}^s_x \circ \hat{\pi}^u_x=
    \hat{\pi}^u_x \circ\hat{\pi}^s_x = 0 \ .
$$
The operators
$\hat{Q}_x:\mathcal{B} \to \mathcal{B}$ and
$\hat{P}_x:\mathcal{B} \to \mathcal{B}$ may be defined by
$\hat{Q}_x^m=\hat{L}_x^m \circ \hat{\pi}^s_x$ and
by $\hat{P}_x^m=\hat{L}_x  \circ \hat{\pi}^u_x$.
Then{\rm ,} there exists $C>1$ such that
$$
\left \|\hat{Q}_x^m \right \|_\mathcal{B}   \leq    C \hat{\theta}^m\hbox{ and }
\left \|\hat{P}^m \right \|_\mathcal{B}    \geq   C^{-1} \lambda^m
$$
for all   $x \in \mathbb{K}$ and for all $m \in \mathbb{N}$.
\end{enumerate}
\end{corollary}

\Proof
First, we observe that   for all $x \in \mathbb{K}$
   and
for all $v\in \mathcal {B}$, we can write
  $v= (v- \hat{\sigma}_x (v) \mathbf{u}_x)+ \hat{\sigma}_x (v)  
\mathbf{u}_x$,
where $v- \hat{\sigma}_x (v) \mathbf{u}_x \in {\rm Ker} (\hat{\sigma}_x  
)$
and $\hat{\sigma}_x (v) \mathbf{u}_x \in \hat{E}_x^u$.
Since  $\hat{\sigma}_x (\mathbf{u}_x) = 1 \ne 0$,
we obtain that $\mathcal{B} = \hat{E}_x^s \bigoplus \hat{E}_x^u$.

By inequality \eqref{ineqsig}, there exists $C > 0$ such that
\begin{equation}
\label{gedx}
\left \|\hat L_x^{(m)}(v) \right \|_\mathcal{B} \le  C   \hat \theta^m
\|v\|_\mathcal{B} \end{equation}
for all $x \in \mathbb{K}$,
   for all  $v\in {\rm Ker} (\hat{\sigma}_x )$, and for all $m \in  
\mathbb{N}$;
  if  $v\in \mathcal {B} \setminus {\rm Ker} (\hat{\sigma}_x )$ then
there exists $C_v > 0$ such that
$$
\left \|\hat L_x^{(m)}(v) \right \|_\mathcal{B} \ge  C_v \lambda^m   \ .
$$
Therefore,  $\hat L_x({\rm Ker} (\hat{\sigma}_x )) \subset
{\rm Ker} (\hat{\sigma}_{T(x)} )$.
Since $L_x(E^u_x)=E^u_{T(x)}$ implies that
  $\hat{L}_x(\hat{E}^u_x)=\hat{E}^u_{T(x)}$, the  splitting is invariant  
by $\hat{L}_x$.

By inequality \eqref{gedx}, we obtain that property (v)
and  the first inequality  in  property (vi) are   satisfied.
Since $\mathbb{K}$ is compact and  the map
$\mathbb{K} \to \mathbb{R}^+$
which associates $\|\mathbf{u}_x\|_{\mathcal{B}}$ to each $x$ is  
continuous,
there is $C > 1$ such that
  \begin{equation}
\label{uuuuu}
C^{-1}  \|v\|_{\mathcal{A}} <
    \|v\|_{\mathcal{B}}
<C  \|v\|_{\mathcal{A}}
\end{equation}
for all $x \in \mathbb{K}$ and for all $v\in E^u_x$.
Thus,   property (iv) and the second  inequality  in property (vi)
follow from \eqref{fggtrs} and \eqref{uuuuu}.
\hfq

\section{Extending the invariant splitting for renormalization}
\label{hfgrtyresd}

Our aim in this section is to show that the invariant splitting on the  
limit
set $\mathbb{K}$ of the  operator $T$ given by Theorem~\ref{lyubhyp},  
which is
an iterate of the renormalization operator, can be extended to an  
invariant
splitting of the same operator acting in the space of
$C^r$ unimodal maps. Given the abstract results of the previous section,
namely Theorem~{\ref{11ttt}} and Corollary~\ref{corrob}, all we have to  
do is
find the appropriate compatible spaces and the corresponding  
compatibility
constants. More precisely, we shall prove the following theorem.

\begin{theorem}
\label{dfgdsfgg}
Let $T$ and $\mathbb{K}$ be as above{\rm ,} and let $\lambda$ be the expansion
constant satisfying \eqref{fggtrs}.
\begin{enumerate}
\item For all $\alpha>0$ the pair of spaces
$(\mathbb{A}^{2+\alpha}, \mathbb{A}^0)$ is
$1$-compatible with $(T,\mathbb{K})$.
\item For all $1<\rho<\lambda$ there exists $\alpha>0$ sufficiently  
small such
that $(\mathbb{A}^{2-\alpha}, \mathbb{A}^0)$ is
$\rho$-compatible with  $(T,\mathbb{K})$.
\end{enumerate}
\end{theorem}

The path towards the proof of this theorem (presented in \S  
\ref{sec:deriv})
leads us to perform what amounts to a {\it spectral analysis} of the  
formal
derivative of the renormalization operator, which in turn calls for  
certain
estimates on the geometry of the post-critical set of each map in the
limit set of renormalization.
We have the following explicit formula for the
derivative $L_f=DT(f)$ of $T$ at $f\in\mathbb{K}$:
\begin{eqnarray*}\label{deetee}
 DT(f)v&=&
\frac {1}{\lambda_f}
\sum_{j=0}^{p-1}Df^j(f^{p-j}(\lambda_f x))v(f^{p-j-1}(\lambda_f x))  \\
&&
+\frac {1}{\lambda_f}
[x(Tf)'(x)-Tf(x)]
\sum_{j=0}^{p-1}Df^j(f^{p-j}(0))v(f^{p-j-1}(0))
\ ,
\end{eqnarray*}
where as before $\lambda_f=f^p(0)$ for some positive integer $p=p(f,N)$.
We observe that the operator $L_f$ extends naturally to
each of the spaces $\mathbb{A}^\gamma$ for $\gamma \ge 0$.

Properties $\mathbf{A1}$, $\mathbf{A2}$ and $\mathbf{A3}$
of Definition {\ref{robust}} are easily verified in
our setting. Property $\mathbf{A4}$ follows from a general result of  
H\"older
spaces that can be proved via smoothing operators. Hence, the heart of  
the
matter is verifying property $\mathbf{A5}$. This is where the geometric  
scaling
properties of the invariant Cantor set of a map in $\mathbb{K}$ become
important -- see \S \ref{sec:bdedgeom}.
We follow  Davie's observation  that $L_f^{(m)}$ is a special sort of  
operator
-- what we call an $L$-operator -- which is amenable to analysis.
The verification of the fifth property
  (with $(\mathcal{B},\mathcal{C})=(\mathbb{A}^\gamma,\mathbb{A}^0)$)
  -- presented in \S
\ref{sec:deriv} -- consists in controlling the norm of a
certain positive linear operator $L_{f,\gamma}^{(m)}: \mathbb{A}^0 \to
\mathbb{A}^0$ associated to $L_f^{(m)}$
(see Lemma \ref{spectral}).
Using the
bounded distortion properties of $f \in \mathbb{K}$ and the geometry of  
the
invariant set of $f$, we show that the exponential growth rate of the  
$C^0$
norm of $L_{f,\gamma}^{(m)}$
  is bounded by some $\mu<\lambda$ if $\gamma=2-\alpha$ with $\alpha>0$  
small
enough and is bounded by some $\mu<1$ if $\gamma=2+\alpha$ with  
$\alpha>0$.

\Subsec{\older\ norms and $L$-operators}
\label{sec:holder}
First we define what we mean by an $L$-operator, and to each such
operator $L$ we associate another operator
$L_\gamma$, acting on continuous functions.
  Then, we use local H\"older  estimates
to control the  norm of compositions
$L_{m} \circ \cdots \circ L_2 \circ L_1$
of $L$-operators $L_i$ by the norm of
$(L_{m} \circ \cdots \circ L_2 \circ L_1)_\gamma$.

\begin{definition}
An $L$-operator is a {\it bounded linear operator}     $L:C^\gamma(I) \to
C^\gamma(I)$ that can be written in the form
$$
L v (x)=\sum_{i=1}^n \phi_i(x) v (\psi_i(x)),
$$
where $\phi_i \in C^{\gamma_1}(I)$
and $\psi_i \in C^{\gamma_2}(I)$ are maps
such that $\psi_i(I) \subset I$ for $i=1,\dots ,n$, and where
$\gamma_1 > 0$ and $\gamma_2 \ge 1$ are such that  $0< \gamma <
\gamma_1, \gamma_2$.
\end{definition}

\begin{example}
For  all $f \in \mathbb{K}$ and all $i \ge 0$, the {\it formal derivative}
$L_f=DT(T^i(f))$ is an $L$-operator.
\end{example}

An $L$-operator $L$ as above yields  a positive,  bounded linear
operator $L_\gamma:C^0(I) \to C^0(I)$ defined by
$$
L_\gamma v(x) = \sum_{i=1}^n|\phi_i(x)\|D\psi_i(x)|^\gamma
  v (\psi_i(x)) \ .
$$

A straightforward computation yields the following result.

\begin{lemma}
\label{trivial}
If $L_1,L_2: C^\gamma(I) \to C^\gamma(I)$ are
$L$-operators{\rm ,} then $(L_1\circ L_2)_\gamma=L_{1,\gamma}\circ
L_{2,\gamma}$.
\end{lemma}

We remind the reader that a function $\varphi:I \to I$ is {\it
$\alpha$-\older\ continuous},
for a fixed $0 < \alpha <1$, if
there is $c>0$ such that
$ |\varphi(x)-\varphi(y)| \le c |x-y|^\alpha$ for all $x,y \in I$.
Let  $C^{\alpha}(I)$ be the Banach space of all
$\alpha$-H\"older continuous  real functions on $I$,
with norm
$$
\|\varphi\|_\alpha =\max\left\{\|\varphi\|_0\, , \,\sup_{x\ne y}
  \frac{|\varphi(x)-\varphi(y)|}{|x-y|^\alpha}\right\}\ .
$$
Let $C^{k+\alpha}(I)$ be the Banach space of all real functions on $I$
for which the $k^{\rm th}$ derivative is $\alpha$-H\"older continuous, with
norm
$$
\|\varphi\|_{k+\alpha}=\max\{\|\varphi\|_0,\|D^k\varphi\|_{\alpha}\}.
$$

\begin{lemma} \label{spectral}
Let $L_i:C^\gamma(I) \to C^\gamma(I)$ be a sequence of
$L$-operators{\rm ,} and assume that there exist constants $\mu>0$ and $C>0$
such that for all  $n${\rm ,}
\begin{equation}
\label{exp}
\left\| (L_n\circ \cdots\circ L_2\circ L_1)_\gamma  \right\|_0 \leq
C\mu^n
\ .
\end{equation}
Then for all $\rho>\mu$ and all $\varepsilon >0$ there exist $m>0$ and
$K>0$ such that for all $ v \in C^\gamma(I)${\rm ,}
$$
\left\| L_m\circ \cdots \circ L_2\circ L_1 ( v ) \right\|_\gamma \leq
\max
\left \{   \epsil \rho^m \| v \|_\gamma   , K \| v \|_0 \right \} \ .
$$
\end{lemma}

To prove the above proposition,  we will use local H\"older estimates  
for\break
$L$-operators given in our next lemma.
For  each $\eta>0$ and each $\varphi  \in  C^{\alpha}(I)$, we consider  
an
associated semi-norm
$$
\|\varphi\|_{\alpha,\eta} =\sup_{0<| x - y|< \eta}
  \frac{|\varphi(x)-\varphi(y)|}{|x-y|^\alpha}\ .
$$
The corresponding semi-norm  of $\varphi \in
C^{k+\alpha}(I)$  for $k>0$ is $\|\varphi\|_{k+\alpha,\eta}=
\|\varphi^{(k)}\|_{\alpha,\eta}$.

\begin{lemma}
\label{holder}
Let  $L:C^\gamma(I) \to C^\gamma(I)$ be an $L$-operator
as defined above.
\begin{enumerate}
\item For every $\varepsilon >0${\rm ,} there exists $\eta >0$ such that
$$
\|L v \|_{\gamma,\eta} \le (\varepsilon+   \|L_{\gamma}\|_{C^0(I)} )
\| v \|_{\gamma}.
$$
\item For every $\varepsilon >0$ and $0< \xi < \gamma${\rm ,} there is
$\eta>0$ such that
$$
\|L v \|_{\xi,\eta} \le \varepsilon \| v \|_{\gamma}.
$$
\end{enumerate}
\end{lemma}
\Proof
See Lemmas 1 and 2 in \cite{davie}.
\Endproof

{\it Proof of Lemma} \ref{spectral}.
Choosing $m$ such that $C\mu^m < \epsil \rho^m / 8$,  we have
$$
M=\left\|L_{m,\gamma}\circ\cdots \circ L_{1,\gamma}\right\|_0<
\frac{\varepsilon\rho^m}{8}
\ .
$$
By Lemma~\ref{holder}, given $\varepsilon'=\varepsilon\rho^m/8$, there
exists $\eta>0$ such that
$$
\left\|L_{m}\circ\cdots \circ L_{1}( v )\right\|_{\gamma,\eta}\leq
\left(\varepsilon'+M\right)\| v \|_\gamma\leq
\frac{\varepsilon\rho^m}{4}\| v \|_\gamma
\ .
$$
Taking $K=8k!\|L_{m}\circ\cdots \circ L_{1}\|_0/\eta^\gamma$, writing
$\gamma=k+\alpha$, where $k$ is an integer and $0< \alpha < 1$, and
using interpolation of norms (see Lemma 4 in \cite{davie}), we deduce  
that
\begin{eqnarray*}
\left\| L_{m}\circ\cdots \circ L_{1}( v ) \right\|_\gamma
& \leq &
4 \max
\left \{  \left \|L_{m}\circ\cdots \circ L_{1}( v )\right  
\|_{\gamma,\eta}  ,
\frac{2k!}{\eta^\gamma}   \left \|L_{m}\circ\cdots \circ L_{1}(v)
\right \|_0 \right \} \\
& \leq &
\max
\left \{   \epsil \rho^m \| v \|_\gamma  ,  K  \| v \|_0 \right \} \ .
\end{eqnarray*}
\vglue-18pt\Endproof
 
\Subsec{Bounded geometry}
\label{sec:bdedgeom}
Our aim in this section is to prove two crucial propositions concerning  
the
geometry of the invariant Cantor set of an infinitely renormalizable  
map in
the limit set of renormalization.
They are important not only in the proof of Theorem \ref{dfgdsfgg},
but also in the proof (presented in \S \ref{sec:robust}) that the  
renormalization
operator is robust (in the sense of \S \ref{sec:stable}).

We recall our notation.
For each $f\in\mathbb{K}$, let $\mathcal{I}_f\subseteq I$ be the
closure of the postcritical set of $f$ (the Cantor
attractor of $f$). For each $k\ge 0$, we can write
$$
R^k f (x) = \frac{1}{\lambda_k} \cdot f^{p_k} (\lambda_k x)
$$
where $p_k=\prod_{i=0}^{k-1} p(R^if)$ and
$\lambda_k=\prod_{i=0}^{k-1}  \lambda (R^if)$.
Recall that  the renormalization intervals
$\Delta_{0,k}=[-|\lambda_k|,|\lambda_k|]
\subset [-1,1]$, and
$\Delta_{i,k}=f^i(\Delta_{0,k})$ for
$i=0,1, \dots, p_k-1$.
The collection
${\mathbf C}_k=\{\Delta_{0,k}, \dots, \Delta_{p_k-1,k}\}$
consists of pairwise disjoint intervals.
Moreover,
$\bigcup \{\Delta:\Delta\in {\mathbf C}_{k+1}\} \subseteq \bigcup
\{\Delta:\Delta\in {\mathbf C}_k\}$ for all $k \ge 0$ and we have
$$
\mathcal{I}_f=
\bigcap _{k=0}^{\infty}
\bigcup_{i=0}^{p_k-1}
\Delta_{i,k} \ .
$$

In our first proposition, $f$ is a normalized, symmetric quadratic  
unimodal
map, infinitely renormalizable, sufficiently smooth (say $C^2$)
for Sullivan's real bounds to
be true for $f$. But there are no restrictions on the combinatorics.
We shall use the general fact, due to Guckenheimer \cite{guck}, that  
among
those renormalization intervals at the $k^{\rm th}$ level the one that  
contains
the critical point of $f$ (namely, $\Delta_{0,k}$) is the largest (up to
multiplication by a constant). This can be seen as follows. First
suppose that $f$ is also $S$-unimodal. If $n>0$ is such that $f^n(x)$
belongs to the interval with endpoints $-x,x$ but $f^j(x)$ does not,
for all $1\le j<n$, then $|Df^n(x)|>1$ -- this uses the fact that $f$
has negative Schwarzian. From this it follows that if $J\subset
[-1,1]$ is an interval that does not contain the critical point, whose
iterates $f^j(J)$ are pairwise disjoint for $0\le j\le n$, such that
$f^n(J)$ lies in the convex-hull of $J$ union its symmetric while the
previous iterates $f^j(J)$, $1\le j <n$, do not, then $|f^n(J)|>|J|$.
Hence, if $f$ is renormalizable, symmetric and $S$-unimodal then
at each renormalization level the interval that contains the critical
point is the largest. If we drop the negative Schwarzian hypothesis,
the same is true up to a multiplicative constant. This is because
every sufficiently deep renormalization of $f$ already has negative  
Schwarzian
derivative.

\begin{proposition}
\label{zaf}
For each $\alpha > 0$ there exist constants $C_0$ and $0 < \mu < 1$
such that
\begin{equation}
\label{estimate}
\sum_{i=0}^{p_k-1} \frac{|\Delta_{i,k}|^{2+\alpha}}{|\Delta_{i+1,k}|}
\le C_0 \mu^k \ .
\end{equation}
\end{proposition}

\Proof Let $\ell(\Delta_{i,k})$ be the level of $\Delta_{i,k}$, i.e.,
the largest integer $j$ such that $\Delta_{i,k}\subseteq  
\Delta_{0,j}\setminus
\Delta_{0,j+1}$. Let $d_{i,k}$ be the distance from $\Delta_{i,k}$ to  
zero
(the critical point). Using that $\Delta_{i,k}$ has space around itself
we see that for all $i\neq 0$ and all $x\in\Delta_{i,k}$ we have
$d_{i,k}\le |x|\le Kd_{i,k}$, where $K>1$ is a constant that depends  
only on
the real bounds. Hence $K^{-1}\le |x|/|y|\le K$ whenever
$x,y\in\Delta_{i,k}$. These facts are implicitly used in the estimates  
below.

Now, we have $|\Delta_{i,k}|/|\Delta_{i+1,k}|=1/|f'(x_{i,k})|$ for some
$x_{i,k}\in\Delta_{i,k}$. Since the critical point is quadratic, we have
$|f'(x_{i,k})|\ge C_1|x_{i,k}|$, and so
\begin{equation}
\label{manip}
\frac{|\Delta_{i,k}|}{|\Delta_{i+1,k}|}\leq \frac{1}{C_1|x_{i,k}|}
\ .
\end{equation}
Therefore, for all $0\le j\le k-1$,
\begin{equation*}
\begin{split}
\sum_{\ell(\Delta_{i,k})=j} \frac{|\Delta_{i,k}|^{2}}{|\Delta_{i+1,k}|}
& \;\le\;
C_1^{-1}\sum_{\ell(\Delta_{i,k})=j}\frac{|\Delta_{i,k}|}{|x_{i,k}|}\; 
\le\;
C_2\sum_{\ell(\Delta_{i,k})=j}\int_{\Delta_{i,k}}\frac{dx}{|x|}\\
&\;\le\;
C_3\int_{\Delta_{0,j}\setminus\Delta_{0,j+1}}\frac{dx}{|x|}\;\le\;
C_4\log\frac{|\Delta_{0,j}|}{|\Delta_{0,j+1}|}
\ .
\end{split}
\end{equation*}
With these estimates, and using the fact
proved above that $|\Delta_{0,k}|\ge
C_5|\Delta_{i,k}|$ for all $0\le i\le p_k-1$, we see that
\begin{equation*}
\begin{split}
\sum_{i=0}^{p_k-1}  
\frac{|\Delta_{i,k}|^{2+\alpha}}{|\Delta_{i+1,k}|}&\;\le\;
C_6\max{|\Delta_{i,k}|^{\alpha}}\left( 1+ \sum_{j=0}^{k-1}
\log\frac{|\Delta_{0,j}|}{|\Delta_{0,j+1}|}\right)\\
& \le C_6|\Delta_{0,k}|^{\alpha}
\left(1+\log{\frac{1}{|\Delta_{0,k}|}}\right)\\
& \le K_{\alpha}|\Delta_{0,k}|^{\alpha/2}
\ ,
\end{split}
\end{equation*}
where $K_{\alpha}$ is a positive constant depending on $\alpha$.
This proves {\eqref{estimate}} because $|\Delta_{0,k}|$ decays  
exponentially
with $k$,  with uniform rate depending only on the real bounds.
\Endproof

In addition to Proposition \ref{zaf} --- valid for maps with arbitrary
combinatorial type --- we shall need also an estimate that seems  
specific for
maps with {\it bounded} combinatorial type, namely Proposition \ref{zef}
below; first, a couple of lemmas.

For each $f\in {\mathbb{K}}$, let $d_f$ be the infimum of all positive  
numbers
$s$ such that
$$
\sum _{j=0}^{p_k-1}
|\Delta_{j,k}|^s \to 0 ~~ {\rm as} ~~ k \to \infty
\ .
$$
It is possible to prove, using some thermodynamic
formalism, that $d_f$ agrees with the Hausdorff dimension of
$\mathcal{I}_f$, but we will not need this fact.
Let $0<D<1$ be the supremum of $d_f$ as $f$ ranges through $\mathbb{K}$.

\begin{lemma}\label{hausdorff}
For each $s>D$ there exist $C_s>0$
and $0<\eta_s<1$
such that for all $f \in \mathbb{K}${\rm ,}
$$\sum _{j=0}^{p_k-1}
|\Delta_{j,k}|^s < C_s \eta_s^k \ .
$$
\end{lemma}

\Proof
Apply  bounded geometry and the  compactness of $\mathbb{K}$.
\Endproof

Next, let us define
$$
S_{j,k}(f;s) =
\sum_{\ell(\Delta_{i,k})=j}|\Delta_{i,k}|^s
\ ,
$$
for $j=0,1, \dots ,k-1$, all  $k\ge 0$, and all $f\in\mathbb{K}$.

\begin{lemma}\label{betomega}
For each $s>D$ and each $f \in \mathbb{K}$ we have $S_{j,k}(f;s)\le
C_s\lambda_j^s\eta_s^{k-j}${\rm ,} where $C_s>0$ and $0<\eta_s<1$ are the
constants of Lemma {\rm \ref{hausdorff},} and $\lambda_j=\lambda(f,j)$ is as  
in
\eqref{scalings}.
\end{lemma}

\Proof Using renormalization, we see that
$S_{j,k}(f;s)=\lambda_j ^{s} S_{0,k-j}(R^j(f);s)$.
 From Lemma \ref{hausdorff}, we know that
$$
S_{0,k-j}(R^j(f);s)\le  C_s \eta_s^{k-j} \ .
$$
The result follows.
\hfq

\begin{proposition} \label{zef}
For each $\mu>1$ close to one{\rm ,} there exist $0<\alpha <1-D$ close to  
zero and $C>0$
such that for all $f \in \mathbb{K}${\rm ,}
$$
\sum_{i=0}^{p_k-1} \frac{|\Delta_{i,k}|^{2-\alpha}}{|\Delta_{i+1,k}|}
\le C \mu^k \ .
$$
\end{proposition}

\Proof Using \eqref{manip} and the fact that $|x_{i,k}|\geq  
\lambda_{j+1}$
   when $\ell(\Delta_{i,k})=j$, we have
\begin{eqnarray*}
\sum_{i=0}^{p_k-1} \frac{|\Delta_{i,k}|^{2-\alpha}}{|\Delta_{i+1,k}|}
&=& \sum_{i=0}^{p_k-1} \frac{|\Delta_{i,k}|}{|\Delta_{i+1,k}|}\,
|\Delta_{i,k}|^{1-\alpha}\\
&\leq& C\left(|\Delta_{0,k}|^{-\alpha}+\sum_{j=0}^{k-1}
   \lambda_j^{-1}S_{j,k}(f;1-\alpha)\right)
\ .
\end{eqnarray*}
If $1-\alpha>D$ then, by Lemma \ref{betomega} with $s=1-\alpha$,  
$$
\sum_{i=0}^{p_k-1} \frac{|\Delta_{i,k}|^{2-\alpha}}{|\Delta_{i+1,k}|}
\leq CC_{1-\alpha}\sum_{j=0}^{k}\lambda_j^{-\alpha}\eta_{1-\alpha}^{k-j}
\leq K_\alpha\mu_{\alpha}^k
\ ,
$$
where $K_\alpha>0$ and
$\mu_\alpha=\max\{\lambda(f,1)^{-\alpha}:f\in\mathbb{K}\}$ depend on
$\alpha$. But if $\alpha$ is small enough we will have  
$\mu_\alpha<\mu$, and
this completes the proof.
\hfq

\begin{remark}
\label{515151}
By a continuity argument and the real bounds,
we can prove that Propositions \ref{zaf} and \ref{zef}
remain true for maps $\tilde{f} \in \mathbb{U}^{\,4}$
sufficiently close to $\mathbb{K}$.
More precisely, for each $k > 0$ there exists
$\varepsilon_k > 0$ such that for all $f \in \mathbb{K}$
and all $\tilde{f} \in \mathbb{U}^4$ with
$\left \|\tilde{f}-f \right \|_{C^4(I)} < \varepsilon_k$,
the map $\tilde{f}$ is $k$-times renormalizable,
and the statements of both propositions hold for $\tilde{f}$.
This will be used in \S \ref{sec:frechet} only for real analytic maps  
in an
open neighborhood of $\mathbb{K}$ in $\mathbb{A}$ .
\end{remark}

\vglue-12pt
\Subsec{Spectral estimates}
\label{sec:deriv}
In this section we prove  Theorem \ref{dfgdsfgg}.
Fixing $f\in\mathbb{K}$ and considering the Banach space $\mathcal{A}$
given by Theorem \ref{lyubhyp}, we recall that the Fr\'echet derivative
$L_f=DT(f):\mathcal{A}\to \mathcal{A}$ is given by formula
{\eqref{deetee}}. It is clear from that formula that $L_f$ extends to a
bounded linear operator $\hat{L}_f:\mathbb{A}^0\to \mathbb{A}^0$, and  
moreover
$\hat{L}_f(\mathbb{A}^s) \subseteq \mathbb{A}^s$ for all
$s\geq 0$ (because $f$ is analytic).

We want to verify the compatibility properties of
Definition~{\ref{robust}} for the spaces $\mathcal{B}=\mathbb{A}^s$
and $\mathcal{C}=\mathbb{A}^0$ when $s$ is close to (but different
from) $2$. Properties $\mathbf{A1}$ and  $\mathbf{A2}$ are clearly
satisfied, while $\mathbf{A3}$ follows from Lemma \ref{333}.
Property $\mathbf{A4}$ is a consequence of the following simple fact
about H\"older spaces (see \cite{Hor}).

\begin{lemma} There exists $\Delta>1$ such that $\mathbb{A}\cap
\mathbb{A}^s(\Delta)$ is $C^0$-dense in $\mathbb{A}^s(1)$.
\end{lemma}

\Proof By Theorem A.10 on page 43 of \cite{Hor}, there exist  a
family $S_t$, $t>0$, of smoothing operators preserving even
functions and $C\geq 1$ such that, for all $v \in \mathbb{A}^s(1)$,
we have $\|S_tv\|_{C^s} \leq C$ and $\|v- S_tv\|_{C^0}\leq C t^s$.
By the Stone-Weierstrass theorem, for all small $0<\varepsilon < C$
there is a polynomial $w_t$ with real coefficients and vanishing at
zero such that $\|w_t-S_tv\|_{C^s} < \varepsilon$. Now let
$v_t(x)=\frac{1}{2}(w_t(x)+w_t(-x))$, so that $v_t\in \mathbb{A}$
and we still have $\|v_t-S_tv\|_{C^s}<\varepsilon$. Then
$\|v_t\|_{C^s} < \|S_t(v)\|_{C^s} + \varepsilon< 2C$ on one hand,
while $\|v_t-v\|_{C^0} < \varepsilon + C t^s$ on the other hand. For
$t$ small enough, this gives $\|v_t-v\|_{C^0} < 2\varepsilon$ with
$v_t \in \mathbb{A}^s(2C)$. \phantom{overthere} \Endproof

Hence all that remains is to check that property $\mathbf{A5}$ is  
satisfied. By
Lemma \ref{spectral}, this will be the case provided we can control
the $C^0$ norms of $\hat{L}_{f,s}^{(m)}$.  We
shall prove this now, with the help of Propositions \ref{zaf} and  
\ref{zef}.

Recall that for each $m \ge 1$ the operator
$\hat{L}_f^{(m)}$ is an $L$-operator and its associated positive,
bounded linear operator $\hat{L}_{f,s}^{(m)}:\mathbb{A}^0 \to  
\mathbb{A}^0$
is given by
\begin{equation}
\label{E8}
\hat{L}_{f,s}^{(m)}(v)\;=\;
  \frac {1}{\lambda_k}
\sum_{j=0}^{p_k-1}|Df^j(f^{p_k-j}(\lambda_kx))|
|\lambda_k Df^{p_k-j-1}(\lambda_kx)|^s
v(f^{p_k-j-1}(\lambda_kx)) \ ,
\end{equation}
where $k=mN$ (recall that $T=R^N$).
Now we have the following fact  coming from bounded geometry:
\begin{equation}
\label{asd1}
|Df^j(f^{p_k-j}(\lambda_kx))|  \asymp
\frac{|\Delta_{0,k}|}
{|\Delta_{p_k-j,k}|} \ ,
\end{equation}
for all $0 \le j \le p_k-1$.
Since $|Df(\lambda_kx)|\le  C \lambda_k$
  for some   constant  $C>0$ independent of $k$ and
uniform in  $f \in {\mathbb{K}}$,
and $|\Delta_{0,k}| = 2\lambda_k$ we have
$$
| Df^{p_k-j-1}(\lambda_kx)| \le  C |\Delta_{0,k}|  
|Df^{p_k-j-2}(f(\lambda_kx))|
\ .
$$
Again, by bounded geometry, for all $0 \le j \le p_k-2$
$$
|Df^{p_k-j-2}(f(\lambda_kx))| \asymp
\frac{|\Delta_{p_k-j-1,k}|}
{|\Delta_{1,k}|} \ ,
$$
and so
\begin{equation}
\label{asd2}
|Df^{p_k-j-1}(\lambda_kx)|\le
C     |\Delta_{0,k}|
\frac{|\Delta_{p_k-j-1,k}|}
{|\Delta_{1,k}|} \ .
\end{equation}
Using \eqref{asd1}
and \eqref{asd2}
in \eqref{E8},
we see that
$$
\|\hat{L}_{f,s}^{(m)}\|\;\leq\;
\frac{C}{\lambda_k} \left(
  \frac{ \lambda_k^s  |\Delta_{0,k}| }{|\Delta_{1,k}|}
+
\sum _{j=0}^{p_k-2}
   \lambda_k^{2 s}
\frac{|\Delta_{0,k}|}{|\Delta_{p_k-j,k}|}
\frac{|\Delta_{p_k-j-1,k}|^s}{|\Delta_{1,k}|^s}
\right).$$
But $|\Delta_{0,k}|=2 \lambda_k $ and since the critical point
of $f$ is quadratic,
$|\Delta_{1,k} | \asymp |\Delta_{0,k}|^2 \asymp \lambda_k^2$.
Therefore, we arrive at
\begin{equation}
\label{E10}
\|\hat{L}_{f,s}^{(m)}\|\;\leq\;
C_1\sum _{j=0}^{p_k-1}
\frac{|\Delta_{p_k-j-1,k}|^s }{|\Delta_{p_k-j,k}|}
\ .
\end{equation}
The proof of part (i) of Theorem~\ref{dfgdsfgg} now follows from
Proposition~\ref{zaf}, while the proof of part (ii) is a consequence of
Proposition~\ref{zef}. This ends the proof of  Theorem \ref{dfgdsfgg}.

\section{The local stable manifold theorem}
\label{sec:stable}

In this section we isolate those features of the renormalization
operator that are essential for the promotion of ``hyperbolicity''
from the Banach space $\mathbb{A}$ of Theorem \ref{lyubhyp} to the
space  $\mathbb{U}^r$. This leads us to the definition of a robust
operator (see \S \ref{frfeww}). Such definition is necessarily
rather technical, since it has to account for the fact that the
renormalization operator is not Fr\'echet differentiable in
$\mathbb{U}^r$. In particular a robust operator acts simultaneously
on four different Banach spaces (corresponding in the case of
renormalization to the space $\mathbb{A}$ given by Theorem
\ref{lyubhyp}, $\mathbb{A}^r$, $\mathbb{A}^s$ and $\mathbb{A}^0$,
where $r >1+s$ and $s$ is close to $2$), and satisfies several
properties. The major goal of this section is to prove a local
stable manifold theorem for robust operators.

\Subsec{Robust operators}
\label{frfeww} Before moving on to a  precise definition of a robust
operator, we give the following informal description. A robust
operator acts on four Banach spaces $\mathcal{A} \subset \mathcal{B}
\subset \mathcal{C} \subset \mathcal{D}$. In the smaller space
$\mathcal{A}$ it acts smoothly and has a hyperbolic basic set
$\mathbb{K}$. The pairs of spaces $(\mathcal{B},\mathcal{D})$ and
$(\mathcal{C},\mathcal{D})$ are compatible with $(T,\mathbb{K})$,
and in particular the invariant hyperbolic splitting for
$\mathbb{K}$ in $\mathcal{A}$ extends to an invariant hyperbolic
splitting for $\mathbb{K}$ in $\mathcal{B}$. Viewed as a map from
$\mathcal{B}$ into $\mathcal{C}$, a robust operator is $C^1$. It
also satisfies a uniform Gateaux differentiability condition in
$\mathcal{C}$ for points and directions in $\mathcal{B}$. Finally,
as an operator in $\mathcal{B}$, it is reasonably well-approximated
by the extension of its derivative at a point of $\mathbb{K}$ in
$\mathcal{A}$ to a bounded linear operator in $\mathcal{B}$. It will
take us considerable effort (see \S \ref{sec:robust}) to verify that
the renormalization operator indeed satisfies all these conditions.

Let $T:\mathcal{O}\to\mathcal{A}$ be a $C^2$ operator having a compact
hyperbolic basic set $\mathbb{K}$ with the property that the unstable
subspace of the $DT$-invariant splitting of the tangent space at each  
point of
$\mathbb{K}$ is one-dimensional.
By standard invariant manifold theory (see~\cite{HP}), we know that
for all
$g\in \mathbb{K}$ the local unstable manifold $W^u(g)$ of
$T$ at $g$ exists and is $C^2$. In particular, we can
find a $C^2$ parametrization
$$
t\mapsto u_g(t)\in W^u(g)\subseteq \mathcal{A}
$$
varying continuously with $g$ such that
$\mathbf{u}_g=u_g'(0)$ is a unit vector.
We define a $C^2$ function
$t\mapsto \hat\delta_g(t)$ by
$$
T(u_g(t))\;=\;u_{T(g)}(\hat\delta_g(t)) \ .
$$
This function also
varies continuously with $g$ and
$\hat\delta_g(t)={\delta}_g t + O(t^2)$ for some
${\delta}_g>0$.
Recall that by hyperbolicity of $\mathbb{K}$, there exist $C_0>0$ and
$\lambda>1$ such that for every $g \in \mathbb{K}$ and every $m \ge 1$,
\begin{equation}
\label{fggtrs1}
\delta_{T^{m-1}(g)} \cdots \delta_g >  C_0 \lambda^m \ .
\end{equation}

\vglue8pt

\begin{definition}\label{def4}
Let $\mathcal{A}\subseteq\mathcal{B}\subseteq\mathcal{C}
\subseteq\mathcal{D}$ be Banach spaces, where each inclusion is a
compact linear operator. Let $\mathcal{O}_\mathcal{A}\subseteq
\mathcal{A}$ and $\mathcal{O}_\mathcal{B}\subseteq \mathcal{B}$ be
open sets in their respective spaces such that
$\mathcal{O}_\mathcal{A}\subseteq \mathcal{O}_\mathcal{B}$.
Let $\mathbb{K} \subset \mathcal{O}_\mathcal{A}$ be a hyperbolic
basic set of a $C^2$ operator $T:\mathcal{O}_\mathcal{A} \to
\mathcal{A}$. We  say that $T$ is {\it robust}  with respect to
$(\mathcal{B},\mathcal{C},\mathcal{D})$ if it has an extension to an
operator $  T:\mathcal{O}_\mathcal{B} \to \mathcal{B}$ that
satisfies the following conditions.
\begin{enumerate}
\item[$\mathbf{B1.}$] The pair $(\mathcal{B},\mathcal{D})$ is  
$1$-{\it compatible}
with $T$, while the pair $(\mathcal{C},\mathcal{D})$ is\break
$\rho_\mathcal{C}$-{\it compatible} with $(T,\mathbb{K})$ for some
$\rho_\mathcal{C} < \lambda$ (where $\lambda$ is as in
\eqref{fggtrs1}).
\item[$\mathbf{B2.}$]
For each $m>0$, the interior  $\mathcal{O}_\mathcal{B}^{(m)}$ of the
set $\{ f \in \mathcal{O}_{\mathcal{B}}:   T^i(f) \in
\mathcal{O}_{\mathcal{B}},\break \forall i <m\}$ contains $\mathbb{K}$,
and  $ T^m:\mathcal{O}_\mathcal{B}^{(m)}\to \mathcal{C}$ is $C^1$
and its derivative is uniformly continuous in some neighbourhood of
$\mathbb{K}$. Furthermore, for all $f \in \mathcal{A} \cap
\mathcal{O}_\mathcal{B}^{(m)}$ the linear map
$$
D  T^m(f):\mathcal{B}\to \mathcal{C}
$$
extends to a continuous linear operator $L_m:\mathcal{D}\to \mathcal{D}$
that satisfies $L_m(\mathcal{X}) \subseteq \mathcal{X}$,
  for  $\mathcal{X}=
\mathcal {B}, \mathcal{C}$.
\item[$\mathbf{B3.}$]
  For every $m$ there exists $C_{m,1}>1$
  with the property that for each  $g\in\mathbb{K}$ there is an open
set
  $\mathcal{V}_g\subseteq \mathcal{O}_\mathcal{B}$ containing $g$ such
that
for all $f \in \mathcal{V}_g$,
$$
\|DT^m(f)\mathbf{u}_g- DT^m(g)\mathbf{u}_g \|_{\mathcal{C}}\leq
  C_{m,1}  \|f-g\|_{\mathcal{B}}
\ . $$
\item[$\mathbf{B4.}$]
  There exist  $C_1>1$ and $\rho >1$
  with the property that for each  $g\in\mathbb{K}$ there is an open
set
  $\mathcal{V}_g\subseteq \mathcal{O}_\mathcal{B}$ containing $g$ such
that
for all $f_1, f_2\in \mathcal{V}_g$,
$$
\|T(f_1)-T(f_2)-DT(f_2)(f_1-f_2)\|_{\mathcal{C}}\leq
  C_1  \|f_1-f_2\|_{\mathcal{C}} ^{\rho}
\ .
$$
\item[$\mathbf{B5.}$]
  For all $m>0$, there exists $C_{m,2} > 0$, and there exists $\nu_m > 0$
such that for all $g \in {\mathbb K}$ and
for all $f\in \mathcal{B}$ with
$\|f-g\|_{\mathcal{B}} < \nu_m$,
$$
\|DT^m(f)-DT^m(g)\|_{\mathcal{C}}\leq
C_{m,2}  \lambda^m
\ .
$$
Moreover, there exists $m_0>0$ such that for all
$m > m_0$, $C_{m,2} <  C_0/8$
(where $C_0$ and $\lambda$ are as in \eqref{fggtrs1}).
\item[$\mathbf{B6.}$]
  For all $m>0$, there exists $C_{m,3} > 0$
such that   for all $g \in {\mathbb K}$, for all $f \in \mathcal{A}$
with $\|f-g\|_{\mathcal{A}} < \nu_m$ and
for all $v\in \mathcal{B}$ with
$\|v\|_{\mathcal{B}} < \nu_m$,  
$$
\|T^m(f+v)-T^m(f)-DT^m(g)v\|_{\mathcal{B}}\leq
C_{m,3} \|v\|_{\mathcal{B}}
\ .
$$
Moreover, there exists $m_0>0$ such that
for all $m > m_0$, $C_{m,3} < 1/4\/$.
\end{enumerate}
\end{definition}

\begin{example}
\label{hmgjghsd} As one might expect, the main example of a robust
operator is provided by renormalization. We know from Theorem
\ref{lyubhyp} that the renormalization operator $T={R}^N:\mathbb{O}
\to \mathbb{A}$ is hyperbolic over $\mathbb{K}$. We also know that
this renormalization operator is well-defined as  a map  from an
open set of $\mathbb{U}^\gamma$ containing $\mathbb{K}$ into
$\mathbb{U}^\gamma \subset 1+\mathbb{A}^\gamma\cong
\mathbb{A}^\gamma$  for all $\gamma \ge 2$. We will show in Section
\ref{sec:robust} that $T$ is robust with respect to
  the spaces ${\mathcal A}=\mathbb{A}$,
${\mathcal B}=\mathbb{A}^r $, ${\mathcal C}=\mathbb{A}^s$
and ${\mathcal D}=\mathbb{A}^0$ whenever $s<2$ is close to $2$ and
$r>s+1$ is not an integer.
\end{example}

\vglue-12pt

\Subsec{Stable manifolds for robust operators}
  We can now formulate a general local stable manifold theorem for
robust operators.

\begin{theorem}
\label{main3}
Let $T:\mathcal{O}_\mathcal{A} \to \mathcal{A}$ be a
$C^k$ with $k\ge 2$ \/{\rm (}\/or real analytic\/{\rm )}\/  hyperbolic
operator over $\mathbb{K} \subset \mathcal{O}_\mathcal{A}${\rm ,}  and
be robust
with respect to $(\mathcal{B},\mathcal{C},\mathcal{D})$.
Then conditions {\rm (i), (ii), (iii)} and {\rm (iv)} of Theorem~{\rm \ref{hpicture}}
hold true for the   operator $T$ acting on $\mathcal{B}$. The local
unstable manifolds are $C^k$ with $k\ge 2$ \/{\rm (}\/or real analytic\/{\rm )}\/
  curves{\rm ,} and the local stable manifolds are
of class $C^1$ and form a $C^0$ lamination.
  \end{theorem}

The proof of this theorem will occupy the rest of
Section~\ref{sec:stable}. In the end, the theorem will follow by
  putting together Corollary \ref{corrob},
Proposition \ref{dfgfdger} and Theorem \ref{sddw}.

\Subsec{Uniform bounds}
Before proceeding  we prove the following
simple bounds that we will use quite often.

\begin{lemma}
\label{l11100}
There exist $\mu_0>0$ and  $1  < \lambda < M$
such that for all $g \in  \mathbb{K}$ and all $t
\in  \mathbb{R}$ with $|t| < \mu_0${\rm ,}
   $u_g(t)$ and $\hat \delta_g(t)$ are well-defined and
\begin{enumerate}
\item
$M^{-1} \lambda^n    <  \delta_{T^{n-1}(g)} \cdots \delta_g     <
M^n$ and
$\left | \hat{\delta}_g(t) \right  |    <     M |t|${\rm ;}
\item
$M^{-1} < \|\mathbf{u}_g\|_{\mathcal B}   <    M $ and
$M^{-1} < \|\mathbf{u}_g\|_{\mathcal C}   <     M ${\rm ;}
\item
$ M^{-1} <    \left \|\sigma_g \right  \|_{\mathcal B} < M $ and
$ M^{-1} <    \left \|\sigma_g \right  \|_{\mathcal C} <M ${\rm ;}
\item
$M^{-1} |t| < \left \|u_g(t)-g \right  \|_{\mathcal C}    <
M |t|$ and
$M^{-1} |t|  <  \left \|u_g(t)-g \right  \|_{\mathcal B}     <   M |t|${\rm ;}
\item
$\left | \sigma_g \left( u_g(t) - g \right ) \right  |
> \frac{1}{2} |t|   $.
\end{enumerate}

\end{lemma}

\Proof
(i) By Definition \ref{hyp} and \eqref{fggtrs1}, there exist  $1  <  
\lambda <
M_1$ such that for all $g \in \mathbb{K}$ and all $n \ge 1$ we have
$M^{-1}_1 \lambda^n < \delta_{T^{n-1}(g)} \cdots \delta_g  < M^n_1$
and also
$\left | \hat{\delta}_g(t) \right  |     <     M_1 |t|$
for  all $|t| \le  \mu_1$
(where $\mu_1>0$ is a uniform constant).
\vskip6pt

(ii) For ${\mathcal X}$ equal to  ${\mathcal B}$ and ${\mathcal C}$,
we have that   $g \mapsto \mathbf{u}_g$ as a map
${\mathbb K} \to {\mathcal X}$ is continuous and does not vanish. Hence,
by   compactness of $\mathbb{K}$ there is $M_2 > 1$ such that
$M_2^{-1} < \|\mathbf{u}_g\|_{\mathcal X}   <     M_2$.

\smallbreak (iii) Since $\sigma_g(\mathbf{u}_g)=1$ and by property $\mathbf{B1}$
in Definition \ref{def4}, the functional~$\sigma_g$ extends
continuously to ${\mathcal X}$ and there is $M_3 > 1$ such that
      $M_3^{-1} <     \left \|\sigma_g \right  \|_{\mathcal X}    < 
M_3$.

\smallbreak (iv) Since $t \to u_g (t)$ as a map $\mathbb{R} \to {\mathcal X}$
is   $C^1$ and varies continuously with $g \in {\mathbb K}$, there is
$M_4 > 0$ and $\mu_2 > 0$ such that
$
\|u_g  (t) -
u_g  (s)
\|_{{\mathcal X}}
  \le M_4 |t-s|
$ for all $g \in {\mathbb K}$ and all $t,s$ with $|t| < \mu_2$ and $|s|  
<
\mu_2$.
Moreover, since
$$
\left.
\frac{d}{dt}
u_g  (t)\right | _{t=0}  =
\mathbf{u}_g \ne 0 \,
$$
there exists $M_5 > 0$ and $\mu_3 > 0$ such that
$ |t-s|\le M_5
\|u_g  (t) -
u_g  (s)
\|_{{\mathcal X}}$
for all $g \in {\mathbb K}$ and all $|t| < \mu_3$. Hence (iv) follows by
taking $s=0$ and noting that $u_g(0)=g$.

\smallbreak(v) This follows from (iv) and the fact that $\sigma_g(\mathbf{u}_g)=1$.
\hfq

\def\reals{\mathbb{R}}
\Subsec{Contraction towards the unstable manifolds}
The one-dimensional unstable manifolds of $T$ in $\mathcal{A}$
are embedded in $\mathcal{B}$, and remain invariant.
The first important estimate given by the following lemma shows that
in $\mathcal{B}$ the operator $T$ contracts towards such manifolds.
Therefore, if $T$ is to have unstable manifolds in  $\mathcal{B}$,
these have to coincide with unstable manifolds in $\mathcal{A}$.
In what follows, we fix $g \in \mathbb{K}$ and for simplicity of
  notation we write
$$
\sigma_i=\sigma_{T^i(g)}\ , \
u_i=u_{T^i(g)}\ ,\
{\mathbf u}_i= {\mathbf u}_{T^i(g)}\ ,
$$
and
$$
\delta_i^m = \prod_{j=i}^{m-1}  \delta_{T^{j}(g)}, \
  \hat{\delta}_i^m =\hat{\delta}_{T^{m-1}(g)} \circ  \cdots \circ
\hat{\delta}_{T^{i}(g)} \ .
$$
Set $\mu_0>0$ as in Lemma \ref{l11100}.

\begin{lemma}
\label{L21}
For every $m>0$ there exist
  $0 < \eta_m < \mu_0$ and $B_m> 0$ such that for every $g \in
\mathbb{K}$ and  every $v \in {\mathcal B}$
with  $\|v\|_{{\mathcal B}}<\eta_m$ and $t \in \reals$ with
$|t|<\eta_m${\rm ,}   $\hat{\delta}_0^m (t+\sigma_0(v)) < \mu_0${\rm ,}
   $u_0(t)+v \in {\mathcal O}_{\mathcal B}^{(m)}$ and
$$
\left \|T^m(u_0(t)+v)-u_m \left (\hat{\delta}_0^m (t+\sigma_0(v))
\right )
\right \|_{{\mathcal B}}
<    B_m \|v\|_{{\mathcal B}} \ .
$$
  Furthermore{\rm ,} there is $m_1>m_0$ \/{\rm (}\/where $m_0$ is as in $\mathbf{B6}${\rm )}  
such
  that for all $m > m_1${\rm ,} $B_m < 1/2$.
\end{lemma}

\Proof We prove the second inequality only. The first is proven in
the same way. By  property  $\mathbf{B6}$  in Definition \ref{def4},
there is $m_0 > 0$ such that  for all $m> m_0$, all $v$ with
$\|v\|_{\mathcal{B}} < \nu_m$, and all $t \in \mathbb{R}$ with $|t|
< \nu_m$,
\begin{equation}
\label{www1}
\left \|T^m(u_0(t)+v)- T^m(u_0(t))-DT^m(g) v \right \|_{\mathcal{B}}
\le \frac{1}{4} \|v \|_{\mathcal{B}} \ .
\end{equation}
By property  $\mathbf{B1}$ in Definition \ref{def4},
$(\mathcal{B},\mathcal{D})$ is 1-compatible with $(T,\mathbb K)$.
Hence, by Corollary \ref{corrob}, there exists $m_1> m_0$ such that
for all $m > m_1$ we have
\begin{equation}
\label{www2}
\left \| DT^m (g) v - \delta_0^m  \sigma_0 (v) \mathbf{u}_m
  \right \|_{\mathcal{B}}
   \le \frac{1}{8} \|v \|_{\mathcal{B}} \ .
\end{equation}
Putting \eqref{www1} and \eqref{www2} together we get
\begin{equation}
\label{www3}
\left \|T^m(u_0(t)+v)- T^m(u_0(t))-\delta_0^m  \sigma_0 (v)
\mathbf{u}_m \right \|_{\mathcal{B}}
  \le \frac{3}{8} \|v \|_{\mathcal{B}} \ .
\end{equation}
Now, we know that
$T^m(u_0(t))=u_m (\hat{\delta}_0^m(t))$
and $t \to u_m \circ \hat{\delta}_0^m(t)$  is $C^2$.
Hence,
\begin{eqnarray}
\label{www4}
\left \|u_m \circ \hat{\delta}_0^m (t+\sigma_0(v))- u_m \circ
\hat{\delta}_0^m (t)
\right . &-& \left .\delta_0^m  \sigma_0(v) \mathbf{u}_m  \right  
\|_{\mathcal{B}} \nonumber\\
& \le & c_1 \left( (\sigma_0(v))^2+|t|\sigma_0(v) \right)  \nonumber\\
& \le & c_2 \left(\|v \|_{\mathcal{B}}+|t|  \right) \|v\|_{\mathcal{B}}  
\ .
\end{eqnarray}
Therefore, choosing  $\eta_m < \nu_m$ so small that
$C_2 \eta_m < 1/16$ and putting \ref{www3} and \ref{www4} together,
we see that if  $|t| < \eta_m$ and $\|v\|_{\mathcal{B}} < \eta_m$ then
$$
\left \|T^m(u_0(t)+v)-u_m \circ \hat{\delta}_0^m (t+\sigma_0(v))
\right \|_{\mathcal{B}}
  \le \frac{1}{2} \|v \|_{\mathcal{B}}
$$
as desired.
\hfq

\begin{lemma}
\label{l100}
Let $m_1> 0$ be as in Lemma {\rm \ref{L21}.}
For all  $m > m_1$
  there exist small constants $0 <  \varepsilon_2 < \varepsilon_1 <  
\varepsilon_0$
such that the following holds
for every $\varepsilon < \varepsilon_2$.
For  every $g \in \mathbb{K}$ and  every  $v \in {{\mathcal B}}$ with
   $\|v\|_{{\mathcal B}}<\varepsilon${\rm ,}
the recursive scheme given by $f_0=g+v${\rm ,} $t_0=0${\rm ,} $v_0=v$ and
\begin{eqnarray}
\label{sadsa}
f_{k+1} & = & T^m(f_k),  \\
t_{k+1} & = & \hat{\delta}_{km}^{(k+1)m} (t_k + \sigma_{km}(v_k) ),
\nonumber
\\
v_{k+1} & = & f_{k+1} -  u_{(k+1)m}(t_{k+1})\nonumber
\end{eqnarray}
is well-defined  for all
$k=0, \dots,k_0  -1$
where $$k_0 =k_0(g,f_0) =  \min \{j:|t_j| \ge \varepsilon_1 \}.
$$
For all $k \le   k_0${\rm ,}
\begin{eqnarray}
\label{fgddfg}
\left \|T^{km}(g+v) - u_{km}(t_k) \right \|_{\mathcal{B}}& <&
2^{-k}\|v\|_{{\mathcal B}} ~~\mbox{and}~~\\
\left \|T^{km}(g+v) - T^{km}(g)\right  \|_{\mathcal{B}}& <&
\varepsilon_0  /M_0 \ , \nonumber
\end{eqnarray}
where $M_0= M^{m+2}+ B_1+\dots+B_m${\rm ,} $M$ is as in Lemma {\rm \ref{l11100}}
and $B_1,\dots,B_m$ are as in Lemma {\rm \ref{L21}.}
Furthermore{\rm ,}
\begin{enumerate}
\item  $\varepsilon_1 \le  \left |t_{k_0} \right |<\varepsilon_0${\rm ;}
\item   $\left \|T^{k_0m}(g+v) - T^{k_0m}(g)\right \|_{\mathcal{B}} >
\varepsilon_2${\rm ;}
\item     $\left |\sigma_{k_0m}\left (T^{k_0m}(g+v) - T^{k_0m}(g)
\right )\right | > \varepsilon_2${\rm ;}
\item  $\left \|T^{km+i}(g+v) - T^{km+i}(g)\right  \|_{\mathcal{B}} <
M_0 \left \|T^{km}(g+v) - T^{km}(g)\right  \|_{\mathcal{B}}$
\/{\rm (}\/which is less than $\varepsilon_0${\rm )}
for all $k \le   k_0$ and all $i=0,\dots,m$.
  \end{enumerate}
\end{lemma}

\Proof
For every $g \in \mathbb{K}$, let ${\mathcal B}(g,\varepsilon)$ be the  
open
ball in $\mathcal{B}$
of radius $\varepsilon$ centered at $g$.
Let us fix $m> m_1$ and
choose $\varepsilon_0 < \min\{\mu_0,\eta_1, \dots , \eta_m \}$
  such that all the properties $\mathbf{B1}$ to $\mathbf{B6}$ of  
Definition \ref{def4}
are satisfied in $\bigcup_{g \in \mathbb{K}} {\mathcal  
B}(g,\varepsilon_0)
\subset {\mathcal O}_{\mathcal B}^{(m)}$,
  where
$\mu_0$ is as in Lemma \ref{l11100} and  $\eta_1, \dots ,\eta_m$ are
as  in
Lemma \ref{L21}. Since $m > m_1$,
we have that $B_m < 1/2$ where $B_m$ is as in Lemma \ref{L21}.
Let us take $M>1$   as in Lemma \ref{l11100}.
We choose $0 <\varepsilon_2 < \varepsilon_1 < \varepsilon_0 < \mu$
such that
\begin{eqnarray}
\label{0101}
   \varepsilon_1 & < &  \varepsilon_0 / \left( 3M_0 M^{m+2} \right) \ ,  
 \\
   \varepsilon_2   & < & \varepsilon_1  /(2 + 2M)   \ .\nonumber
\end{eqnarray}
Now  we work by induction on $k$. Let us assume that $f_k$, $t_k$,
and $v_k$ have
been defined so that \eqref{fgddfg} holds.
Hence  $ \left \|f_k -  T^{mk}(g) \right \|_{\mathcal{B}}
   \le   \varepsilon_0 \le \mu $,
and so $f_k \in  {\mathcal O}_{\mathcal B}^{(m)}$ and
$f_{k+1}=T^m(f_k)$ is well-defined.
  Since $|t_k|\le \varepsilon_1$ and $2 M^{m+1} \varepsilon_1  \le
\varepsilon_0 \le \mu$, by Lemma \ref{l11100} and
   \eqref{l11100}, and by \eqref{sadsa}    and   \eqref{0101},
  we have  that $t_ {k+1}$ is well-defined and
\begin{eqnarray}
\label{dddDDd}
\left |t_{k+1} \right | & = &
\left | \hat{\delta}_{km}^{(k+1)m} (t_k + \sigma_{km}(v_k) )  \right
|
   \le    M^m (\left | t_k  \right | +
  \left | \sigma_{km}(v_k)  \right | )  \\
  & \le &   M^m \left( \varepsilon_1 + M \frac{\varepsilon}{2^k}
\right )
   <    2 M^{m+1} \varepsilon_1
    <    \varepsilon_0  \ .\nonumber
\end{eqnarray}
   Thus, by Lemma \ref{l11100},
$u_{(k+1)m}(t_{k+1})$ and
   $v_{k+1}=f_{k+1}-u_{(k+1)m}(t_{k+1})$ are also well-defined.
   By Lemma \ref{L21} and by \eqref{sadsa}, we get
\begin{eqnarray}
\label{eerrewq}
 \qquad\quad \left \|v_{k+1} \right  \|_{\mathcal{B}}
& =  &
\left \|T^m(f_k) -u_{(k+1)m}(t_{k+1}) \right \|_{\mathcal{B}}
 \\
  & =  &
\left \|T^m(v_k+u_{km}(t_k)) -u_{(k+1)m}
\left (\hat{\delta}_{km}^{(k+1)m}
(t_k + \sigma_{km}(v_k) ) \right) \right \|_{\mathcal{B}} \nonumber
\\
& \le &   2^{-1}   \left  \|v_k \right \|_{\mathcal{B}}
  \le  2^{k+1}  \left \|v_0 \right \|_{\mathcal{B}}    \ .\nonumber
\end{eqnarray}
Now, let us estimate $\left \|f_{k+1} -  T^{(k+1)m}(g) \right  
\|_{\mathcal{B}}$.
 From \eqref{sadsa} and  \eqref{eerrewq}, we get
\begin{equation}
\label{fddfdfdW}
\left \|f_{k+1} -  u_{(k+1)m}(t_{k+1}) \right \|_{\mathcal{B}}
\le \left \|v_{k+1} \right \|_{\mathcal{B}} \le \frac
{\varepsilon}{2^k} \ .
\end{equation}
 From Lemma \ref{l11100} and by \eqref{dddDDd}, we obtain
\begin{equation}
\label{fddfdsfdW}
\left \|u_{(k+1)m}(t_{k+1}) - T^{(k+1)m}(g) \right \|_{\mathcal{B}}
\le M \left |t_{k+1} \right | \le
2 M^{m+2} \varepsilon_1 \ .
\end{equation}
Thus, by \eqref{0101}, \eqref{fddfdfdW} and \eqref{fddfdsfdW},
\begin{eqnarray*}
\left \|f_{k+1} -  T^{(k+1)m}(g) \right \|_{\mathcal{B}}
& \le &
\left \|f_{k+1} -  u_{(k+1)m}(t_{k+1}) \right \|_{\mathcal{B}} \\
& & + \left \|u_{(k+1)m}(t_{k+1}) - T^{(k+1)m}(g) \right
\|_{\mathcal{B}} \\
& \le &  \frac {\varepsilon}{2^k}  + 2 M^{m+2} \varepsilon_1
\\
& \le &  3M^{m+2} \varepsilon_1 < \varepsilon_0 / M_0 \ .
\end{eqnarray*}
  This completes the induction.

Now, we must prove (i), (ii), (iii) and    (iv).
  Property (i) follows from  \eqref{dddDDd}.
Let us  prove (ii).
By property (i)  and  Lemma \ref{l11100},
$$ \left \| u_{k_0 m} (t_{k_0} ) -  T^{k_0 m}(g) \right
\|_{\mathcal{B}} \ge M^{-1} \varepsilon_1 \ .$$
By  \eqref{eerrewq}, we get
  $\left \| v_{k_0} \right  \|_{\mathcal{B}}   \le
  \varepsilon / 2^{k_0}$.
Thus, by  \eqref{0101}, we obtain
\begin{eqnarray*}
\left \|T^{k_0 m}(g+v) - T^{k_0 m}(g) \right  \|_{\mathcal{B}}
&  =  &
\left \|u_{k_0 m} (t_{k_0} ) + v_{k_0}  -  T^{k_0 m}(g)  \right
\|_{\mathcal{B}}  \\
  &  \ge  &
\left |
\left \| u_{k_0 m} (t_{k_0} ) -  T^{k_0 m}(g) \right
\|_{\mathcal{B}} -
\left \| v_{k_0} \right  \|_{\mathcal{B}} \right  |\\
  &  \ge  &
  M^{-1} \varepsilon_1 -  \frac{\varepsilon}{2^{k_0}} \\
  & \ge  &
  \varepsilon_2
  \ .
\end{eqnarray*}
Let us prove (iii).
By  property (i)  and Lemma \ref{l11100},
$$\left |\sigma_{k_0 m} (u_{k_0 m} (t_{k_0} - T^{k_0 m} g) ) \right |
\ge
\varepsilon_1 / 2 \ .$$
Using Lemma \ref{l11100} yet again and  \eqref{eerrewq},
we have
$\left |\sigma_{k_0 m} (v_{k_0}) \right  | \le
M    \varepsilon/ 2^{k_0}$. Thus, by \eqref{0101},
\begin{eqnarray*}
\left |\sigma_{k_0 m} \left (T^{k_0 m} (g+v) - T^{k_0 m} g\right)
\right |
&  \ge  &
\left |\  \left |\sigma_{k_0 m} (u_{k_0 m} (t_{k_0} ) ) \right | -
\left |\sigma_{k_0 m} (v_{k_0}) \right  | \ \right  | \\
& \ge  &
\frac{1}{2} \varepsilon_1  - M   \frac{\varepsilon_2}{2^{k_0}} \\
& \ge  &
  \varepsilon_2
  \ .
\end{eqnarray*}
Finally, we  prove (iv). Fix $0 \le k \le k_0$
and  $0 \le i \le m$. Setting $w_k=T^{km}(f)-T^{km}(g)$
  we have by \eqref{fgddfg} that
   $\|w_k\|_{\mathcal B} < \varepsilon_0 / M_0 < \eta_i$
where $\eta_i$ is as in Lemma  \ref{L21}.
Hence $T^{km}(g)+w_k \in {\mathcal O}_{\mathcal B}^{(i)}$
and by  Lemma  \ref{L21},
$$
\left \|T^i (T^{km}(g)+w_k)-u_{km+i}\left
(\hat{\delta}_{km}^{km+i}(\sigma_{km}(w_k)) \right )
\right \|_{\mathcal B} \le B_i \|w_k\|_{\mathcal B} \ .
$$
On the other hand, by Lemma \ref{l11100},  we have
$$
\left \|u_{km+i} \left
(\hat{\delta}_{km}^{km+i}(\sigma_{km}(w_k))\right )-T^{km+i}(g)
\right \|_{\mathcal B}
\le
M^{m+2} \|w_k\|_{\mathcal B} \ .
$$
Therefore,
\begin{eqnarray*}
  \|T^{km+i}(f)-T^{km+i}(g)\|_{\mathcal B}
& \le &
\left \|T^i (T^{km}(g)+w_k)-u_{km+i}\left
(\hat{\delta}_{km}^{km+i}(\sigma_{km}(w_k)) \right )
\right \|_{\mathcal B} \\
& & +
\|u_j \left  (\hat{\delta}_{km}^{km+i}(\sigma_{km}(w_k))\right )-
T^{km+i}(g)\|_{\mathcal B} \\
& \le &  (B_i +  M^{m+2}) \|w_k\|_{\mathcal B}
     \le  \varepsilon_0 \ ,
\end{eqnarray*}
which ends the proof.
\hfq

\Subsec{Local stable sets}
Let us now consider the local stable set
$W^s_\varepsilon(g)$  of $T$ at $g$ in $\mathcal B$ which
consists of all points $f \in  \mathcal{B} (g, \varepsilon)$
such that  for all $n >  0$, we have
$T^n (f)  \in
\mathcal{B}  (T^{n}(g),\varepsilon)$
and
$$
\left \|T^n (f)- T^n (g) \right \|_\mathcal{B}
  \to 0 ~{\rm when}~ n \to \infty \ .
$$
Our aim in this section is to give a finite characterization of
  $W^s_\varepsilon(g)$   and prove that $T$ contracts in the
$\mathcal{B}$-norm exponentially along $W^s_\varepsilon(g)$.
This is done in Lemma \ref{L22} below (see also Remark \ref{gggrad}).

 From now on in this section, we let $m_1$ and
$\varepsilon_0> \varepsilon_1> \varepsilon_2 > \varepsilon$
be as in
Lemma \ref{l100}.
For all sufficiently small  $0<\varepsilon < \varepsilon_2$ and
for all $f \in {\mathcal B}(g,\varepsilon)$,
we let $k_0=k_0(g,f)$ and $t_k=t_k(g,f)$ for  $k=0,\dots,k_0$ be as  in
Lemma \ref{l100}.
We write
${\mathcal B}(g,\varepsilon)=V_{\varepsilon}^-(g) \cup  
V_{\varepsilon}^0(g)
\cup V_{\varepsilon}^+(g)$ where
\begin{eqnarray*}
V_{\varepsilon}^-(g)
& = &
\left \{f\in {\mathcal B}(g,\varepsilon): -\varepsilon_0 <t_{k_0}(g,f)<-
\varepsilon_1 \right \} \ , \\
V_{\varepsilon}^+(g)
& = &
\left \{f \in {\mathcal B}(g,\varepsilon): \varepsilon_1 <t_{k_0}(g,f)<
\varepsilon_0 \right  \} \ , \\
V_{\varepsilon}^0(g)
& = &
{\mathcal B}(g,\varepsilon) \setminus \left ( V_{\varepsilon}^-(g) \cup
V_{\varepsilon}^+(g) \right ) \ .
\end{eqnarray*}

\begin{lemma}
\label{L22}
There exist an integer $m$ and a positive constant $C_2$ with
the following properties.
For all $\varepsilon > 0$   sufficiently small  and for all $g \in
\mathbb{K}${\rm ,} the sets $V_\varepsilon^-(g)$
and  $V_\varepsilon^+(g)$ are open subsets of ${\mathcal  
B}(g,\varepsilon)$
\/{\rm (}\/and so
$V_\varepsilon^0(g)$ is  relatively closed in ${\mathcal  
B}(g,\varepsilon)${\rm ),}
and  for all $f \in V_\varepsilon^0(g)$
\begin{equation}
\label{domingo1}
\left \|T^j(f)-T^j(g) \right\|_{\mathcal B}
\le \varepsilon C_2 2^{- j/m} \ .
\end{equation}
Furthermore{\rm ,} the local stable set $W^s_\varepsilon(g)$ is a
relatively open
subset of
$V_\varepsilon^0(g)$ and
\begin{equation}
\label{domingo2}
  W^s_\varepsilon(g) =
\left  \{f \in V_\varepsilon^0(g) :
\left  \|T^j(f)-T^j(g) \right \|_{\mathcal B} < \varepsilon, ~ {\rm
for} ~ {\rm all} ~0 \le j \le m \log C_2 / \log 2 \right \} \ .
\end{equation}
\end{lemma}

\Proof
The first assertion is a consequence of the definitions of
$V_\varepsilon^-(g)$ and $V_\varepsilon^+(g)$ and Lemma \ref{l100}.
It follows from property (i) of Lemma \ref{l100}  that
$$
V_\varepsilon^0(g)  =
\{f \in {\mathcal B}(g,\varepsilon)  : |t_k(g,f)|<\varepsilon_1 , ~
{\text{for all}} ~ k \ge  0 \} \ .
$$
It also follows from property (ii) of Lemma \ref{l100} that  if $f
\in  {\mathcal B}(g,\varepsilon)$
and $\left |t_{k_0}(g,f) \right |\break \ge \varepsilon_1$  then
$\left \|T^{k_0m}f-T^{k_0m}g \right  \|_{{\mathcal B}} >
\varepsilon$ where $k_0=k_0(g,f)$. This shows that
$W^s_\varepsilon(g) \subset V^0_\varepsilon(g)$, and therefore
\eqref{domingo1} implies \eqref{domingo2}. Furthermore,
   $W^s_\varepsilon(g)$ is a relatively open subset of
$V_\varepsilon^0(g)$.

It remains to show that if $f \in V_\varepsilon^0(g)$ then
\eqref{domingo1}
holds.
Set $1<\lambda< M$ as in Lemma  \ref{l11100}.
Since  $\beta > 2$, by Lemma  \ref{l11100}, there is  $m$
large enough such that
$ {\delta}_{km}^{(k+1)m} \ge M^{-1} \lambda^m >  \beta >2$
for every $k \ge 0$.
By Lemma \ref{l100}, for all $k \ge 0$, we know that $t_k=t_k(g,f)$ and
$v_k=v_k(g,f)$
are well-defined, and satisfy $|t_k| < \varepsilon_1$ and $\|v_k
\|_{\mathcal B} \le \varepsilon 2^{-k}$.
Furthermore,
$
t_{k+1} =\hat{\delta}_{km}^{(k+1)m} \left(t_k + \sigma_{km}(v_k)
\right)
$.
Since ${\delta}_{km}^{(k+1)m}$ is $C^2$
   and  $\|\sigma_{km}\|_{\mathcal B} < M$
(see Lemma \ref{l11100}), there is $c_0 > 1$ so
that
\begin{eqnarray*}
\label{domingo3}
\left |t_{k+1}-{\delta}_{km}^{(k+1)m} t_k \right |
& \le &
\left  |\hat{\delta}_{km}^{(k+1)m} \left(t_k + \sigma_{km}(v_k)
\right) -
  \hat{\delta}_{km}^{(k+1)m}  \left(t_k \right) \right  |     \\
  & &
+
\left  | \hat{\delta}_{km}^{(k+1)m} \left(t_k \right) -
{\delta}_{km}^{(k+1)m}  t_k \right |   \\
& \le &
c_0 \left (\|v_k\|_{\mathcal B} + |t_k|^2 \right ) \\
& \le &
c_0 \left (\varepsilon 2^{-k}+ |t_k|^2 \right ) \ .
\end{eqnarray*}
Hence \eqref{domingo3} gives
$$
|t_k| \le  \varepsilon c_0  \beta^{-1}2^{-k} + \beta^{-1} |t_{k+1}| +
c_0
\beta^{-1} |t_k|^2 \ .
$$
Taking $\varepsilon$ (in Lemma \ref{l100}) so small that $c_0
\beta^{-1} \varepsilon  < 1/2$,
   we get
\begin{equation}
\label{domingo4}
|t_k| \le
2 \left( |t_k|- c_0 \beta^{-1}|t_k|^2 \right )
\le
  \varepsilon 2 c_0 \beta^{-1} 2^{-k} +  2 \beta^{-1}  |t_{k+1}|  \ ,
\end{equation}
for all $k \ge 0$. Since $2 \beta^{-1} < 1$, using induction in
\eqref{domingo4}
and the fact that $t_k$ is bounded, we get
$|t_k| \le  \varepsilon  c_1 2^{-k}$
with $c_1=2 c_0 \beta^{-1} / \left (1-2
\beta^{-1} \right )$  for all $k \ge 0$. Now this estimate together  
with  Lemma
\ref{l11100}
  gives
$$
\|u_{km}(t_k) - T^{km}(g)\|_{\mathcal B} \le M |t_k| \le
\varepsilon   c_1 M
2^{-k} \ .
$$
Hence, using Lemma   \ref{l100} again,  we get
\begin{eqnarray*}
\label{fdsqwash}
\left \|T^{km}(f)-T^{km}(g) \right \|_{\mathcal B}
& \le &
\|v_k\|_{\mathcal B} +
\left  \|u_{km}(t_k) - T^{km}(g) \right \|_{\mathcal B} \nonumber \\
& \le &
\varepsilon   2^{-k} + \varepsilon   c_1 M 2^{-k} = \varepsilon   c_2
2^{-k}
\ .
\end{eqnarray*}
Therefore, by   (iv) in Lemma \ref{l100},
for all  $i \in \{1,\dots, m-1\}$,
$$
\left \|T^{km+i}(f)-T^{km+i}(g) \right \|_{\mathcal B}
\le M_0  \left \|T^{km}(f)-T^{km}(g)\right \|_{\mathcal B}  \le  
\varepsilon c_3 2^{-k}
\ ,
$$
  which ends the proof.
\hfq

\begin{remark}
\label{gggrad}
Note that since the constant $C_2$ is uniform (independent of~$\varepsilon$) in the above lemma,
inequality \eqref{domingo1}  can be improved to
$$ \left \|T^j(f)-T^j(g) \right\|_{\mathcal B}
\le  C' 2^{- j/m} \left  \|f-g \right\|_{\mathcal B}
     \ ,
$$
where $C'=2C_2$. Therefore, we have exponential
contraction in $\mathcal{B}$
(along the local stable sets) in the strong sense.
\end{remark}

\Subsec{Tangent spaces}
Our next goal is to show that $V_{\varepsilon}^0(g)$ is a $C^1$ manifold
provided $\varepsilon$ is sufficiently small.
The first step towards this goal is to find the natural candidate  for
the tangent space  at every point    $f \in V_{\varepsilon}^0(g)$.
This will be accomplished in  Lemma {\ref{L23}} below.
The proof will require the following elementary bootstrapping result.

\begin{lemma}\label{boo}
Let $(a_n)$ be a sequence of real numbers such that{\rm ,} for some
$c_0 >0 $ and all $n\geq 1${\rm ,}
\begin{equation}\label{boo1}
|a_{n+1}|\leq \frac{1}{4}|a_n| +\frac{c_0}{2^n}\sum_{j=1}^{n-1}|a_j|
\ .
\end{equation}
Then $|a_n|\leq c_1 2^{-n}$ for some $c_1>0$ and all $n\geq 1$.
\end{lemma}

\Proof We may assume that $c_0\geq 1$. Let $n_0>0$ be such that
$c_0 n_0/2^{n_0}<1/2$, and set $b=\max_{1\leq j\leq n_0} \{|a_j|\}$.
Then we see by
induction from \eqref{boo1} that $|a_n|\leq b$ for all $n\geq 1$, and
so
\begin{eqnarray*}
|a_{n+1}| & \le & \frac{1}{4}|a_n| +\frac{nbc_0}{2^n}\\
& \le & \frac{1}{4}|a_n| + {bc_0}\left(\frac{3}{4}\right)^n
\ .
\end{eqnarray*}
By induction, this yields $|a_n|\leq (2bc_0)(\frac{3}{4})^n$ for
all $n\ge 1$.
Hence $\sum_{n=1}^\infty  |a_n| \leq 6bc_0$. Using
\eqref{boo1} once more, we
deduce that
$$
|a_{n+1}|\leq \frac{1}{4}|a_n| +\frac{6bc_0^2}{2^{n}} \ ,
$$
for all $n\geq 1$. Again by induction, this gives $|a_n|\leq
(24bc_0^2)2^{-n}$ for all $n\geq 1$, which is the desired result.
\Endproof

\begin{lemma}
\label{L23}
There exist an integer $m${\rm ,} constants  $C_3 ,C_4 > 0$
and $\varepsilon > 0$ small enough with the following properties.
For every $g \in \mathbb{K}$ and for every
$f \in V^0_\varepsilon(g)$
there exists a linear functional $\theta_{f,g} \in {\mathcal C}^*$
with norm bounded from above by $C_3$ and
with the property that
\begin{equation}
\label{AA}
\left \|DT^j(f)v- {\delta}_0^j   \theta_{f,g}(v)\mathbf{u}_j
\right \|_{{\mathcal C}}
\le C_4  \delta_0^j  2^{-j/m}\|v\|_{{\mathcal C}} \ ,
\end{equation}
for all  $v \in {\mathcal C}$
and all $j \ge 1$.
If $g_0,g_1 \in \mathbb{K}$ and
$f \in  V^0_\varepsilon(g_1) \cap  V^0_\varepsilon(g_2)$
  then
$$
  \theta_{f,g_1}|{\mathcal B}= \theta_{f,g_2}|{\mathcal B}.
$$
  Furthermore{\rm ,} the map
$\Psi:\bigcup_{g \in \mathbb{K}} V^0_\varepsilon(g) \to {\mathcal B}^*$
given by $\Psi(f)=\theta_f=\theta_{f,g}|{\mathcal B}$
\/{\rm (}\/where $g$ is any point of $\mathbb{K}$  such that $f \in   
V^0_\varepsilon(g)${\rm )}
is
well-defined and uniformly continuous.
\end{lemma}

\begin{remark}
\label{remmm}
Condition \eqref{AA} entails that
for every $g \in \mathbb{K}$,
  ${\mathcal B}$ is the
direct sum of the one dimensional unstable  subspace
$E^u_g$  with the kernel of $\theta_f$, i.e.\ ${\mathcal B}=E_g^u  \bigoplus {\rm ker}(\theta_f)$,
provided $f$ is sufficiently close to $g$.
To see this,  note that we can write
$$v = \theta_f(v)
(\theta_f(\mathbf{u}_g))^{-1}\mathbf{u}_g +
\left(v -\theta_f(v) (\theta_f(\mathbf{u}_g))^{-1}\mathbf{u}_g  \right)  
\ .$$
Thus,
from the continuity of  $f \mapsto \theta_f$ plus the fact that
$\theta_g(\mathbf{u}_g)\neq 0$ it follows that if $f$ is close to $g$
then
$\mathbf{u}_g$ is transversal to ${\rm ker}(\theta_f)$.
The hyperplane ${\rm ker}(\theta_f)$ is the natural candidate to
be the tangent space of $V_\varepsilon^0 (f)$ at $f$ since it
corresponds to  all vectors which expand under $DT^j(f)$  by
a factor  less than $\delta_0^j$.
\end{remark}

{\it Proof of Lemma} 6.7.  Let $\varepsilon>0$ be small enough such that Lemma \ref{L22}
is satisfied and  $ \varepsilon C_2 < \nu_m$ (where $\nu_m$ is as in
property $\mathbf{B5}$   in Definition \ref{def4} and $C_2$ is as in
Lemma \ref{L22}). Let $R_k =R_{f,k} = \left (\delta_0^{km} \right
)^{-1} DT^{mk}(f)$ and  write $f_k=T^{km}(f)$, and $g_k=T^{km}(g)$
for all $k \ge 0$. Then 
\begin{equation}
\label{domingo21}
R_{k+1}(v) = \left( \delta_{km}^{(k+1)m}\right) ^{-1} DT^m(f_k) R_k
(v) \ .
\end{equation}
Let us take $v \in {\mathcal{C}}$ with $\|v\|_{\mathcal{C}}=1$. We
can write $R_k(v) = \alpha_k u_{km} + w_k$,
  where $\alpha_k \in \mathbb{R}$ and $w_k \in {\mathcal C}$ are
defined recursively by
$\alpha_0=0$, $w_0=v$ and
\begin{eqnarray}
\label{domingo22}
\alpha_{k+1} & = & \alpha_k + \sigma_{km} (w_k)   \\
w_{k+1} & = & \alpha_{k}  \left( \delta_{km}^{(k+1)m}\right) ^{-1}
(DT^m(f_k)-DT^m(g_k)) \mathbf{u}_{km}   \nonumber\\
&{\ }& +
\left( \delta_{km}^{(k+1)m}\right) ^{-1} (DT^m(f_k)-DT^m(g_k)) w_k \nonumber\\
&{\ }& +  \left( \delta_{km}^{(k+1)m}\right) ^{-1}
Q_{km}^{(k+1)m-1} (w_k)  \ . \nonumber
\end{eqnarray}
Now, by Lemma \ref{L22}, 
\begin{equation}
\label{domingo23}
\left \|f_k-g_k \right \|_{\mathcal{B}} \le \varepsilon C_2 2^{-k} \ .
\end{equation}
Since, by property $\mathbf{B3}$  of Definition \ref{def4}, the map
$f \to DT^m(f) \mathbf{u}_{km}$ is Lipschitz at $f=g_k$ (as a map
from ${\mathcal{B}}$ to ${\mathcal{C}}$),  for all $k$
large enough
\begin{equation}
\label{aaa1}
\left \|DT^m(f_k)\mathbf{u}_{km} -DT^m(g_k)  \mathbf{u}_{km}
\right \|_{\mathcal{C}}
\le c_1 \|f_k-g_k\|_{\mathcal{B}} \le c_2 2^{-k} \ .
\end{equation}
By property $\mathbf{B5}$  in Definition \ref{def4} and
\eqref{domingo23}, for all $m$ large enough we  also have
\begin{equation}
\label{aaab3}
\left \|DT^m(f_k)-DT^m(g_k) \right \|_{\mathcal C}  \le
\frac{\delta_{km}^{(k+1)m}}{8} \ .
\end{equation}
Since,  by property $\mathbf{B1}$   of Definition \ref{def4},
$(\mathcal{C}, \mathcal{D})$ is $\rho$-compatible with
$(T,\mathbb{K})$,   by Corollary \ref{corrob}, for all $m$ large
enough,
\begin{equation}
\label{aaa2}
\left \| Q_{km}^{(k+1)m-1} \right \|_{\mathcal C} \le
\frac{\delta_{km}^{(k+1)m}}{8} \ .
\end{equation}
Using Lemma \ref{l11100} and putting
\eqref{aaa1}, \eqref{aaab3} and \eqref{aaa2}
in \eqref{domingo22} we get
\begin{eqnarray}\label{www}
\|w_{k+1} \|_{\mathcal C}  & \le &  \frac{1}{4} \|w_k \|_{\mathcal C}
+ c_3 2^{-k}  |\alpha_k|    \\
&\leq & \frac{1}{4} \|w_k \|_{\mathcal C}  +
                \frac{c_3 M}{2^k}      \sum_{j=0}^{k-1} \|w_j
\|_{\mathcal C}  \ .\nonumber
\end{eqnarray}
 From   \eqref{www} and Lemma \ref{boo},
we deduce that $\|w_k\|_{\mathcal C} < c_4 2^{-k}$.
Thus, by \eqref{domingo22}  we obtain
$|\alpha_{k+1}-\alpha_k| \le c_5 2^{-k}$  for all $k \ge 0$.
Therefore,
$\theta_{f,g}(v) = \lim \alpha_k$ exists
  and
\begin{equation}
\label{aaa4}
\|R_k(v) - \theta_{f,g}(v)  \mathbf{u}_{km}  \|_{\mathcal C}  \le  c_6
2^{-k}  \ ,
\end{equation}
for all $v \in {\mathcal C}$ with $\|v \|_{\mathcal C}=1$.
If $v \in   {\mathcal C}$ and $\|v \|_{\mathcal C} \ne 1$
then 
$\theta_{f,g}(v)=
\|v \|_{\mathcal C} \theta_{f,g}(v/ \|v
\|_{\mathcal C})$.
By \eqref{aaa4} and by Lemma \ref{l11100}, for all $v,w \in {\mathcal  
C}$,
\begin{eqnarray*}
&&\hskip-36pt\left|\theta_{f,g}(v)+\theta_{f,g}(w)-\theta_{f,g}(v+w)\right|  
  \\
  & \le & M \left \|\theta_{f,g}(v) \mathbf{u}_{km} +\theta_{f,g}(w)
\mathbf{u}_{km}-\theta_{f,g}(v+w) \mathbf{u}_{km} \right \|_{\mathcal C}
   \\
& \le &
M  \left \|\theta_{f,g}(v)\mathbf{u}_{km}-R_k(v) \right \|_{\mathcal C}
   + M
  \left \|\theta_{f,g}(w) \mathbf{u}_{km} -R_k(w) \right \|_{\mathcal C}  
\\
&& + M
  \left \|\theta_{f,g}(v+w) \mathbf{u}_{km} -R_k(v+w) \right  
\|_{\mathcal C} \\
& \le &  c_7  2^{-k} \left (\|v\|_{\mathcal C}+\|w\|_{\mathcal C}  
\right)\ .
\end{eqnarray*}
Hence, letting $k$ go to infinity we deduce that $\theta_{f,g}$ is a
linear functional in~${\mathcal C}^*$. Again by  \eqref{aaa4},
$\|\theta_{f,g}\|_ {\mathcal C}$ is uniformly  bounded and
inequality \eqref{AA} is satisfied  for $j=km$. By  \eqref{aaa4} and
by property $\mathbf{B3}$ in Definition \ref{def4}, for $j=km+i$
with $i \in \{1,\dots,m-1\}$,  
\begin{eqnarray}
\label{ghrted} 
&&\qquad  \left \|R_j(v) - \theta_{f,g}(v)  \mathbf{u}_j \right \|_{\mathcal C}\\
&&\quad \qquad\le 
\left \| \left(\delta_{km}^{km+i}\right)^{-1} DT^i (T^{km} f)
( R_k (v) - \theta_{f,g}(v)  \mathbf{u}_{km}) \right \|_{\mathcal C}  
\nonumber \\
&&\qquad\qquad + \left \|\left(\delta_{km}^{km+i}\right)^{-1}\!\!
\left( DT^i (T^{km} f) -  DT^i (T^{km} g) \right)\!
\theta_{f,g}(v)  \mathbf{u}_{km} \right \|_{\mathcal C}  \nonumber \\
&&\quad\qquad \le  c_8  2^{-k}  \left  \|v \right \|_{\mathcal C}
  \ ,\nonumber
\end{eqnarray}
which proves \eqref{AA}.
In particular,
there is $M_0 > 0$ such that
\begin{equation}
\label{hfggfc}
\left \|R_{k,f}(v) \right \|_{\mathcal C} \le M_0 \ ,
\end{equation}
for all $g \in \mathbb K$, all $f \in V_\varepsilon^0(g)$
and  all $v \in \mathcal B$ with $\|v\|_{\mathcal B} =1$.

Let us   prove that the map $f \mapsto \theta_{f,g}|{\mathcal B}$ is
continuous from $V_\varepsilon^0(g)$ into $\mathcal{B}^*$ for every
$g \in \mathbb K$. By property  $\mathbf{B1}$  of Definition
\ref{def4}, for every $k \ge 1$ the functional $\sigma_{km}$ is
continuous on $\mathcal{C}$ and its norm is uniformly bounded. By
property  $\mathbf{B2}$   in Definition \ref{def4},  the map $f
\mapsto R_{k,f}$ is continuous from $\mathcal{B}$ into
$\mathcal{C}$. Hence,  the mapping $V_\varepsilon^0(g) \to
\mathcal{B}^*$ given by $f \mapsto  \sigma_{km} \circ R_{k,f}$ is
also continuous. By \eqref{ghrted},  
\begin{eqnarray}
\label{sdfsrd}
\left |\sigma_{km} \circ R_{k,f}(v) - \theta_{f,g}(v) \right |  & = &
\left | \sigma_{km} \left ( R_k(v) - \theta_{f,g}(v)
\mathbf{u}_{km}\right)\right |   \\
& \le & c_9  2^{-k}\left \|v \right \|_{\mathcal B} \ .\nonumber
\end{eqnarray}
Therefore, the continuous maps $f \mapsto  \sigma_{km} \circ R_{k,f}$
converge uniformly  to $f \mapsto \theta_{f,g}$, which implies that $f  
\to
\theta_{f,g}$ is also a continuous map from $V_\varepsilon^0(g)$ to
$\mathcal{B}^*$.

Let us   prove that $\theta_{f,g}|{\mathcal B}$
  for $f \in \bigcup_{g \in \mathbb K}
V_\varepsilon^0(g)$ does not depend on $g \in \mathbb K$.
We take $f \in \bigcup_{g \in \mathbb K}
V_\varepsilon^0(g)$ and $g_0, g_1 \in \mathbb K$
such that $f \in V_\varepsilon^0(g_0)$
and $f \in V_\varepsilon^0(g_1)$.
By Lemma \ref{L22}, for every $k \ge 1$,
\begin{eqnarray*}
\label{ffseee}
\left \|T^{km} (g_1)- T^{km} (g_0) \right \|_{\mathcal B} & \le &
\left \|T^{km} (g_1)- T^{km} (f)\right  \|_{\mathcal B} +
\left \|T^{km} (f)- T^{km} (g_0)\right  \|_{\mathcal B} \\
  & \le &  c_{10} \varepsilon 2^{-k}  \ .
\end{eqnarray*}
By property $\mathbf{B1}$ of Definition \ref{def4}, the map $g
\mapsto \sigma_g$ from $\mathbb{K}$ into $\mathcal{C}^*$ is
uniformly continuous. Hence, for every $\varepsilon>0$ there is $k_0
> 0$ large enough such that for all $k>k_0$ and all $w \in {\mathcal
C}$ with $\|w\|_{\mathcal C}\le M_0$,
\begin{equation}
\label{fsfdewfdg}
\left |\sigma_{T^{km} (g_1)}(w)-\sigma_{T^{km} (g_0)}(w) \right | \le
\varepsilon/2  \ .
\end{equation}
By \eqref{hfggfc}, \eqref{sdfsrd} and \eqref{fsfdewfdg} and taking $k$  
large
enough, we get
\begin{eqnarray*}
\left  |\theta_{f,g_1}(v)  - \theta_{f,g_0}(v) \right |  & \le &
\left  |\theta_{f,g_1}(v)  - \sigma_{T^{km} (g_1)} \circ R_{k,f}(v)  
\right| \\
&& +
  \left |\sigma_{T^{km} (g_1)} \circ R_{k,f}(v)  -
\sigma_{T^{km} (g_0)} \circ R_{k,f}(v) \right| \\
&& +
\left  |\sigma_{T^{km} (g_0)} \circ R_{k,f}(v) -\theta_{f,g_0}(v)  
\right | \\
  & \le & 2 c_9  2^{-k} +  \varepsilon/2  \le \varepsilon \,
\end{eqnarray*}
for  all $v \in {\mathcal B}$
with $\|v\|_{\mathcal B}=1$.
Thus, $\theta_{f,g_1}(v)=\theta_{f,g_0}(v)$
and so the map $\Psi$ is
well-defined.

Let us   prove that
the map $\Psi$  is uniformly  continuous.
For every $\alpha_0 > 0$,  we choose $k_0>0$ large enough
such that $2 c_9 2^{-k_0} \le \alpha_0 /3$.
Since the map $g \to \sigma_g$  is uniformly continuous,
there is $\alpha_1 >0$ small enough such
  that for all $g_0,g_1 \in \mathbb K$
with  $\left \|g_1-g_0 \right \|_{\mathcal C} < \alpha_1 $ and all
  $w \in {\mathcal B}$
with $\|w\|_{\mathcal C} \le M_0$,
\begin{equation}
\label{jghhjgdg}
\left |\sigma_{g_1}(w)-\sigma_{g_0}(w) \right | \le \alpha_0/3 \ .
\end{equation}
Choose $k_1>k_0$ large enough such that $\varepsilon C_2
2^{-k_1} \le \alpha_1/3$ where $C_2>0$ is the constant of Lemma
\ref{L22}. Since  $T:{\mathcal O}_B \to {\mathcal C}$ is a $C^1$
operator, by property $\mathbf{B2}$ of Definition \ref{def4}, (and
compactness of $\mathbb K$), there is $\alpha_2>0$ small enough
  such that for all
$f_0 \in V_\varepsilon^0(g_0)$ and
  $f_1 \in V_\varepsilon^0(g_1)$
with $\|f_1-f_0\|_{\mathcal B} < \alpha_2$ we obtain that
$\|T^{k_1m}(f_1)-T^{k_1m}(f_0)\|_{\mathcal C} \le \alpha_1/3$.
Hence, by Lemma \ref{L22},
\begin{eqnarray}
\label{kuyfdvd}&&\\
\left \|T^{k_1m}(g_1)\!-\!T^{k_1m}(g_0) \right \|_{\mathcal C}
& \le &
\left \|T^{k_1m}(g_1)\!-\!T^{k_1m}(f_1)  \right \|_{\mathcal C}
 \!  +
\left \|T^{k_1m}(f_1)\!-\!T^{k_1m}(f_0)  \right \|_{\mathcal C}\nonumber  \\
&& +
\left \|T^{k_1m}(f_0)\!-\!T^{k_1m}(g_0) \right \|_{\mathcal C} \nonumber \\
& \le &  2 \varepsilon C_2 2^{-k_1} + \alpha_1/3 \le \alpha_1 \ .\nonumber
\end{eqnarray}
By \eqref{jghhjgdg} and \eqref{kuyfdvd},  
\begin{equation}
\label{vbfmnh}
\left| \sigma_{T^{k_1m} (g_1)} \circ R_{k_1,f_1}(v) -
\sigma_{T^{k_1m} (g_0)} \circ R_{k_1,f_1}(v) \right |
  \le \alpha_0 / 3 \ .
\end{equation}
Using property $\mathbf{B2}$  of Definition \ref{def4}, choose $0 <
\alpha_3 <  \alpha_2$ small enough
  such that for all $f_1 \in  \bigcup_{g \in \mathbb K}
V_\varepsilon^0(g)$
with $\left \|f_1-f_0 \right \|_{\mathcal B} < \alpha_3$,
$$\left  \| R_{k_1,f_1}(v) -  R_{k_1,f_0}(v)\right  \|_{\mathcal C} \le  
\alpha_0 / (3M)\ ,
$$
 where $v \in {\mathcal B}$ with $\|v\|_{\mathcal B} =1$ and
$M$ is as in  Lemma \ref{l11100}.
Hence,
\begin{equation}
\label{gfghtgs}
\left |\sigma_{T^{k_1m} (g_0)} \circ R_{k_1,f_1}(v) -
\sigma_{T^{k_1m} (g_0)} \circ R_{k_1,f_0}(v) \right | \le
  \alpha_0 /3\ .
\end{equation}
By \eqref{sdfsrd},
\eqref{vbfmnh} and \eqref{gfghtgs}, we obtain that
\begin{eqnarray*}
\left |\theta_{f_1} (v) -\theta_{f_0} (v)\right | & \le &
\left  |\theta_{f_1}(v)  - \sigma_{T^{k_1m} (g_1)} \circ  
R_{k_1,f_1}(v)\right | \\
&& +
\left |\sigma_{T^{k_1m} (g_1)} \circ R_{k_1,f_1}(v) -
\sigma_{T^{k_1m} (g_0)} \circ R_{k_1,f_1}(v)\right | \\
&& +
\left |\sigma_{T^{k_1m} (g_0)} \circ R_{k_1,f_1}(v) -
\sigma_{T^{k_1m} (g_0)} \circ R_{k_1,f_0}(v)\right | \\
&& +
\left |\sigma_{T^{k_1m} (g_0)} \circ R_{k_1,f_0}(v) -
\theta_{f_0} (v) \right |  \\
& \le &
2 c_9 2^{-k_1} + 2 \alpha_0 /3 \le \alpha_0  \ .
\end{eqnarray*}
  Therefore, the map $\Psi$ is uniformly  continuous and Lemma 6.7 is proved.
\hfq

\Subsec{The main estimates}
Besides aiming at proving that the local stable set  is a $C^1$  
manifold,
we want to show
that the local hyperbolicity picture holds (in $\mathcal B$) near
$\mathbb{K}$.
In other words we want to show that if the iterates $T^{km}(f_1)$ of a  
point
$f_1 \in  {\mathcal B}(g,\varepsilon)$ remain in ${\mathcal  
B}(T^{km}(g),\varepsilon)$
for a long time, that is for $k=0,1,\dots,N$ with $N$ large,
then $f_1$ has to be very close to a point $f_0$ on the stable set
$W^s_\varepsilon(g)$ at the outset, and in the end $T^{Nm}(f_1)$
has to be very close to the unstable manifold  
$W_\varepsilon^u(T^{Nm}(g))$.

To prove these facts, we consider in this section (see Lemma \ref{L28})
an intermediate time $l$ for which we can find a good quantitative
estimate for the point on the unstable manifold
$W_\varepsilon^u(T^{lm}(g))$ that best approximates $T^{lm}(f_1)$.
This estimate is provided by the value of $\theta_{f_0}(f_1-f_0)$,
and its most important consequence
is obtained when $f_1$ also belongs to the local
stable set $W^s_\varepsilon(g)$. In this case we prove
an inequality of the form
$\left | \theta_{f_0}(f_1-f_0) \right | \le
C \|f_1-f_0\|_{\mathcal  B}^{1+\tau}$
(see Lemma \ref{fffcccooo}).
As we shall see in \S \ref{dfgfredvgv},
this is precisely what is needed to show that the tangent space
to the stable set at $f_0$ varies continuously with~$f_0$.

In this section we  fix  $m$ large enough
and $\varepsilon_0 > \varepsilon_1 > \varepsilon_2$ small enough such
that Lemmas  \ref{l100},
  \ref{L22} and \ref{L23} are satisfied
  for all $\varepsilon < \varepsilon_2 $ suficiently small.

\begin{lemma}
\label{L26}
There exist constants $C_5, C_6,\varepsilon >0$ with the
following property. For
all $g \in \mathbb{K}${\rm ,} all
$f_0 \in V_{\varepsilon}^0(g)${\rm ,}
and all $f_1\in {\mathcal B}(g,\varepsilon)$   such that
$\|f_1-f_0\|_{\mathcal{C}}\leq C_5 \left(\delta_0^n\right)^{-1}${\rm ,}
\begin{equation}
\label{E2}
\left  \|T^k(f_1)-T^k(f_0)-DT^k(f_0)\left (f_1-f_0\right ) \right
\|_{{\mathcal C}}
  \le C_6  \left(\delta_k^n\right)^{-\rho} \  ,
\end{equation}
for all $0\le k \le \hat k_0(g,f_1)${\rm ,}
where
$$
\hat k_0(g,f_1)= \min \left  \{j \in \{0,\dots, n\}:  
\|T^j(f_1)-T^j(g)\|_{{\mathcal
B}} \ge \varepsilon_0 \right \}
$$ and
$\rho >1$ is as in property $\mathbf{B4}$ of Definition {\rm \ref{def4}.}
\end{lemma}

\Proof
  By Lemmas   \ref{l11100} and \ref{L23},
  there are $c_0, c_1 > 0$ and $\lambda > 1$ such that
\begin{equation}
\label{sssppp}
\left  \|DT^{k-i}\left (T^i f_0\right  )\right \|_{\mathcal C}
\le c_0 \delta_i^k \ \text{and} \
\delta_i^k > c_1 \lambda^{k-i}
\end{equation}
for  all $0 \le i < k$.
Define a sequence $v_i \in {\mathcal C}$
as follows: $v_0 = T (f_1)-T (f_0)$ and
$$
v_i=T^i(f_1)-T^i(f_0)-DT\left (T^{i-1}(f_0)\right )
\left(T^{i-1}(f_1)-T^{i-1}(f_0) \right) \ ,
$$
   for all $0 < i \le k$.
Hence,
\begin{equation}
\label{eqw5ese}
T^i(f_1)-T^i(f_0)-DT^i(f_0)\left (f_1-f_0\right ) =
\sum_{j=1}^i DT^{i-j}\left (T^jf_0 \right)v_j \ .
\end{equation}
Applying property $\mathbf{B4}$ of Definition \ref{def4}, we get
\begin{eqnarray}
\label{aacccaa}
\|v_{i+1}\|_{{\mathcal C}}
& \le & c_2 \left \|T^i(f_1)-T^i(f_0)\right \|_{\mathcal C} ^\rho
\nonumber \\
  & \le & c_2 \left
\|DT^i(f_0)\left (f_1-f_0 \right ) +
\sum_{j=1}^i DT^{i-j}(T^jf)v_j
\right  \|_{\mathcal C}   ^\rho \ .
\end{eqnarray}
Let us first choose $C_6 > 0$  such that
\begin{equation}
\label{eeesss}
C_6^{1-\rho}  >   c_2 \left( \frac{ 2 c_0  c_1^{1-\rho}}
{(\delta_i)^\rho (1-\lambda^{1-\rho})} \right)^\rho \ ,
\end{equation}
for all $0<i \le k$, and then choose $C_5>0$ such that
\begin{equation}
\label{eeesss3}
   C_5^\rho   <   \frac{C_6 (\delta_0)^\rho}{c_2} \  , \qquad
   C_5   <  
\frac{c_1^{1-\rho} C_6  (\delta_i)^{1-\rho}}
{1-\lambda^{1-\rho} } \ .
\end{equation}
Let us  prove inductively for $i=1,2,\dots,k$ that
$\left \|v_i \right \| _{{\mathcal C}}
  \le C_6 \left( \delta_i^n \right)^{-\rho}$.
Using   inequality \eqref{aacccaa} and \eqref{eeesss3}, we get
$$ \|v_1\|_{{\mathcal C}}
   \le  c_2 \left \|  f_1-f_0   \right \|_{{\mathcal C}}^\rho
   \leq  \frac{c_2 C_5^\rho}{(\delta_0)^\rho}
(\delta_1^n)^{-\rho}
   \le  C_6 (\delta_1^n)^{-\rho} \ .$$
Using the inequalities \eqref{sssppp}, \eqref{aacccaa},
  \eqref{eeesss} and  \eqref{eeesss3}, we
get
  \begin{eqnarray}
  \label{jdeied7}
  \|v_{i+1}\|_{{\mathcal C}}
  & \le & c_2 \left( C_5 c_0   \delta_0^i \left( \delta_0^n
\right)^{-1}
  + \sum_{j=1}^i c_0 \delta_j^i C_6 \left( \delta_j^n \right)^{-\rho}
  \right)^{\rho}  \\
  & \le & c_2 c_0^{\rho}
  \left(
  \frac{C_5}{\delta_i} +
  \frac{C_6}{\left( \delta_i \right)^{\rho}}
  \sum_{j=1}^i \left( \delta_j^i \right)^{1-\rho}
  \right)^{\rho}
  \left( \delta_{i+1}^n \right)^{-\rho} \nonumber  \\
& \le & c_2 c_0^{\rho}
  \left(
  \frac{C_5}{\delta_i} +
  \frac{  c_1^{1-\rho}C_6}{\left( \delta_i \right)^\rho
\left(1-\lambda^{1-\rho}\right) }
  \right)^{\rho}
  \left( \delta_{i+1}^n \right)^{-\rho}\nonumber  \\
& \le & C_6 \left( \delta_{i+1}^n \right)^{-\rho} \ ,\nonumber
\end{eqnarray}
which ends the induction.
Thus, using   \eqref{sssppp} and \eqref{jdeied7}
in \eqref{eqw5ese},
we   get
\begin{multline*}
  \left \|T^k(f_1)-T^k(f)-DT^k(f)\left (f_1-f_0 \right
) \right \|_{{\mathcal C}}
   \leq  
\sum_{i=1}^k   c_0 \delta_i^k
    C_6 \left( \delta_i^n \right)^{-\rho}    \\
 \le  
c_0 C_6 \left( \delta_k^n \right)^{-\rho}
\sum_{i=1}^k \left( \delta_i^k \right)^{1-\rho}
\;  \leq \;  \frac{c_0 c_1^{1-\rho}C_6 }{1-\lambda ^{1-\rho}}
\left( \delta_i^n \right)^{-\rho}
\ .
\end{multline*}
This proves the lemma.
\hfq

\begin{lemma}
\label{L27}
Let $C_5, C_6,  \varepsilon >0${\rm ,} $\rho >1$ and $\hat k_0(g,f_1)$
be as in Lemma {\rm \ref{L26}.}
There exist $C_7,C_8 > 0$    such that
for all $g \in \mathbb{K}${\rm ,} all $f_0 \in V_{\varepsilon}^0(g)$
and all $f_1\in {\mathcal B}(g,\varepsilon)$   such
that $\|f_1-f_0\|_{\mathcal{C}}\leq C_5 \left(\delta_0^n\right)^{-1}${\rm ,}
\begin{multline}
\label{E3}
 \left \|T^k(f_1)-T^k(g)    -\delta_0^k
  \theta_{f_0} \left( f_1-f_0 \right){\mathbf u}_k \right
\|_{{\mathcal C}}   \\
   \le C_6
\left(\delta_k^n \right)^{-\rho} +  {\varepsilon}  C_7 2^{-k/m} +  C_8  
2^{-k/m}
\left(\delta_k^n \right)^{-1}  \ ,
\end{multline}
  for all $k \le \hat k_0(g,f_1)$.
\end{lemma}

\Proof
  By Lemma \ref{L23},    
\begin{equation}
\label{E4}
\left \|DT^k(f_0) \left( f_1-f_0 \right) -  \delta_0^k
\theta_{f_0}\left( f_1-f_0 \right) {\mathbf u}_k \right \|_{{\mathcal
  C}} \le C_4 2^{-k/m} \left(\delta_k^n \right)^{-1}   \ .
\end{equation}
By  Lemma \ref{L22}, we obtain that
\begin{equation}
\label{E455}
\left \|T^k(f_0)- T^k(g)
\right \|_{{\mathcal C}} \le c_0
\left \|T^k(f_0)- T^k(g)
\right \|_{{\mathcal B}} \le  {\varepsilon} c_0 C_2 2^{-k/m} \ .
\end{equation}
  Combining \eqref{E2}, \eqref{E4} and \eqref{E455}, we get \eqref{E3}.
\hfq

\begin{definition}
Given   $g \in \mathbb{K}$ and  $p \ge 1$
we denote by  $l=l(g,p)$  the {\it smallest integer} such that
\begin{equation}
\label{fghqw}
  \left( \delta_{lm}^{pm} \right)^\rho \le  2^l \ ,
\end{equation}
where $\rho >1$ is as in property $\mathbf{B4}$ of Definition
\ref{def4}.
\end{definition}

\begin{lemma}
\label{L278}  
\begin{enumerate}
\item
There exist    $0 < \mu_0 <  \mu_1 < 1$ with
the property that       $\mu_0 p \le l=l(g,p) \le \mu_1 p$
for all $g \in \mathbb{K}$ and all $p \ge 1$.
\item
There exists  $0 <\tau_1 < 1$  such that
for all  $g \in \mathbb{K}$
and all  $f_0,f_1 \in {\mathcal B}(g,\varepsilon)${\rm ,}
if $\|f_0-f_1\|_{\mathcal C} \ge C_5
\left(\delta_0^{(p+1)m}\right)^{-1}$
then $\left(\delta_{l  m}^{pm}\right)^{-1}
  \le C_9 \|f_0-f_1\|_{\mathcal C}^{\tau_1}${\rm ,}
where $C_9$ depends only  upon $C_5>0$.
\end{enumerate}
\end{lemma}

\Proof
Let us prove part (i).
Set $1 < \lambda < M$   as in Lemma \ref{l11100}.
Then, by \eqref{fghqw},
$$
\rho \log M^{-1} +({pm}-lm) \rho \log \lambda \;<\; \rho \log
\delta_{lm}^{pm} \;<\; l   \log 2 \;<\; l m \log 2 \ .
$$
Hence,
    $$
\left (1+
\frac{\log 2} {\rho \log \lambda}
\right) lm \;\geq\;
\frac{\log M^{-1}} {\log \lambda} +{pm} \;>\;
\frac{pm}{2} \ ,
$$
for all $p$ such that  ${pm} > N_0= \max \{ 2m, \left | 2 \log M^{-1} /  
\log
\lambda \right | \}$.
Thus, taking
$$
\mu_0' = 2^{-1} \left (1+
\frac{\log 2} {\rho \log \lambda}
\right)^{-1} > 0 \ ,
$$
  we get $ \mu_0' {p}  \le l$ for all such values of $p$.
By \eqref{fghqw} and by Lemma \ref{l11100} there exists
  a uniform constant $0<c_0 \le 1$ such that
$c_0 2^l \le  \left( \delta_{lm}^{pm} \right)^\rho$
and so
$$
\log c_0 + l \log 2 < \rho \log  \delta_{lm}^{pm}  \le \rho ({pm}-lm)  
\log M \ .
$$
  Letting $\alpha= \log 2 / \left ( \rho \log M \right ) > 0$
and $\beta= \log c_0 /  \left ( \rho \log M \right )$
  we get
$$
lm \left(1 + \frac{\alpha}{m}\right) \le pm-\beta
\le pm \left (1+\frac{\alpha}{2m} \right)
$$
    for all ${p} > -2 \beta/\alpha $.
Thus, taking
$$
0 < \mu_1' = \frac{2m+\alpha}{2(m+\alpha)}   < 1
$$
we obtain that $lm \le \mu_1' {pm}$ for
all such values of $p$.
Since $\delta_g$ varies continuously with $g$ in the
compact set ${\mathbb K}$, we can extend the previous results to all
$p \ge 0$ for some $\mu_0 \le \mu_0'$ and $\mu_1 \ge \mu_1'$.

Let us prove part (ii).
Take $0 <\tau_1=(1-\mu_1) \log M/\log
\lambda <1$.
Then, by Lemma \ref{l11100},  
\begin{eqnarray*}
\left(\delta_{lm}^{pm}\right)^{-1}
  & \le & c_0 \lambda^{-(p-l)m} \le
c_0 \lambda^{-(1-\mu_1)pm} \\
& \le & c_0 M^{-\tau_1 p m} \le c_1 \omega
\left( \delta_0^{(p+1)m} \right)^{-\tau_1}\\
& \le & c_1 \|f_0-f_1\|_{\mathcal C}^{\tau_1} \ .
\end{eqnarray*}
\vglue-24pt
\Endproof

\begin{lemma}
\label{L28}
There exist $\varepsilon > 0$ sufficiently small and $C_{10} >0$
such that the following holds
for  $g \in \mathbb{K}${\rm ,}
       $f_0 \in V_\varepsilon^0(g)$
and   $f_1\in {\mathcal B}(g,\varepsilon)$.
If $p$ is the largest integer such that
$$\|f_1-f_0\|_{\mathcal{C}}\leq C_5
\left(\delta_0^{pm}\right)^{-1}$$
then $l=l(g,p) \le k_0=k_0(g,f_1)$
  and so $t_{l}=t_{l}(g,f_1)$ is well-defined
\/{\rm (}\/where  $l$ is as in Lemma {\rm \ref{L278},}
$k_0$ and  $t_{l}$ are  as in Lemma  {\rm \ref{l100},} and
$C_5  >0$ is   as in {\rm Lemma \ref{L27}).}
Furthermore{\rm ,}
\begin{equation}
\label{efgert}
  \left  |t_{l}- \delta_{0}^{lm} \theta_{f_0}(f_1-f_0) \right |
    \le   C_{10} \left( \delta_{lm}^{pm} \right)^{-\rho} 
\end{equation}
where   $\rho >1$ is   as in Lemma {\rm \ref{L27}.}
\end{lemma}

\Proof Let $\hat k_0=\hat k_0(g,f_1)$ be as in Lemma \ref{L27}.
Let us prove that $l \le k_0$.
By (iv) in Lemma  \ref{l100},
$m k_0 < \hat k_0$. Hence it is enough to prove that
$ \min \left \{l m, \hat k_0 \right \} \le m k_0$.
Let  $\varepsilon > 0$ be small enough
such that Lemmas \ref{L26} and \ref{L27}
are satisfied.
Let us show that $lm \le  \hat k_0$.
  By inequality  \eqref{E3}, for all $k$ such that $ m k \le \hat k_0$  
we have
\begin{multline}
\label{gedwedfg}
 \left |\sigma_{km} \left ( T^{km}(f_1)-T^{km}(g) \right  
)\right |
  \\
  \le
   \|\sigma_{km} \|_{\mathcal C}
  \left( \delta_0^{km} \left |\theta_f \left (f_1-f_0 \right )\right |+
  c_0 \left( \delta_{km}^{pm} \right)^{-1}
  + \varepsilon c_1 2 ^{-k} \right ) \ .
\end{multline}
By Lemma \ref{l11100} and Remark \ref{remmm},
  there is $M_1 >1$ such that
  $M_1^{-1} \le \|\sigma_{km}\|_{\mathcal C}\break \le M_1
  $ and
   $M_1^{-1} \le \|\theta_{f_0}\|_{\mathcal C} \le M_1$.
Since by   Lemma \ref{l11100}, we have
$\left(\delta_{km}^{pm}   \right)^{-1} \le M \lambda^{-(p-k)m}$
we deduce that
\begin{equation}
\label{gedwedfg3}
  \delta_0^{km} |\theta_f  \left (f_1-f_0 \right )|   \le
\left(\delta_{km}^{pm}   \right)^{-1}
\|\theta_f\| _{\mathcal C}
\le  M M_1 \lambda^{-(p-k)pm} \ .
\end{equation}
By  Lemma  \ref{L278}, there is $0 <  \mu_1 < 1$
such that for all $p>0$ and all ${k} \le l$ we have
  ${p}-{k} \ge {p}-l  \ge (1- \mu_1) {p}$.
Now, we make    $\varepsilon > 0$   small  enough
(and so    $p$   large enough) such that
the following inequalities are satisfied
  \begin{eqnarray*}
(c_0+M_1) M \lambda^{-(1-\mu_1)pm} & <  &
\frac{ \varepsilon_2}{2   \|\sigma_{km} \|_{\mathcal C}  }  \ ,
   \\
  \varepsilon  c_1 2 ^{-k}  & <  &
  \frac{ \varepsilon_2}{4   \|\sigma_{km} \|_{\mathcal C}  }  \ ,
\end{eqnarray*}
for  all $k$ such that ${km}   \le \min \{lm, \hat k_0\}$.
  Therefore, for all such $k$, combining \eqref{gedwedfg} and  
\eqref{gedwedfg3}
we deduce that
\begin{equation}
\label{wwwgggwww}
|\sigma_{km} \left (T^{km}(f_1)-T^{km}(g) \right)| < \varepsilon_2  \ .
\end{equation}
Since $f_1 \in   {\mathcal B}(g,\varepsilon)$ and \eqref{wwwgggwww}
reverses the inequality (iii) in Lemma \ref{l100},
we obtain that $\min \{lm, \hat k_0\} \le m k_0$, and so $l \le k_0$.

Now, let us prove \eqref{efgert}.
Since $l \le    k_0$, by \eqref{sadsa} and \eqref{fgddfg}
in  Lemma \ref{l100},
there is $t_l=t_{l}(g,f_1)$ such that
\begin{equation}
\label{E6}
\|T^{lm}(f_1)-u_{lm}(t_l)\|_{{\mathcal B}}
  \le \varepsilon 2^{-l}
\le   \varepsilon \left(\delta_{lm}^{pm} \right)^{-\rho} \ .
\end{equation}
Since  $lm  < \hat k_0$,
by lemmas \ref{L27} and \ref{L278} we get
$$
\|T^{lm}(f_1)- T^{lm}(g) -
   s_l \mathbf{u}_{lm} \|_{{\mathcal C}}
\le c_2 \left(\delta_{lm}^{pm} \right)^{-\rho} \ ,
$$
where $s_l=\delta_0^{lm} \theta_{f_0}(f_1-f_0)$.
   Thus, using \eqref{E6},  we obtain that
$$
\|u_{lm}(t_l)-  T^{lm}(g)
-s_l \mathbf{u}_{lm}\|_{{\mathcal C}}
  \le c_3 \left(\delta_{lm}^{pm} \right)^{-\rho} \ .
$$
  Since $t \to u_{lm}(t)$ is $C^2$
as a map ${\mathbb R} \to {\mathcal C}$, 
\begin{eqnarray*}
\| u_{lm}(s_l)  -  T^{lm}(g)
  -s_l \mathbf{u}_{lm}\|_{{\mathcal C}}
  &  \leq  &   c_4  s_l^2 \\
&  =  &   c_4  \left | \delta_0^{lm} \theta_{f_0}(f_1-f_0) \right
|^2 \\
  &  \le   &
c_5 \left (\delta_{lm}^{pm} \right )^{-2}
\ .
\end{eqnarray*}
Therefore,
\begin{eqnarray*}
\lefteqn{\|u_{lm}(t_l)-u_{lm}(s_l)\|_{{\mathcal C}}}\\[4pt]
&&\ \ \ \ \   \le    \|u_{lm}(t_l)-  T^{lm}(g)
-s_l \mathbf{u}_{lm}\|_{{\mathcal C}} +
\| u_{lm}(s_l)  -  T^{lm}(g)
  -s_l \mathbf{u}_{lm}\|_{{\mathcal C}} \nonumber \\[4pt]
&&\ \ \ \ \   \le     c_3 \left (\delta_{lm}^{pm} \right )^{-\rho}+
c_5 \left (\delta_{lm}^{pm} \right )^{-2} \nonumber  \\[4pt]
&&\ \ \ \ \   \le     c_6 \left (\delta_{lm}^{pm} \right )^{-\rho}
\ ,
\end{eqnarray*}
because $1<\rho <2$.
Hence, applying Lemma \ref{l11100}, we get
\vglue12pt
\hfill $
\displaystyle{|t_l-s_l|\le M^{-1}
\|u_{lm}  (t_l) -
u_{lm}  (s_l)
\|_{{\mathcal C}}
\le c_7 \left (\delta_{lm}^{pm} \right )^{-\rho} \ .}
$ 
\Endproof

\begin{lemma}
\label{fffcccooo} \label{P24} \label{jkht}
There exist constants $\tau, \varepsilon, C>0$
with the following properties\/{\rm :}\/ for all $g \in \mathbb K$ and
  all  $f_0,f_1 \in V^0_\varepsilon(g)${\rm ,}
$$
\left | \theta_{f_0}(f_1-f_0) \right | \le
C \left \|f_1-f_0 \right \|_{\mathcal  B}^{1+\tau}
\ .
$$
\end{lemma}

\Proof
We shall in fact prove a stronger inequality, with the  
$\mathcal{C}$-norm
replacing the  $\mathcal{B}$-norm.
Let $\varepsilon>0$ be so small  that Lemmas \ref{L26}  to \ref{L28}   
are
satisfied ($\varepsilon>0$ will be made even smaller in the course
of the argument).
Let $p$ be such that
$C_5 \delta_0^{-(p+1)m} <   \|f_0-f_1\|_{\mathcal C} \le  C_5  
\delta_0^{-pm}$
where  $C_5>0$ is as in Lemma \ref{L26}.
As in Lemma \ref{l100}, set
  $k_0=k_0(g,f_1)$, $t_j=t_j(g,f_1)$ and $v_j=v_j(g,f_1)$
for all $0 \le j \le k_0$.
Also, let  $l=l(g,p)$ be as in Lemma \ref{L278}.
By Lemma \ref{L28}, we have
$l \le k_0$ and so $t_l$ is well-defined.
Thus, applying Lemmas \ref{l100}, \ref{L22} and \ref{L278}, we get
  \begin{eqnarray*}
\left \|u_{lm}(t_l)-u_{lm}(0)\right \|_{{\mathcal B}} & \le &
  \left \|u_{lm}(t_l)-T^{lm}(f_1)\right \|_{{\mathcal B}}  +
  \left \|T^{lm}(f_1)-T^{lm}(g)\right \|_{{\mathcal B}} \nonumber \\[4pt]
   & \le &    \varepsilon c_0 2^{-l} \le  \varepsilon c_0
   \left(\delta_{lm}^{pm}\right)^{-\rho} \ .
\end{eqnarray*}
Hence, by Lemma \ref{l11100} we see that $|t_l| \le
c_1 \left(\delta_{lm}^{pm}\right)^{-\rho}$.
Let us write $t_j = \alpha_j \left(\delta_{jm}^{pm}\right)^{-1}$
  for $l \le j \le k_0(g,f_1)$.
Recalling that $\delta_{jm}^{(j+1)m} > \beta > 2$ for all $j$ and
using  Lemma \ref{L278}, we have
\begin{equation}
\label{gggggg}
  \left( \delta_{lm}^{pm} \right )^{-1}  \le \beta^{-(1- \mu_1)p}
\ \text{and} \
\left(\delta_{lm}^{pm}\right)^{-(\rho-1)/2}    \le \beta^{-\tau_2 p}
\end{equation}
where $\tau_2=(1-\mu_1)(\rho-1)/2$.
  Hence, making $\varepsilon>0$   smaller if necessary
  (and so $p$ large enough), we get
\begin{equation}
\label{rgertrte}
\alpha_l <  4^{-1} \left(\delta_{lm}^{pm}\right)^{-(\rho-1)/2}
<  4^{-1} \beta^{-\tau_2 p}
< \varepsilon_1/2 \ .
\end{equation}
By   Lemma \ref{l100}, we have $\|v_j\|_{\mathcal B} \le \varepsilon  
2^{-j}$
and $t_{j+1}= \hat \delta_{jm}^{(j+1)m} (t_j+ \sigma_{jm}(v_j))$.
  Since $t \mapsto \hat \delta_{jm}^{(j+1)m} (t)$ is $C^2$
as a map ${\mathbb R} \to {\mathcal C}$,
we deduce that
\begin{eqnarray*}
\label{seseWw}
\left| t_{j+1} - \delta_{jm}^{(j+1)m}  t_j \right| & \le &
c_2 \left (|t_j|^2 + \|v_j\|_{\mathcal B} \right)   \\[4pt]
     & \le &    c_2 \left( |t_j|^2  + \varepsilon 2^{-j} \right)  
\nonumber \\[4pt]
   & \le &   c_2   \left( |t_j|^2  +  \varepsilon 2^{-(j-l)}
  \left( \delta_{lm}^{pm} \right )^{-\rho} \right )  \ .
\end{eqnarray*}
Therefore, 
\begin{equation}
\label{eeezzzaaa}
| \alpha_{j+1}       - \alpha_j |
   \le
  c_2 \left( \alpha_j^2 \beta^{-(p-j+1)}
+  \varepsilon (2\beta)^{-(j-l)} \left( \delta_{lm}^{pm} \right
)^{-(\rho-1)} \right) \ .
\end{equation}
We prove  that   $k_0(g,f_1) > p$.
To do this, we need to show that $|t_j| < \varepsilon_1$
for all $j \le p$.
We prove by induction a slightly stronger statement, namely,
that $\alpha_j \le 2^{-1} \left(\delta_{lm}^{pm}\right)^{-(\rho-1)/2}<
\varepsilon_1$ for all $j=l,\dots,p$.
This is certainly satisfied for $j=l$, as can be seen from  
\eqref{rgertrte}.
Suppose it is satisfied for $\alpha_i$ for all $i=l,\dots,j$.
Using  \eqref{gggggg} and \eqref{eeezzzaaa}, and  making $\varepsilon>0$
even smaller
(and thus $p$ large enough), we get
\begin{eqnarray}
\label{fererfefef}\qquad\quad
| \alpha_{j+1}       - \alpha_l |  & \le &
\sum_{i=l}^j  | \alpha_{i+1}       - \alpha_i |  \\[4pt]
  & \le & \frac{1}{4}
\left( \sum_{i=l}^j
  \left(  c_2   \beta^{-(p-i+1)}+
4 c_2 \varepsilon (2\beta)^{-(i-l)} \right) \right)
\left( \delta_{lm}^{pm} \right )^{-(\rho-1)} \nonumber \\[4pt]
  & \le & \frac{1}{4}
  \left( \frac{   c_2 \beta^{-1}}{ 1- \beta^{-1}}
  + \frac{ 4c_2 \varepsilon}
  { 1- (2\beta)^{-1}} \right) \beta^{-\tau_2 p}
\left(\delta_{lm}^{pm}\right)^{-(\rho-1)/2} \nonumber  \\[4pt]
  & \le &\frac{1}{4}
   \left(\delta_{lm}^{pm}\right)^{-(\rho-1)/2} \ .\nonumber
  \end{eqnarray}
Since  $\alpha_l \le 4^{-1}
\left(\delta_{lm}^{pm}\right)^{-(\rho-1)/2}$, we deduce that
   $\alpha_{j+1} \le   2^{-1}
\left(\delta_{lm}^{pm}\right)^{-(\rho-1)/2} < \varepsilon_1$
(in particular  $j+1 \le k_0(g,f_1)$)
  which ends the induction.

Now set $s_j=s_j(g,f_0,f_1)= \delta_0^{jm} \theta_{f_0}(f_1-f_0)$ for  
all $j$.
Let us estimate $|t_p-s_p|$.
By Lemma  \ref{L28} and the above estimates on $\alpha_j$'s,
       we have
       \begin{equation}
       \label{sffswewehy}
       \left |t_p-s_p \right |   \le
      \left  |\alpha_p-\alpha_l \right | +
       \delta_{lm}^{pm}  \left |t_l - s_l\right |
        \le  c_3 \left(\delta_{lm}^{pm}\right)^{-(\rho-1)/2} \ .
  \end{equation}
On the other hand, from Lemmas \ref{l100}, \ref{L22} and \ref{L278},
we also know that
  \begin{eqnarray*}
\left \|u_{pm}(t_p)-u_{pm}(0)\right \|_{{\mathcal B}} & \le &
  \left \|u_{pm}(t_p)-T^{pm}(f_1)\right \|_{{\mathcal B}}  +
  \left \|T^{pm}(f_1)-T^{pm}(g)\right \|_{{\mathcal B}} \nonumber \\[4pt]
   & \le &    \varepsilon c_4 2^{-p}\ .
\end{eqnarray*}
Hence, again by Lemma \ref{l11100},
we have $|t_p| \le \varepsilon c_5 2^{-p}$.
Since $p \ge l$, we deduce from Lemma \ref{L278} that
\begin{equation}
\label{gfiuy}
|t_p| \le    \varepsilon c_5 2^{-l} \le
\varepsilon c_5 \left(\delta_{lm}^{pm}\right)^{-\rho}
\ .
\end{equation}
But   Lemma \ref{L278}, also tells us that
there exists  $\tau_1>0$ such that
$\left(\delta_{lm}^{pm}\right)^{-1} \le c_6 \|f_0-f_1\|_{\mathcal  
C}^{\tau_1}$.
Moreover,   $\left(\delta_0^{pm}\right)^{-1} \le c_7
\|f_0-f_1\|_{\mathcal C}$ by hypothesis.
  Therefore, combining these facts with   \eqref{sffswewehy} and
\eqref{gfiuy}, we get at last
\begin{eqnarray*}
\left| \theta_{f_0}(f_1-f_0) \right| & \le &
\left(\delta_0^{pm}\right)^{-1}  |s_p|   \\
  & \le &
\left(\delta_0^{pm}\right)^{-1} \left( |t_p|+ |t_p-s_p| \right) \\
  & \le &
\left(\delta_0^{pm}\right)^{-1}
\left( \varepsilon c_5 \left(\delta_{lm}^{pm}\right)^{-\rho}
+ c_3
\left(\delta_{lm}^{pm}\right)^{-(\rho-1)/2} \right) \\
& \le & c_8 \|f_0-f_1\|_{\mathcal C}^{1+ \tau_1(\rho-1)/2}
\ ,
\end{eqnarray*}
which finishes the proof.
  \hfq

\Subsec{The local stable sets are graphs}
\label{dfgfredvgv}
We shall prove now that the local stable set of every $g_0 \in \mathbb  
K$
in a sufficiently small neighborhood of $g$ is the graph of a function
defined over ${\rm Ker} \theta_g \cap {\mathcal B}$
(and with values on the one-dimensional subspace
$\mathbb{R} \mathbf{u}_g \subset \mathcal{B}$).
The idea is to show that every ``vertical line'' of
the form $f+\mathbb{R} \mathbf{u}_g$ with $f$ close to $g$
cuts the local stable set $W^s_\varepsilon (g_0)$
exactly at one point.
All other points in the same vertical
line escape exponentially fast away from
$W^s_\varepsilon (g_0)$ under iteration by $T$
and the time $k_0m$ each such point $f$ takes to escape is
logarithmic on the reciprocal of its distance to $W^s_\varepsilon  
(g_0)$.
Moreover, $T^{k_0m}(f)$ will be exponentially close (in $k_0$)
to $W^u_\varepsilon(T^{k_0m}(g_0))$.

\begin{proposition}
\label{dfgfdger}
There exist $0<\alpha_0, \alpha_1, \alpha_2< \varepsilon${\rm ,}
  $0<\mu_0 < \mu_1$ and $M_0>1$ with the following properties.
If $g_0  \in \mathbb K$
  and $g  \in \mathbb K$ is such that
  $\|g -g_0\|_{\mathcal B} < \alpha_0${\rm ,}
then for  every
$v \in {\rm Ker}~ \theta_{g_0} \cap \mathcal B$
with $\|v\|_{\mathcal B} < \alpha_1${\rm ,}
there exists $-\alpha_2/2<\tau(g ,v)<\alpha_2/2$ such that
\begin{enumerate}
\item
$f_{\tau(g ,v)}= g_0 +v+\tau(g,v) {\mathbf u}_{g_0} \in
   W^s_\varepsilon(g) \subset V_\varepsilon^0(g)${\rm ;}
\item
$f_{t}= g_0 +v+t {\mathbf u}_{g_0} \in
V_\varepsilon^+(g)$ for all
$\tau(g,v)<t< \alpha_2${\rm ;}
\item
$f_{t}= g_0 +v+t {\mathbf u}_{g_0} \in
V_\varepsilon^-(g)$ for all
  $-\alpha_2<t< \tau(g,v)${\rm ;}
\item
$-\mu_0 \log (|t-\tau(g,v)|)\le k_0(g,f_t) \le - \mu_1 \log  
(|t-\tau(g,v)|)${\rm ,}
where  $k_0(g,\phi_t)$ is as in Lemma {\rm \ref{l100}.}
\end{enumerate}
\end{proposition}

\Proof
Let $\varepsilon>0$ be sufficiently small such that Lemmas \ref{L26}   
to \ref{L28}
are satisfied and $0< \varepsilon' < \varepsilon$ such that Lemma  
\ref{L22}
is satisfied.
Let $M>0$ be as in Lemma \ref{l11100} and
take positive numbers $\alpha_1$ and $\alpha_2$ such that
\begin{equation}
\label{btgrd}
0< 8 \alpha_1 M < \alpha_2 \ \text{and} \
  \alpha_1+ 2\alpha_2 M < \varepsilon'/2 \ .
\end{equation}
Take $g \in  {\mathbb K}$ and
  $f  \in V_{\varepsilon}(g_0)$ with
  $\|f-g\|_{\mathcal B} <\varepsilon'/2$.
Let  $v \in {\rm Ker}~ \theta_{g_0} \cap \mathcal B$
with $\|v\|_{\mathcal B} < \alpha_1$,
and $t \in \mathbb R$ with $2 M \|v\|_{\mathcal B}  < |t| < 2\alpha_2$.
By the second inequality in \eqref{btgrd}, we have $\phi_t=f+v+t  
\mathbf{u}_{g_0} \in
{\mathcal B}(g,\varepsilon)$  and $\|\phi_t-g\|_{\mathcal B}  
<\varepsilon'$
for all $|t| < 2 \alpha_2$.
Now, we have the following claim.

\demo{\scshape Claim}
{\it The family $\phi_t$ satisfies the following property}\/{\rm :}\/
\begin{equation}
\label{fggfgthy}
  \left \{
\begin{array} {ll}
\phi_t \in V_\varepsilon^+(g), &  ~~~{\rm if}~~~
2 M \|v\|_{\mathcal B} < t < 2 \alpha_2  \\
\phi_t \in V_\varepsilon^-(g), &  ~~~{\rm if}~~~
- 2 \alpha_2 < t < - 2 M \|v\|_{\mathcal B} \ .
\end{array} \right .
\end{equation}
\Enddemo

To prove this claim,
let $C_5>0$ be as in Lemma \ref{L26} and let $p$ be such that
$C_5 \delta_0^{(p+1)m} <   \|\phi_t-f\|_{\mathcal C} \le  C_5  
\delta_0^{pm}$.
Set $k_0=k_0(g,\phi_t)$, $t_j=t_j(g,\phi_t)$ and $v_j=v_j(g,\phi_t)$
for all $0 \le j \le k_0$ as in Lemma \ref{l100}.
Set $s_j=s_j(g,f,\phi_t)= \delta_0^{jm} \theta_{f}(\phi_t-f)$.
Set  $l=l(g,p)$ as in Lemma \ref{L278}.
By Lemma  \ref{l11100} and Remark \ref{remmm}, there exist
$c_0 > 1$  and $\alpha_0>0$ sufficiently small
such that
if $\|g-g_0\|_{\mathcal B} < \alpha_0$
    then
\begin{equation}
\label{pppppp}
c_0^{-1} |t| \le  \left|\theta_{f}(\phi_t-f) \right| \le c_0 |t| \ ,
\end{equation}
(when $\|f-g\|_{\mathcal B} < \varepsilon'/2$
and  $\varepsilon'>0$ made  smaller if necessary).
Since $t \mathbf{u}_{g_0}= \phi_t-f + v$ and $2 M \|v\|_{\mathcal B}  <  
|t|$,
  by Lemma \ref{l11100} there is $c_1 > 1$ such that
\begin{equation}
\label{dffefgrth}  c_1^{-1}   \left(\delta_0^{pm} \right) ^{-1}
\le  |t|
  \le   c_1  \left(\delta_0^{pm} \right)^{-1} \ .
  \end{equation}
Hence, by \eqref{pppppp}, we obtain that
$$
  c_2^{-1} \left(\delta_0^{pm} \right) ^{-1}
\le \left|\theta_{f}(\phi_t-f) \right|
  \le  c_2 \left(\delta_0^{pm} \right)^{-1} \ .
$$
Thus,
\begin{equation}
\label{dddfgrth}
  c_3^{-1} \left(\delta_{lm}^{pm} \right) ^{-1}
\le  \left|s_l \right| \le c_3 \left(\delta_{lm}^{pm} \right)^{-1} \ .
  \end{equation}
   Recall  that $\delta_{jm}^{(j+1)m} > \beta > 2$.
By Lemma \ref{L278} we get
  $\left( \delta_{lm}^{pm} \right )^{-1}  \le \beta^{-(1- \mu_1)pm}$.
  Let us suppose from now on that $\theta_{f}(\phi_t-f)$ is positive and  
so $s_l>0$.
Hence, by  Lemma  \ref{L28}, by  \eqref{dddfgrth}
  and making $\alpha_1$ and $\alpha_2$ smaller if necessary
     (and so $p$ large enough),
we obtain that
\begin{eqnarray}
\label{ddwwwdfgrth1}
  t_l  & \ge &  s_l  - |t_l-s_l|   \ge
  s_l\left(1 - c_4 \left( \delta_{lm}^{pm} \right )^{-(\rho-1)} \right )  
  \\
&\ge & s_l\left(1 - c_4 \beta^{-(\rho-1)(1- \mu_1)pm} \right )  
\nonumber  \\
& > &  s_l/2  > 0  \ .\nonumber
    \end{eqnarray}
  Thus,   $t_l$ is positive and so it has the same sign as  
$\theta_{f}(\phi_t-f)$.
By induction on  $j=l,\dots, k_0(g_0,\phi_t)$  let us show that
$t_{j+1} \ge  t_j$ and so
that each $t_{j}$ is positive as well.
  By Lemma \ref{l100}, $\|v_j\|_{\mathcal B} \le \varepsilon 2^{-j}$
and $t_{j+1}= \hat \delta_{jm}^{(j+1)m} (t_j+ \sigma_{jm}(v_j))$.
  Since $t \mapsto \hat \delta_{jm}^{(j+1)m} (t)$ is $C^2$
as a map ${\mathbb R} \to {\mathcal C}$,
we obtain that
$\left| t_{j+1} - \delta_{jm}^{(j+1)m}  t_j \right| \le
c_5 \left (|t_j|^2 + \|v_j\|_{\mathcal B} \right)$.
Thus,
\begin{eqnarray}
\label{sdjkkj}
t_{j+1} & \ge & \beta t_j - c_5 |t_j|^2- c_6 \varepsilon 2^{-j}
\\
& \ge &   t_j (\beta - c_5 |t_j|) -c_6 \varepsilon 2^{-j} \ .\nonumber 
\end{eqnarray}
Let $\varepsilon_1>0$ be as given by Lemma \ref{l100}
and recall that $|t_j| < \varepsilon_1$ and $\varepsilon <  
\varepsilon_1$.
Since $\beta>1$ and by taking  $\varepsilon_1>0$
  sufficiently small, there is $\tau'>0$ with the property
that
$\beta - c_5 |t_j| > \beta - c_5 \varepsilon_1  > 1+ 2 \tau'$.
By  \eqref{dddfgrth} and \eqref{ddwwwdfgrth1}, we get
\begin{eqnarray}
\label{fdsdakh}
t_j (\beta - c_5 |t_j|-1- \tau' ) & > &  t_j \tau' \ge t_l \tau'
>   s_l  \tau' /2
> c_7 \left(\delta_{lm}^{pm} \right) ^{-1}  \ .
\end{eqnarray}
By Lemma \ref{L278},  
\begin{equation}
\label{fdewj}
     2^{-j} \le    2^{-l}
\le    \left(\delta_{lm}^{pm} \right) ^{-\rho} \ .
\end{equation}
Putting together \eqref{sdjkkj}, \eqref{fdsdakh} and
\eqref{fdewj}, we deduce that
\begin{equation}
\label{fdewj333}
t_{j+1}  \ge   (1+\tau') t_j +
  c_7   \left(\delta_{lm}^{pm} \right) ^{-1}
- c_6 \varepsilon \left(\delta_{lm}^{pm} \right) ^{-\rho} \ .
\end{equation}
Making  $\alpha_2$
  sufficiently small (and so $p$ large enough)
and recalling from  Lemma  \ref{L278} that  $l$ is a fraction of $p$,
we obtain  that
$
  c_7   \left(\delta_{lm}^{pm} \right) ^{-1}
- c_6 \varepsilon \left(\delta_{lm}^{pm} \right) ^{-\rho} \ge 0
$. Thus, by  \eqref{fdewj333}, 
  \begin{equation}
\label{khjert}
t_{j+1} > (1+ \tau') t_j
\end{equation}
which implies that $t_{j+1}$ has the same sign as $t_j$
and that $\phi_t \in V_\varepsilon^+(g)$.
If we suppose that  $\theta_{f_0}(\phi_t-f)$ is negative,
the proof that $t_l$ is negative and that
  $t_{j+1} < (1+ \tau') t_j$     follows in the same way
  for all $j=l,\dots, k_0(g_0,\phi_t)$
and so $\phi_t \in V_\varepsilon^-(g)$.
Therefore    \eqref{fggfgthy}
is satisfied and   the claim is proved.
\vskip6pt

Let us now prove the assertions of the lemma.
Take  $f=g_0$  and consider the family
$\phi_t=g_0+v+t \mathbf{u}_{g_0}$.
Since  $2M \|v\|_{\mathcal B} < 2M \alpha_1 < \alpha_2/4$,
the claim tell us that
\begin{equation*}
  \left \{
\begin{array} {ll}
\phi_t  \in V_{\varepsilon}^+ (g), &  ~~~{\rm if}~~~
\alpha_2 /4 \le t< 2\alpha_2  \\
\phi_t \in V_{\varepsilon}^- (g), &  ~~~{\rm if}~~~
-\alpha_2 <t \le -\alpha_2 /4 \ .
\end{array} \right .
\end{equation*}
Thus, by Lemma \ref{L22},  there is at least one value
$-\alpha_2/4 < \tau(g,v) < \alpha_2/4$ such that
$\phi_{\tau(g,v)}=g_0+v+ \tau(g,v)  \mathbf{u}_{g_0} \in  
V_{\varepsilon}^0 (g)$.

Next, take $f=\phi_{\tau(g,v)}$, and define a new family
$\psi_t=\phi_{\tau(g,v)}+t \mathbf{u}_{g_0}$.
Using the claim again, this time to the family $\psi_t$ (for which
$v=0$), we obtain that
\begin{equation*}
  \left \{
\begin{array} {ll}
\psi_t  \in V_{\varepsilon}^+ (g), &  ~~~{\rm if}~~~
0<t<\alpha_2/2  \\
\psi_t \in V_{\varepsilon}^- (g), &  ~~~{\rm if}~~~
- \alpha_2/2<t<0 \ .
\end{array} \right .
\end{equation*}
Therefore,  $\tau(g,v)$ is the only value  of $t \in \mathbb R$
between
$-2\alpha_2$ and $2\alpha_2$
such that
$\phi_t  \in V_{\varepsilon}^0 (g)$.
Since $\|\phi_{\tau(g,v)}-g\|_{\mathcal B} < \varepsilon'$ we
deduce from Lemma \ref{L22} that $\phi_{\tau(g,v)} \in W^s_\varepsilon  
(g)$.
This proves  assertions (i), (ii) and (iii).

Let us now prove assertion (iv).
Set  $k_0=k_0(g,\psi_t)$.
Using \eqref{khjert} and then \eqref{ddwwwdfgrth1},
we have
\begin{equation}
\label{geopop}
\left|t_{k_0}\right|  \ge  (1+ \tau')^{k_0-l} |t_l|
\ge    (1+ \tau')^{k_0-l} |s_l|/2   \ .
\end{equation}
Combining \eqref{dffefgrth} and \eqref{dddfgrth}, we see that
\begin{equation}
\label{gffweeffe}
   |s_l|   \ge c_3^{-1} \delta_{0}^{lm} \left( \delta_{0}^{pm}  
\right)^{-1}
\ge  c_8  \beta^l  |t|     \ .
\end{equation}
Taking $\tau ''=\min\{1+\tau', \beta\}$  and putting
\eqref{gffweeffe} back into \ref{geopop} we get
\begin{equation}
\label{fdgfcr}
\left|t_{k_0}\right| \ge  c_9 (\tau '')^{k_0} |t| \ .
\end{equation}
By Lemma \ref{l100}, there are $0< \varepsilon_1 < \varepsilon_0$
such that  $ \varepsilon_1 \le |t_{k_0}| \le \varepsilon_0$.
Thus, by \eqref{fdgfcr}, there is $\mu_0>0$ such that
$k_0 \ge  - \mu_0 \log (t)$.
By Lemma \ref{l100},  $\left|t_{k_0} \right|
\le  c_9 \beta^{k_0}  |t|$
and so there is  $\mu_1 \ge \mu_0$ such that
$k_0 \le  - \mu_1 \log (t)$,  which proves assertion (iv).
\hfq

\Subsec{Proof of the local stable manifold theorem}
It  will follow from  Theorem \ref{sddw}
in this section  that for every $g\in \mathbb K$
the local stable manifold at $g$ is a $C^1$ submanifold varying
continuously with $g$.
The proof of this theorem will use the following basic fact of  
calculus.

\begin{lemma}
\label{ghgtaw}
Let $X,Y$ be Banach spaces{\rm ,} let $x_0 \in X$ and consider a map
$\xi:B_X(x_0,\varepsilon) \to Y$ whose image in $Y$
falls within $B_Y(\xi(x_0),\varepsilon)$.
Suppose there is a bounded linear operator $L:X \to Y$ such that
for all $v \in X$ with $\|v\|_X \le \varepsilon${\rm ,}
\begin{equation}
\label{fewew}
\|\xi(x_0+v)-\xi(x_0)-L(v)\|_Y \le
c_0 \left(\|v\|_X+\|\xi(x_0+v)-\xi(x_0)\|_Y \right)^{1+\tau}
\end{equation}
where $c_0>0$ and $\tau>0$.
If $c_0(2\varepsilon)^\tau < 1$ then $\xi$ is differentiable at $x_0$
and $D\xi(x_0)=L$.
\end{lemma}

\Proof
Say $\|L(v)\|_Y \le a \|v\|_X$ for some $a>0$.
Noting that $$\|v\|_X +\|\xi(x_0+v)-\xi(x_0)\|_Y < 2 \varepsilon\ ,$$
we have from \eqref{fewew} that
$$ \|\xi(x_0+v)-\xi(x_0)\|_Y \le
\left( a+c_0(2\varepsilon)^\tau \right)
  \|v\|_X
+ c_0(2\varepsilon)^\tau \|\xi(x_0+v)-\xi(x_0)\|_Y
$$
whence
$$
\|\xi(x_0+v)-\xi(x_0)\|_Y \le
\frac{a+c_0(2\varepsilon)^\tau}{1-c_0(2\varepsilon)^\tau}
\|v\|_X
=c_1 \|v\|_X \ .
$$
Putting this back into the right-hand side of \eqref{fewew}
we get
$$
\|\xi(x_0+v)-\xi(x_0)-L(v)\|_Y \le
c_2 \left(\|v\|_X \right)^{1+\tau}
$$
and therefore $D \xi (x_0)$ exists and equals $L$.
\Endproof\vskip4pt

For every $g \in \mathbb K$ and $\alpha_1>0$,
let us consider  the following sets
\begin{eqnarray*}
E_{g,\alpha_1} & = & \left \{v \in {\rm ker}~ \theta_{g}:
  \|v \|_{\mathcal B }  < \alpha_1 \right \} \ , \\
F_{g} & = & \left \{g+t
\mathbf{u}_g: t \in {\mathbb R}  \right \} \ , \\
G_{g,\alpha_1}  & = &  \left \{g+v+ t
\mathbf{u}_g: v \in  E_{g,\alpha_1}~~~{\rm and}~~~
t \in {\mathbb R}  \right \} \ .
\end{eqnarray*}

\begin{theorem}
\label{sddw}
Set $0<\alpha_0<\alpha_1< \varepsilon$ and $\tau(g,v)$ as in Proposition~{\rm \ref{dfgfdger}.}
  For every  $g_0\in\mathbb{K}$ and
every $g \in \mathbb K$ with $\|g-g_0\|_{\mathcal B} < \alpha_0${\rm ,}
the map ${\xi}_{g}:E_{g_0,\alpha_1} \to F_{g_0}$ given by
${\xi}_{g}(v)=g_0+\tau(g,v){\mathbf u}_{g_0}$ is well-defined
and has the following properties\/{\rm :}\/
\begin{enumerate}
\item
The graph of $\xi_{g}$
is equal to  $W^s_\varepsilon(g) \cap G_{g_0,\alpha_1}${\rm ;}
\item
$\xi_{g}$ is  $C^1$ and varies continuously with $g$.
\end{enumerate}
\end{theorem}

\Proof  Set   $\alpha_1< \alpha_2< \varepsilon$ and
  $0<\mu_0 < \mu_1$   as in Proposition \ref{dfgfdger}.
By Proposition \ref{dfgfdger}, the map ${\xi}_{g}:E_{g_0,\alpha_1} \to  
F_{g_0}$
is well-defined and assertion (i) is satisfied.
Let $\hat{\xi}_{g}:E_{g_0,\alpha_1} \to \mathbb R$ be given by
$\hat{\xi}_{g}(v)=\tau(g,v)$ where $\tau(g,v)$ is
given by (i) in Proposition \ref{dfgfdger}.
To prove assertion (ii),  it is enough to show that
$\hat{\xi}_{g}:E_{g_0,\alpha_1} \to \mathbb R$  is
$C^1$ and varies continuously with $g$.
Let $v_1, v_2 \in E_{g_0,\alpha_1}$
  and set
\begin{equation}
\label{ggrrredd}
f_1=g_0+v_1+\hat{\xi}_{g}(v_1){\mathbf u}_{g_0}
\end{equation}
  and
$f_2=g_0+v_2+\hat{\xi}_{g}(v_2){\mathbf u}_{g_0}$. By Lemma
\ref{fffcccooo}, we get
\begin{multline*}
  \left| \theta_{f_1}(v_2-v_1) +
\theta_{f_1}(\mathbf{u}_{g_0})(\hat{\xi}_{g}(v_2)-\hat{\xi}_{g}(v_1) )  
\right |
   =  \left| \theta_{f_1}(f_2-f_1)\right |  \\
  \le  
c_0\left( \| v_2-v_1 \|_{\mathcal B} +
| \hat{\xi}_{g}(v_2)-\hat{\xi}_{g}(v_1) | \right
)^{1+\tau}\ .
\end{multline*}
By Lemma \ref{L23}, and taking
$\varepsilon>0$ sufficiently small,
we have that $\theta_{f_1}(\mathbf{u}_{g_0})$ is uniformly bounded away  
from $0$.
Therefore, by Lemma \ref{ghgtaw},
we deduce  that $\hat{\xi}_{g}$  is differentiable at
every $v_1$ with derivative given by
\begin{equation}
\label{ggrrredd3}
D \hat{\xi}_{g} (v_1)  =-(\theta_{f_1}(\mathbf{u}_{g_0}) )^{-1}
\theta_{f_1} \ .
\end{equation}
 From Lemma \ref{L23}, $\theta_{f_1}$ varies continuously with $f_1$
and so $D \hat{\xi}_{g}(v_1)$ also varies continuously in a  
neighborhood of
$v_1$. Hence,  $\hat{\xi}_{g}$ is a $C^1$ map.

Let us check that  $\hat{\xi}_{g}$ varies continuously with $g$ in the  
$C^1$ sense and,
more precisely, that the map $\mathbb{K} \cap \mathcal{B}(g_0,\alpha_0)  
\to
C^1(E_{g_0,\alpha_1},\mathbb{R})$ given by $g \mapsto \hat{\xi}_{g}$
is continuous. Taking into account that $D\hat{\xi}_{g}$
is given by \eqref{ggrrredd3} and that $f_1$ is given by  
\eqref{ggrrredd},
  and that by Lemma
\ref{L23} the map $f_1 \mapsto \theta_{f_1}$ is uniformly continuous
(as a map into $\mathcal{B}^*$), we see that it suffices to prove 
$g  \mapsto \hat{\xi}_{g}$ is continuous as a map into
$C^0(E_{g_0,\alpha_1},\mathbb{R})$.

To do this, let $v \in E_{g_0,\alpha_1}$ be such that
$g=g_0+v+\hat{\xi}_{g}(v) \mathbf{u}_{g_0}$, take $g_1 \in \mathbb{K}$
with $\|g_1-g_0\|_{\mathcal{B}} < \alpha_0$ and let
$w \in E_{g_0,\alpha_1}$ be such that
$g_1=g_0+w+ \hat{\xi}_{g_1}(w) \mathbf{u}_{g_0}$.
Now, we have the following claim.

\demo{\scshape Claim}
{\it There exist $c_1 > 0$ and $0< \gamma < 1$ such that}
\begin{equation}
\label{hfgedsg}
c_1^{-1} | \hat{\xi}_{g_1}(w)-\hat{\xi}_{g}(w)|^{1/\gamma}
\le  | \hat{\xi}_{g_1}(z)-\hat{\xi}_{g}(z)|
\le c_1 | \hat{\xi}_{g_1}(w)-\hat{\xi}_{g}(w)|^{\gamma} \ ,
\end{equation}
{\it for all} $z \in E_{g_0,\alpha_1}$.
\Enddemo

Let us assume this claim for a moment.
Its geometric meaning is that the distances between
corresponding points of the graphs of
$\hat{\xi}_{g}$ and $\hat{\xi}_{g_1}$ along the vertical fibers
$\left(\{z\} \times F_{g_0}\right )$ are uniform.
We want to control such distances in terms of
$\|g_1-g\|_{\mathcal{B}}$.
The above claim tell us that
it is enough to control
$|\hat{\xi}_{g_1}(w)-\hat{\xi}_{g}(w)|$.
Hence, write
\begin{eqnarray*}
g-g_1 & = & v-w +
\left(\hat{\xi}_{g}(v)-\hat{\xi}_{g_1}(w)\right) \mathbf{u}_{g_0} \\
& = & v-w +
\left(a+b \right) \mathbf{u}_{g_0}
\end{eqnarray*}
where $a=\hat{\xi}_{g}(w)-\hat{\xi}_{g_1}(w)$ and
$b=\hat{\xi}_{g}(v)-\hat{\xi}_{g}(w)$.
Since $\hat{\xi}_{g}$ is $C^1$,
we have $|b| \le c_2 \|v-w\|_\mathcal{B}$.
On the other hand, since
$\mathcal{B}={\rm Ker}\,\theta_{g_0} \bigoplus \mathbb{R}  
\mathbf{u}_{g_0}$
is a splitting into closed subspaces, there exists a constant $c_3>0$
such that
$$
\max \left\{\|v-w\|_{\mathcal{B}},|a+b| \right \} \le c_3
\|g-g_1\|_{\mathcal{B}} \ .
$$
But then
\begin{eqnarray*}
|a| & \le &
\|g-g_1\|_{\mathcal{B}}+ \|v-w\|_{\mathcal{B}} +|b| \\
& \le & (1+c_3+c_2 c_3) \|g-g_1\|_{\mathcal{B}} \ .
\end{eqnarray*}
Hence
$\left|\hat{\xi}_{g}(w)-\hat{\xi}_{g_1}(w)\right|
\le c_4 \|g-g_1\|_{\mathcal{B}} $,
and given the claim this proves that
$g \mapsto \hat{\xi}_{g}$ is indeed continuous.

\vskip6pt
Finally, let us prove the claim.
For each $z \in E_{g_0,\alpha_1}$,
let $h = g_0 + z +  \hat{\xi}_{g_1}(z) \mathbf{u}_{g_0}$. Set
\begin{alignat*}{5}
t_k'&=t_k(g,g_1), &\qquad  t_k''&=t_k(g,h),  \\
u_{k}'&=u_{T^{k m}(g)}(t_{k}'), &\qquad  u_{k}''&=u_{T^{k m}(g)}(t_{k}'') \ ,
\end{alignat*}
as given by Lemma  \ref{l100}, and set (also as in that lemma)
\begin{eqnarray}
\label{hfgrfx}
k_0' & = & k_0(g,g_1)= \min \{j:|t_j'| \ge \varepsilon_1 \} \\
k_0'' & = & k_0(g,h)= \min \{j:|t_j''| \ge \varepsilon_1 \} \ . \nonumber
\end{eqnarray}
Applying Lemma  \ref{l100}, we obtain, for all $k  \le \min  
\{k_0',k_0''\}$,
the estimates
\begin{eqnarray}
\label{htyrt}
\|T^{k m}(g_1)- u_{k }'\|_{\mathcal B}
  & \le & 2^{-k}  \|g_1-g\|_{\mathcal B} \\
\|T^{k m}(h)- u_{k }'' \|_{\mathcal B}
   & \le &  2^{-k}  \|h-g\|_{\mathcal B}  \ . \nonumber
\end{eqnarray}
Since $h \in W_\varepsilon^s(g_1)$,  we also have,
by  Lemma \ref{L22},
\begin{eqnarray}
\label{dffderh}
\|T^{k m}(h)-T^{k m}(g_1)\|_{\mathcal B}
  \le \varepsilon c_5 2^{-k}  \ .
\end{eqnarray}
Combining  \eqref{htyrt} and \eqref{dffderh}, we get
$$
\|u_{k }'- u_{k }''\|_{\mathcal B}    \le  c_6 2^{-k} \ .
$$
Hence, by Lemma \ref{l11100},  we get
\begin{equation}
\label{fghghfrte}
|t_k(g,g_1)-t_k(g,h)|  \le c_7 2^{-k} \ ,
\end{equation}
for all $k \le \min\{k_0',k_0''\}$.
Using \eqref{khjert} together with \eqref{fghghfrte},
we deduce that there exists a uniform
constant $c_8>0$ such that
\begin{equation}
\label{htghtr}
k_0'-c_8 \le k_0'' \le k_0'+c_8 \ .
\end{equation}
On the other hand, applying  (iv) in Proposition \ref{dfgfdger}, we  
also have
\begin{eqnarray}
\label{fdgrteer}\qquad\quad
- \mu_0 \log (|\tau(g_1,w)-\tau(g,w)|) & \le &  k_0' \le
- \mu_1 \log (|\tau(g_1,w)-\tau(g,w)|)   \\
- \mu_0 \log (|\tau(g_1,z)-\tau(g,z)|) & \le &  k_0'' \le
- \mu_1 \log (|\tau(g_1,z)-\tau(g,z)|)  \  .\nonumber
\end{eqnarray}
Combining  \eqref{htghtr} and  \eqref{fdgrteer}
and noting that
$$
\tau(g_1,w)= \hat{\xi}_{g_1}(w)\ , \
\tau(g_1,z)= \hat{\xi}_{g_1}(z) \ ,
\tau(g,w)= \hat{\xi}_{g}(w)\ \text{and} \
\tau(g,z)= \hat{\xi}_{g}(z) \ ,
$$
we get at last
  $$ c_9^{-1}
\left | \hat{\xi}_{g_1}(w)- \hat{\xi}_{g}(w) \right |^{\mu_0/\mu_1}
  \le
\left | \hat{\xi}_{g_1}(z)- \hat{\xi}_{g}(z)\right| \le
  c_9
\left | \hat{\xi}_{g_1}(w)- \hat{\xi}_{g}(w) \right|^{\mu_1/\mu_0} 
$$
for some $c_9 > 1$. This proves the claim with $\gamma = \mu_0/\mu_1 < 1$
and $c_1= c_9$.
\hfq

\begin{remark}
Note that by Proposition \ref{dfgfdger} there exists a uniform
$0<\tilde{\varepsilon} < \varepsilon$ such that
$W^s_{\tilde{\varepsilon}}(g) \subset G_{g_0,\alpha_1}$
for all $g_0 \in \mathbb K$ with $\|g_0-g\|_{\mathcal B} < \alpha_0$.
\end{remark}

\section{Smooth holonomies}
  \label{gholfff}

In the previous section we proved that a robust operator $T$ has
$C^1$ local stable manifolds through each point of its hyperbolic basic  
set
$\mathbb K$, and that such manifolds form a $C^0$ lamination
(near each point of $\mathbb K$).
A natural question that may be asked at this point is this:
how smooth is the holonomy of this lamination?
To answer this question
  we shall assume that there exists a homeomorphism
$H:\Theta^{\mathbb Z} \to \mathbb K$ of a finite-type
shift space onto $\mathbb K$ which conjugates the two-sided
full shift $\sigma:\Theta^{\mathbb Z} \to \mathbb K$
to our robust operator $T$ restricted to $\mathbb K$.
Under this topological assumption,
and an additional geometric assumption concerning the unstable
manifolds of points in the attractor
--both of which are satisfied by the renormalization operator--
we shall prove below
that the holonomies of the local stable laminations are $C^{1+\theta}$
for some $\theta > 0$.

\Subsec{Smooth holonomies for robust operators}
Let ${\mathbb K} \subset \mathcal O_A$ be a hyperbolic basic
set of a $C^2$ operator $T:\mathcal O_A \to \mathcal A$
which is  topologically conjugated to a two-sided shift
of finite type.
For $\varepsilon_0>0$ small enough and
for every $g \in \mathbb K$  let $t \to u_g(t)$
be a parametrization  of the local unstable manifold
$W^u_{\varepsilon_0} (g)$. Set
$$
{\mathbb K}_g =
{\mathbb K} \cap \overline{W^u_{\varepsilon_0}(g)} ~{\rm and}~
  K_g=u_g^{-1}({\mathbb K}_g) \  .
$$
Let
$$\Sigma_{\dots,\theta_{k-1},\theta_k}=
\left \{ \left ( \theta_j' \right )   \in \Theta^{\mathbb{Z}}:
  \theta_j'= \theta_j ~~~{\rm for}~~~{\rm all}  j \le k \right\} \ .
$$
If $H(\Sigma_{\dots,\theta_k}) \cap {\mathbb K}_g \ne \emptyset$
then denote by $\Delta_{\dots,\theta_k}(g)$ the smallest interval in  
$\mathbb R$
such that
$u_g(\Delta_{\dots,\theta_k}(g)) \supset H(\Sigma_{\dots,\theta_k})
\cap {\mathbb K}_g$.
Let $\mathbf{C}_k (g)$ be the set of all these intervals
$\Delta_{\dots,\theta_k}(g)$.

\begin{definition}
We say that the local unstable manifolds $W^u_{\varepsilon_0} (g)$
have {\emph geometry  bounded by $\alpha>0$} if
for every $g \in \mathbb K$,
$K_g$ has   geometry bounded by  $\alpha>0$ with respect
to the collection $(\mathbf{C}_k(g))_{k\ge 0}$
(in the sense of \S \ref{grtsfd}).
\end{definition}

Now suppose in addition that the operator $T$
is robust with  respect to the Banach spaces
$(\mathcal{B},\mathcal{C},\mathcal{D})$. By Theorem \ref{main3}, the  
local
stable manifolds of $T$ in $\mathcal{B}$ form a $C^0$ lamination
Let   $F:[-\mu_0, \mu_0]  \to {\mathcal B}(g,\varepsilon)$  be a
  $C^2$ curve transversal   to  the stable
lamination. Let $K_F$ be the set of all values
$r \in [-\mu_0, \mu_0]$ such that
$$f_r=F(r) \in \bigcup_{g_0 \in \mathbb{K}\cap W^u_{\varepsilon_0}(g_0)}
  W^s_{\varepsilon_0}(g_0)\ .$$
The {\it holonomy map}
$\phi_F: F(K_F) \to W^u_{\varepsilon_0}(g)$
associates to each $f_r$  the point $\phi_F(f_r)$
such that $f_r \in  W^s_{\varepsilon_0} (\phi_F(f_r))$.
In local coordinates, the holonomy map $\phi_F$
is given by   $\psi_F:K_F \to K_g$  where
$\psi_F(t)= u_g \circ \phi_F \circ F^{-1}$
and $K_F,K_g \subset \mathbb R$.
The $C^2$ curve $F:[-\mu_0, \mu_0]  \to {\mathcal B}(g,\varepsilon)$
is an {\it ordered transversal   to  the stable
foliation} if $F$ is transversal   to  the stable
foliation,
  $\phi_F: F(K_F) \to W^u_{\varepsilon_0}(g)$
extends to $F([-\mu_0, \mu_0])$
as a homeomorphism $\hat \phi_F$ over its image such that
$\phi_F( F(K_F))= \hat \phi_F( F(K_F)) \cap {\mathbb K}$.

We note that, by Remark \ref{remmm} and by  Theorem \ref{sddw},
  there is $\varepsilon_1 < \varepsilon_0$
small enough such that
a $C^2$ transversal to  $W^s_{\varepsilon_1}(g)$
in a point $f$ is an ordered transversal   to  the stable
foliation in a small neighborhood of $f$.

\begin{theorem}
\label{trfef}
Let ${\mathbb K} \subset \mathcal O_A$ be a hyperbolic basic
set of a $C^2$ operator $T:\mathcal O_A \to \mathcal A$
which is robust with  respect to  
$(\mathcal{B},\mathcal{C},\mathcal{D})$.
Suppose  that there exits $\varepsilon_0>0$
such  that the local unstable manifolds $W^u_{\varepsilon_0} (g)$
of $ g \in {\mathbb K}$ have bounded geometry.
There exists $0< \varepsilon < \varepsilon_0$ with the  property that
for every  $C^2$ ordered transversal $F:[-\mu_0, \mu_0]  \to {\mathcal  
B}(g,\varepsilon)$
     to  the stable
foliation  in $\mathcal B${\rm ,}
  the holonomy $\phi_F: F(K_F) \to W^u_{\varepsilon_0}(g)$
has a $C^{1+\theta}$ diffeomorphic  extension to
$F([-\mu_0, \mu_0])$
for some $\theta >0$.
\end{theorem}

\begin{example}
\label{fggjggh}
As we know from Theorem \ref{lyubhyp}, the
renormalization operator
$T=R^N:\mathbb{O}_\mathbb{A} \to \mathbb{A}$ is
hyperbolic over $\mathbb{K}$. As we shall see in Theorem \ref{rob},
  $T$ is robust
with respect to $(\mathbb{A}^r,\mathbb{A}^s,\mathbb{A}^0)$
provided  $s>s_0$ with $s_0$ sufficiently close to $2$ and  $r > s+1$  
is not
an integer.
By Theorem \ref{sullii},  there is  a  two-sided full shift
   topologically conjugated to
$T|\mathbb{K}$.
By Lemmas 9.3 and 9.6 respectively on pages 403 and 405
of Lyubich's paper \cite{lyubich}, there is $\alpha>0$
such that
the local unstable manifolds $W^u_{\varepsilon_0} (g)$
  have geometry  bounded by $\alpha$.
Hence the renormalization operator
$T$ satisfies the hypotheses of Theorem \ref{trfef}.
\end{example}

\def\cO{{\mathcal O}}
\def\varpsi{{\psi}}

In what follows, the notation
$A = \cO(B)$ means that
$\mu_1^{-1} B \le A \le \mu_1 B$ and
the notation $A=B(1 \pm \cO(C))$ means that
$B(1-\mu_2 C) \le A \le B(1+\mu_2 C)$ for some constants $\mu_1>1$ and  
$\mu_2 > 0$.

The proof of Theorem \ref{trfef} will be a consequence of
the following lemmas.

\begin{lemma}
\label{tr1}
For every  $C^2$ curve  $F:[-\mu_0, \mu_0]  \to {\mathcal  
B}(g,\varepsilon)$
  transversal   to  the stable
foliation and for
  all  $r,t \in [-\mu_0, \mu_0]$ such that
$r<t${\rm ,}
\begin{eqnarray}
\label{svsrfew}
\left \|f_t-f_r \right \|_{{\mathcal X}} & =  & \cO(|t-r|)\  ,\\
\left | \theta_{f_r}(f_t-f_r) \right |  & =  &   \cO(|t-r|)  \ , \nonumber
\end{eqnarray}
and for all  $s,r,t \in [-\mu_0, \mu_0]$ such that
$s<r<t${\rm ,}
\begin{eqnarray}
\label{atr1}
\frac{ \left \|f_t-f_r \right  \|_{{\mathcal X}}    }
{     \left  \|f_s-f_r \right \|_{{\mathcal X}}    }
& = &
\frac{ |t-r|}
{ |s-r|}   (1 \pm \cO(|t-s|))\ ,\\
   \frac{\left \|\theta_{f_r}(f_t-f_r) \right \|_{\mathcal X}}
{\left \|\theta_{f_r}(f_s-f_r) \right \|_{\mathcal X}}
& = & \frac{ |t-r|} { |s-r|} (1 \pm \cO(|t-s|))  \ ,  \nonumber 
\end{eqnarray}
where $\mathcal X  \in \{B,C,D\}$.
\end{lemma}

\Proof
By Lemma \ref{l11100}, there are $\nu_1,\nu_2 > 0$ such that  for all  
$r \in  K_f$,
$\left  \|{\mathbf u}_r  \right\|_{\mathcal X} >\nu_1$ and
$\left |\theta_{f_r} ({\mathbf u}_r) \right |> \nu_2$.
Since $F$ is $C^2$,  
\begin{eqnarray*}
f_t-f_r & =  & (t-r) {\mathbf u}_r  \pm \cO(|t-r|^2) , \\
\theta_{f_r}\left (f_t-f_r \right ) & =  &
(t-r) \theta_{f_r} \left ({\mathbf u}_r \right )   \pm
\cO(|t-r|^2) \ ,
\end{eqnarray*}
and so \eqref{svsrfew} follows.
Taking $s<r<t$, we get
\begin{eqnarray*}
\frac{\left  \|f_t-f_r \right \|_{\mathcal X}}{ \left \|f_s-f_r \right  
\|_{\mathcal X}}
& = &
\frac{\left \|{\mathbf u}_r \right\|_{\mathcal X}|t-r| (1 \pm  
\cO(|t-r|))}
{ \left  \|{\mathbf u}_r \right \|_{\mathcal X} |s-r| (1 \pm  
\cO(|s-r|))} \\
& = &
\frac{ |t-r|} { |s-r|} (1 \pm \cO(|t-s|))  \ .
\end{eqnarray*}
The second estimate in \eqref{atr1} is obtained in similar fashion.
\Endproof

In what follows, it will be more convenient to denote
$\phi_F(f_r)$ by $g_{\psi_F(r)}$. We will also work with a fixed  $0 <
\varepsilon <\varepsilon_0$ for which  Lemma \ref{L26} holds.

\begin{lemma}
\label{tr3}
Set  $l=l(g_{\psi_F(r)},p)$ as in Lemma {\rm \ref{L278}.}
Let    $F:[-\mu_0, \mu_0]  \to {\mathcal B}(g,\varepsilon)$
be a  $C^2$ curve transversal   to  the stable
foliation.  For all $p>0$   sufficiently large and
   all $s,r,t \in K_F$  such that
$$
|t-r|= \cO\left(\left(\delta_{0}^{pm}\right)^{-1}\right)\ \text{ and }\
|s-r|= \cO\left(\left(\delta_{0}^{pm}\right)^{-1}\right)\ ,
$$
we have
\begin{eqnarray}
\label{atfgfgr3}
\left \|T^{lm}(f_t)-T^{lm}(f_r) \right \|_{\mathcal C} & = &
\cO\left(\left(\delta_{lm}^{pm}\right)^{-1}\right)\ , \\
\|T^{lm}(f_s)-T^{lm}(f_r)\|_{\mathcal C} & = &
\cO\left(\left(\delta_{lm}^{pm}\right)^{-1}\right)\ , \nonumber\\
  \frac{\left \|T^{lm}(f_t)-T^{lm}(f_r) \right \|_{\mathcal C}}
{\left \|T^{lm}(f_s)-T^{lm}(f_r) \right \|_{\mathcal C}}
& = & \frac{|t-r|}{|s-r|}
\left(1 \pm  
\cO\left(\left(\delta_{lm}^{pm}\right)^{-(\rho-1)}\right)\right)   \ .\nonumber
\end{eqnarray}
\end{lemma}

\Proof
Using Lemma \ref{l11100} and \eqref{svsrfew}, we get
$$
\left| \theta_{f_r}(f_t-f_r)  \right|    =
\cO\left(\left(\delta_{0}^{pm}\right)^{-1}\right)\ \text{ and }\
\left| \theta_{f_r}(f_s-f_r)  \right|   =
\cO\left(\left(\delta_{0}^{pm}\right)^{-1}\right)  \ .
$$
Thus, taking  $p$ sufficiently large and using
Lemmas \ref{L27} and \ref{L278}
we deduce that
\begin{eqnarray}
\label{fvfdgr}
\left \|T^{lm}(f_t)-T^{lm}(f_r) \right \|_{\mathcal C}
& =  &  \left|
\delta_{0}^{lm}  \theta_{f_r}(f_t-f_r)  \right|
\pm  \cO  \left( \left( \delta_{lm}^{pm}  \right)^{-\rho}  \right)  
\\
&= &
   \cO  \left( \left( \delta_{lm}^{pm}  \right)^{-1}  \right)
\pm  \cO  \left(   \left(   \delta_{lm}^{pm}  \right)^{-\rho}  \right)   
\nonumber \\
&= &  \cO  \left( \left( \delta_{lm}^{pm}  \right)^{-1} \right) \ .\nonumber
\end{eqnarray}
Similarly,
$ \left \|T^{lm}(f_s)-T^{lm}(f_r) \right \|_{\mathcal C}  =
\cO\left(\left(\delta_{lm}^{pm}\right)^{-1}\right) $.
This proves the first two inequalities
in \eqref{atfgfgr3}. Combining \eqref{svsrfew}
with \eqref{fvfdgr},  we see that
  \begin{eqnarray*}
  \frac{\left  \|T^{lm}(f_t)-T^{lm}(f_r)\right \|_{\mathcal C}}
{\left  \|T^{lm}(f_s)-T^{lm}(f_r)\right \|_{\mathcal C}}
& = &
  \frac
{ \left|
\delta_{0}^{lm}  \theta_{f_r}(f_t-f_r)  \right|
\pm  \cO  \left(  \left( \delta_{lm}^{pm}  \right)^{-\rho} \right)}
{ \left|
\delta_{0}^{lm}  \theta_{f_r}(f_s-f_r)  \right|
\pm  \cO  \left(  \left( \delta_{lm}^{pm}  \right)^{-\rho}\right) } \\
& = &
  \frac
{ \left|
   \theta_{f_r}(f_t-f_r)  \right|
\left (1
\pm \cO  \left(\left( \delta_{lm}^{pm}  \right)^{-(\rho-1)} \right)
  \right) }
{ \left|
   \theta_{f_r}(f_s-f_r)  \right|
  \left (1
\pm \cO  \left(\left( \delta_{lm}^{pm}  \right)^{-(\rho-1)} \right)   
\right) }\ .
\end{eqnarray*}
Therefore, by Lemma \ref{tr1}, we
get
$$
  \frac{\left \|T^{lm}(f_t)-T^{lm}(f_r)  \right \|_{\mathcal C}}
{\left \|T^{lm}(f_s)-T^{lm}(f_r)  \right \|_{\mathcal C}}
   =     \frac
{  |t-r|  }
{  |s-r| }
  \left (1  \pm \cO
\left( \left( \delta_{lm}^{pm}  \right)^{-(\rho-1)}
  \right)\right)  \ ,
$$
and this proves the last inequality in \eqref{atfgfgr3}.
\Endproof

\begin{lemma}
\label{tr5}
Set  $l=l\left (g_{\psi_F(r)},p \right )$ as in Lemma {\rm \ref{L278}.}
  Let  $F:[-\mu_0, \mu_0]  \to {\mathcal B}(g,\varepsilon)$
be a $C^2$ curve  transversal   to  the stable
foliation. For every $s \in K_F$ and
$s'=\psi_F(s) \in K_g${\rm ,}
$$
\left \|T^{lm}(f_s)-T^{lm}(g_{s'}) \right \|_{\mathcal C}
\le  \cO\left(\left(\delta_{lm}^{pm}\right)^{-\rho}\right) \ .
$$
Furthermore{\rm ,} for all $p$ large enough and all
    $s',r',t' \in K_g$
such that
$s=\psi_F^{-1}(s')${\rm ,} $r=\psi_F^{-1}(r')${\rm ,} $t=\psi_F^{-1}(t') \in K_f${\rm ,}
$$|t'-r'|  = \cO\left(\left(\delta_{0}^{pm}\right)^{-1}\right)\ \text{  
and }\
|s'-r'| = \cO\left(\left(\delta_{0}^{pm}\right)^{-1}\right)\ ,
$$
we have
\begin{eqnarray*}
   \frac{\left \|T^{lm}(f_t)-T^{lm}(f_r)\right \|_{\mathcal C}}
{\left \|T^{lm}(f_s)-T^{lm}(f_r)\right \|_{\mathcal C}}
  & =  &
  \frac{\left \|T^{lm}(g_{t'})-T^{lm}(g_{r'})\right \|_{\mathcal C}}
{\left \|T^{lm}(g_{s'})-T^{lm}(g_{r'})\right \|_{\mathcal C}}
\left(1  \pm  
\cO\left(\left(\delta_{lm}^{pm}\right)^{-(\rho-1)}\right)\right)  \ .
\end{eqnarray*}
\end{lemma}

\Proof
By Lemmas \ref{L22} and \ref{L278},  
$$
\left \|T^{lm}(f_s)-T^{lm}(g_{s'}) \right  \|_{\mathcal C}
\le C_3 \varepsilon 2^{-l}  \le
\cO\left(\left(\delta_{lm}^{pm}\right)^{-\rho}\right) \ .
$$
Thus, applying Lemma \ref{tr3} to the transversal
given by the local unstable manifold $\{g_t\}$ we get
\begin{eqnarray*}
  \frac{\left \|T^{lm}(f_t)-T^{lm}(f_r) \right \|_{\mathcal C}}
{\left \|T^{lm}(f_s)-T^{lm}(f_r) \right \|_{\mathcal C}}
& = &
  \frac{\left \|T^{lm}(g_{t'})-T^{lm}(g_{r'}) \right \|_{\mathcal C}
  \left(1 \pm
  \cO\left( \left(\delta_{lm}^{pm}\right)^{-(\rho-1)}\right)    \right )  
}
{\left \|T^{lm}(g_{s'})-T^{lm}(g_{r'}) \right \|_{\mathcal C} \left(1  
\pm
  \cO\left( \left(\delta_{lm}^{pm}\right)^{-(\rho-1)}\right)    \right )  
} \\[5pt]
& = &
\frac{\left  \|T^{lm}(g_{t'})-T^{lm}(g_{r'})\right \|_{\mathcal C} }
{\left  \|T^{lm}(g_{s'})-T^{lm}(g_{r'}) \right \|_{\mathcal C} }
\left(1 \pm
  \cO\left( \left(\delta_{lm}^{pm}\right)^{-(\rho-1)}\right)    \right )  
  \ .
\end{eqnarray*}
\vglue-19pt
\Endproof

\vglue12pt

 {\it Proof of Theorem} \ref{trfef}.
Let   $p>0$ be so large   that Lemmas \ref{tr3} and \ref{tr5}
are satisfied and  let $t,s,r,t', s',  r'$  be as in Lemma \ref{tr5}.
First, a claim.

\demo{\sc Claim} {\it We have}
$$
|t-r| = \cO \left ( \left ( \delta_0^{ p m}
\right )^{-1} \right ) \ \text{ and }\
|s-r| = \cO \left ( \left ( \delta_0^{p m}
\right )^{-1} \right ) \ .
$$
\smallskip

Assuming this claim for a moment, we  finish the proof of
  Theorem \ref{trfef} as follows.
Set  $l=l(g_{r'},p)$ as   in Lemma \ref{L278}.
By Lemmas \ref{L278} and \ref{tr1}, there is $0<\tau_1<1$ such that
$\delta_{lm}^{pm} \le \cO \left(|t'-r'|^{\tau_1} \right)$.
Therefore, by Lemmas \ref{tr3} and \ref{tr5} we deduce that
\begin{eqnarray}
\label{dfferd}
\frac{|t-r|}{|s-r|}
& = &
\frac{\left \|T^{lm} (f_t) - T^{lm} (f_r) \right \|_{\mathcal C} }
{\left \|T^{lm} (f_s) - T^{lm} (f_r) \right \|_{\mathcal C} }
  \left(1\pm \cO \left(\left( \delta_{lm}^{pm} \right)^{-(\rho-1)}
\right)\right)  \\
& = &
\frac{\left  \|T^{lm} (g_{t'}) - T^{lm} (g_{r'}) \right \|_{\mathcal C}  
}
{\left  \|T^{lm} (g_{s'}) - T^{lm} (g_{r'}) \right \|_{\mathcal C} }
   \left(1\pm \cO \left(\left( \delta_{lm}^{pm} \right)^{-(\rho-1)}
  \right) \right) \nonumber \\
& = &
\frac{|t' - r'|}
{|s' - r'|}  \left(1 \pm  \cO \left( |t' - r'|^{-(\rho-1)\tau_1}  
\right) \right)\ .\nonumber
\end{eqnarray}
Since $\mathcal{K}_g$ has bounded geometry
  and by Theorem 9.5 on page 549 of \cite {MS},
  the inequalities \eqref{dfferd}
imply that the  map $\psi_F$ has a $C^{1+\theta}$
diffeomorphic extension to ${\mathbb R}$ for some $0<\theta <1$.

Let us now prove the claim.
Let $\hat p$ be such that
$|t-r| = \cO \left ( \left ( \delta_0^{ \hat p m}
\right )^{-1} \right )$.
All that we have to show is that
\begin{equation}
\label{ggewew}
\left |\hat p - p \right  | \le \cO(1)\ .
\end{equation}
Set  $\hat l= \hat l(g_{r'},\hat p )$
  as   in Lemma \ref{L278}.
By Lemmas \ref{L27} and  \ref{L28}, and the estimates in  
\eqref{svsrfew},
for $k \le \min\{l,\hat l\}$,
\begin{eqnarray}
\label{rrefsddc} \qquad\quad
\left  \|T^{km} (g_{t'}) - T^{km} (g_{r'})\right \|_{\mathcal C} &  =  &
\left |\delta_0^{km} \theta_{g_r} (g_{t'}-g_{r'}) \right |
\pm \cO \left(\left( \delta_{km}^{p m} \right)^{-\rho}
+ 2^{-k} \right)\ , \\
\left \|T^{km} (f_t) - T^{km} (f_r)\right \|_{\mathcal C}  &  =  &
\left |\delta_0^{km} \theta_{f_r} (f_t-f_r) \right |
\pm \cO \left(\left( \delta_{km}^{\hat p m} \right)^{-\rho}
+ 2^{-k} \right)  \ . \nonumber
\end{eqnarray}
By Lemma \ref{L22},
we get (for all $k \le \min \left \{l,\hat l \right \}$)
\begin{eqnarray}
\label{fwefca}
  \left \|T^{km} (f_t) - T^{km} (f_r) \right \|_{\mathcal C} & =  &
  \left \|T^{km} (g_{t'}) - T^{km} (g_{r'}) \right \|_{\mathcal C}
\pm  \cO \left(2^{-k} \right) \ .
\end{eqnarray}
Let us consider separately the case (i) where $p \le \hat p$
and the case (ii) where $p \ge \hat p$.

\demo{Case {\rm (i)}} Here $l \le \hat l$ and by
Lemma \ref{L278} we get
\begin{equation}
  \label{ekuhgg}
\left( \delta_{lm}^{\hat p m} \right)^{-1} \le
\left( \delta_{lm}^{pm} \right)^{-1}
2^{-l} \le \cO
  \left(\left( \delta_{lm}^{pm} \right)^{-\rho}\right) \ .
\end{equation}
By Lemma \ref{tr3} applied to the transversal given by the local
unstable manifold $\{g_t\}$, 
\begin{equation}
  \label{eghghrg}
  \left  \|T^{lm} (g_{t'}) - T^{lm} (g_{r'}) \right \|_{\mathcal C}
= \cO
  \left(\left( \delta_{lm}^{pm} \right)^{-1}\right) \ .
  \end{equation}
  On the other hand, by \eqref{fwefca}, 
$$
\left \|T^{lm} (f_t) - T^{lm} (f_r)\right \|_{\mathcal C}  =
\left \|T^{lm} (g_{t'}) - T^{lm} (g_{r'})\right \|_{\mathcal C}  \pm  
\cO \left(2^{-l}\right) \
.
$$
But $2^{-l}$ is much smaller than
$ ( \delta_{lm}^{pm})^{-1}$.
  Hence by   \eqref{eghghrg},
  \begin{equation}
  \label{efefefwwe}
  \left \|T^{lm} (f_t) - T^{lm} (f_r) \right \|_{\mathcal C}  =
\cO
  \left(\left( \delta_{lm}^{pm} \right)^{-1}\right) \ .
  \end{equation}
Thus, by \eqref{rrefsddc} and \eqref{ekuhgg}, we deduce that
\begin{eqnarray*}
\left |\delta_0^{lm} \theta_{f_r} (f_t-f_r) \right | &  = &
\left \|T^{lm} (f_t) - T^{lm} (f_r)\right \|_{\mathcal C}
\pm \cO
  \left(\left( \delta_{lm}^{\hat p m} \right)^{-\rho}\right)
\pm \cO  \left( 2^{-l}\right) \\
  &  = &
\cO
  \left(\left( \delta_{lm}^{p m} \right)^{-1}\right)
\pm \cO
  \left(\left( \delta_{lm}^{\hat p m} \right)^{-\rho}\right)  \ .
\end{eqnarray*}
Since $\left( \delta_{lm}^{p m} \right)^{-1}$
is much larger than $\left( \delta_{lm}^{\hat p m} \right)^{-\rho}$, it
follows that
$$
\left|\delta_0^{lm} \theta_{f_r} (f_t-f_r) \right | =
  \cO
  \left(\left( \delta_{lm}^{pm} \right)^{-1}\right)  \ .
$$
This shows that
$$
\left |\theta_{f_r} (f_t-f_r) \right | =
\cO
  \left(\left( \delta_{0}^{pm} \right)^{-1}\right) \ .
  $$
Therefore, by Lemma \ref{tr1},  we get
  $\left \|f_t-f_r\right \|_{\mathcal C} =
\cO
  \left(\left( \delta_{0}^{pm} \right)^{-1}\right)$,
  which in turn implies \eqref{ggewew}.

\demo{Case {\rm (ii)}} Here $\hat{l} \le   l$,  and
Lemma \ref{L278} tells us that
\begin{equation}
  \label{yegegrrh}
\left( \delta_{\hat{l}m}^{pm} \right)^{-1} \le
\left( \delta_{\hat{l}m}^{\hat p m} \right)^{-1}
2^{-\hat{l}} \le \cO
  \left(\left( \delta_{\hat{l}m}^{\hat p m} \right)^{-\rho}\right) \ .
\end{equation}
Applying Lemma \ref{tr3}, we get
\begin{equation}
  \label{yhygh}
   \left \|T^{\hat{l}m} (f_t) - T^{\hat{l}m} (f_r)\right \|_{\mathcal C}
= \cO
   \left(\left( \delta_{\hat{l}m}^{\hat p m} \right)^{-1}\right) \ .
  \end{equation}
On the other hand, by \eqref{yegegrrh} and \eqref{rrefsddc}, 
$$
  \left  \|T^{\hat{l}m} (f_t) - T^{\hat{l}m} (f_r)\right \|_{\mathcal C}
   =   \cO
  \left( \left |\delta_0^{\hat{l} m} \theta_{f_r} (f_t-f_r)
  \right| \right ) \pm
\cO
  \left(\left( \delta_{\hat{l}m}^{\hat p m} \right)^{-\rho}\right)
\ .
$$
But  $ \left(\delta_{\hat{l}m}^{\hat p m} \right)^{-\rho}$
is much smaller than
$ \left(\delta_{\hat{l}m}^{\hat p m} \right)^{-1}$.
Hence, by  \eqref{yhygh}, we get
\begin{equation}
  \label{rtregh}
\left \|T^{\hat{l}m} (f_t) - T^{\hat{l}m} (f_r)\right \|_{\mathcal C}
   =   \cO
  \left( \left |\delta_0^{\hat{l} m} \theta_{f_r} (f_t-f_r)
  \right| \right ) \ .
\end{equation}
By   \eqref{fwefca}, we have also
$$
  \left \|T^{\hat{l}m} (f_t) - T^{\hat{l}m} (f_r)\right \|_{\mathcal C}
   =
\left \|T^{\hat{l}m} (g_{t'}) - T^{\hat{l}m} (g_{r'})\right  
\|_{\mathcal C}
\pm \cO  \left(2^{-\hat{l}}\right) \ .
$$
But $2^{-\hat{l}} \le \cO
\left( \left(\delta_{\hat{l}m}^{\hat p m} \right)^{-\rho} \right)$.
Hence, again by \eqref{yhygh}, we deduce that
\begin{equation}
  \label{rffwfh}
\left \|T^{\hat{l}m} (f_t) - T^{\hat{l}m} (f_r) \right \|_{\mathcal C}
   =   \cO
  \left( \left
\|T^{\hat{l}m} (g_{t'}) - T^{\hat{l}m} (g_{r'})\right \|_{\mathcal C}    
\right ) \ .
\end{equation}
By Lemma \ref{tr1} and  \eqref{rrefsddc}, 
$$
\left \|T^{\hat{l}m} (g_{t'}) - T^{\hat{l}m} (g_{r'})\right  
\|_{\mathcal C}
=  \cO \left( \left ( \delta_{\hat{l}m}^{pm}  \right)^{-1} \right) +
  \cO \left( \left ( \delta_{\hat{l}m}^{\hat{p}m}  \right)^{-\rho}  
\right)
\ .
$$
Therefore, by  \eqref{yhygh} and \eqref{rffwfh},
$$
\left \|T^{\hat{l}m} (g_{t'}) - T^{\hat{l}m} (g_{r'})\right  
\|_{\mathcal C}
=  \cO \left( \left ( \delta_{\hat{l}m}^{pm}  \right)^{-1} \right) \ .
$$
But then, using \eqref{rtregh} and \eqref{rffwfh} once more, we deduce  
that
$$
\left |\theta_{f_r} (f_t-f_r) \right | =
\cO
  \left(\left( \delta_{0}^{pm} \right)^{-1}\right) \ .
  $$
Therefore, by Lemma \ref{tr1},  we get
  $\left \|f_t-f_r \right \|_{\mathcal C} =
\cO
  \left(\left( \delta_{0}^{pm} \right)^{-1}\right)$
  which in turn implies \eqref{ggewew}.
\hfill\qed

\section{The renormalization operator is robust}
\label{sec:robust}

 From the very beginning, our main goal is to show that
the renormalization operator is ``hyperbolic'' in $\mathbb{U}^r$,
provided $r$ is sufficiently large.
More precisely, we want to establish Theorem \ref{main} and Corollary
\ref{mainmain}. We have already at our disposal an abstract theorem
(Theorem \ref{main3}) showing that any robust operator is indeed
``hyperbolic''. Hence, our work has been essentially reduced to showing  
that the
renormalization operator $T$, or any one of its powers, is robust (see  
Theorem
\ref{rob} below). The proofs of Theorem \ref{main} and  
Corollary
\ref{mainmain}  will be given in \S \ref{klrlrlll}.

We emphasize the important role played by the geometric estimates of
\S \ref{sec:bdedgeom} in the verification of properties
$\mathbf{B5}$ and $\mathbf{B6}$ of a robust operator (Definition
\ref{def4}) for an iterate of the renormalization operator (see \S
\ref{sec:frechet}). Properties $\mathbf{B2}$, $\mathbf{B3}$, and
$\mathbf{B4}$ are relatively straightforward consequences of the
properties of the composition operator studied  in \S
\ref{sec:composition} and are proved in
  \S  \ref{sec:axiom5}  and \S    \ref{sec:gateaux} .

In this section, we shall prove the following result (see \S  
\ref{sec:rob}).

\begin{theorem}
\label{rob}
Let $T:\mathbb{O}\to \mathbb{A}$ be the renormalization operator given
by Theorem {\rm \ref{lyubhyp}}. If $s>s_0$ with $s_0<2$ sufficiently close to  
$2$
and $r>s+1$ not an integer{\rm ,}  then  $T$ is
a robust operator with respect to
$(\mathbb{A}^r,\mathbb{A}^s,\mathbb{A}^0)$.
\end{theorem}

We shall present in the sequel complete proofs of all the estimates  
that are
necessary for establishing the above result, carefully checking all the
properties of robustness along the way.

In our estimates we will often concern ourselves with a power $T^m$ of  
$T$.
For each $m \ge 1$, let $\mathbb{O}_m \subseteq \mathbb{O}$
be the  (open)  set of those $f$'s which are $mN$ times renormalizable.
Then $T^m$ is well-defined in $\mathbb{O}_m$ and we can write
$$
T^m(f)= \frac{1}{\lambda_f} \cdot f^p \circ \Lambda_f \ ,
$$
where $p=p(f,mN)$, $\lambda_f=f^p(0)$, and
$\Lambda_f:x \mapsto \lambda_f x$ is the linear scaling.
Note that $p$ (and hence $\lambda_f$ and $\Lambda_f$)
depends on $m$, but if $m$ is held  fixed then $p$ is a locally
constant function of $f \in \mathbb{O}_m$.
To keep track of the dependence of constants on $m$, we shall denote
by $K$ those constants that may depend on $m$,
and by $c$ those that are independent of $m$.

Likewise, we define $\mathbb{O}^r_m$ to be an open set
in $\mathbb{U}^r$  containing $\mathbb{K}$, all of whose elements are  
$mN$
times renormalizable, so that $T^m=R^{mN}:\mathbb{O}^r_m \to  
\mathbb{U}^r$ is
well-defined.

\Subsec{A closer look at composition}
\label{sec:composition}
 From a differentiable viewpoint, composition is a notoriously
ill-behaved operation. Such bad behavior is the source of most
technical difficulties arising in this work. Fortunately, some
positive results lie at hand. For example, it is well-known that if
$r$ is a positive integer then composition, viewed as an operator from
$C^r\times C^{r-1}$ into $C^{r-1}$, is a $C^1$ map (see \cite{franks}).  
We
shall need not only this result but also a less well-known  
generalization of
it for H\"older spaces: if $r-1>s\ge 1$ are {\it real} numbers then
composition, as an operator from $C^r\times C^s$ into $C^s$, is a
$C^1$ map (we can say a little bit more -- see Proposition \ref{comp}
below). For related results on the smoothness of composition,
see \cite{deLlave}.

Before we can prove this fact, some auxiliary results are in order.
In what follows, our definition of the $C^r$ norm of $\varphi\in
C^r(I)$ is this: for $r=k+\alpha$ with $k\in \mathbb{N}$
and $0\le \alpha <1$  we write
$$
\|\varphi\|_{C^r}=\max\{\|\varphi\|_0,\|\varphi'\|_0, \dots
,\|\varphi^{(k)}\|_0; \|\varphi^{(k)}\|_\alpha\}
\ .
$$
For $r=k+\text{Lip}$ we define $\|\varphi\|_{C^r}$ as above with  
$\alpha=1$.
This norm is equivalent to the one introduced earlier in  
\ref{sec:holder}, and
has the advantage that $\|\varphi\|_{C^r}=\max\{\|\varphi\|_{C^0},
\|\varphi'\|_{C^{r-1}}\}$ whenever $r\ge 1$. This allows us to prove
certain estimates by induction on $k$, which will be very useful
later.

\begin{lemma} \label{conta0} Given $0 \le \alpha < 1$ and $0 \le
\epsil \le  1- \alpha${\rm ,}
let $w \in C^{\alpha+\epsil}(I)${\rm ,}
  $\varphi, \varpsi \in C^{1}(I,I)$ and    
$||\varpsi-\varphi||_{C^{1}}\leq 1$.
\begin{enumerate}
\item If $\epsil > 0$ then there exists
$K=K \left(\|\varpsi\|_{C^1} \right)>0$
  such that
$$
\left\|w\circ  \varphi -w\circ \varpsi \right\|_{C^\alpha} \leq
K  \|w\|_{C^{\alpha+\epsil}}
   \|\varphi-\varpsi\|_{C^1}^\epsil  \ .
$$
\item  If $\epsil = 0$ then there exist
$c>0$ and $K=K \left(\|\varpsi\|_{C^1}  \right)>0$
     such that
\begin{eqnarray*}
\|w\circ  \varphi -w\circ \varpsi \|_{C^\alpha} & \le  &
  c     \|w\|_{C^\alpha} \|\varpsi'\|_{C^0}^\alpha \\
&&+ K  \|w\|_{C^{\alpha}}
   \|\varphi-\varpsi\|_{C^1}^\alpha  \ .
\end{eqnarray*}
\end{enumerate}
\end{lemma}

\Proof
Let us start proving part (i) of this lemma.
By the mean value theorem, we obtain
\begin{equation}
\label{aqaz5}
\|w\circ  \varphi -w\circ \varpsi\|_{C^0} \le
\|w\|_{C^{\alpha+\epsil}} \|\varphi-\varpsi\|_{C^0}^{\alpha+\epsil}.
\end{equation}
  If  $|y-x| \le \| \varphi-\varpsi \|_{C^\alpha}$ then
\begin{eqnarray*}
|w \circ \varphi (y) - w \circ \varphi (x)|  & \le &
c_0 \|w\|_{C^{\alpha+\epsil}}  \| \varphi\|_{C^1}^{\alpha + \epsil}
  \| \varphi-\varpsi \|_{C^\alpha}^\epsil |y-x|^\alpha \ ,\\
| w \circ \varpsi (y) - w \circ \varpsi (x) |  & \le &
c_1 \|w\|_{C^{\alpha+\epsil}}  \| \varpsi\|_{C^1}^{\alpha + \epsil}
  \| \varphi-\varpsi \|_{C^\alpha}^\epsil |y-x|^\alpha \ .
\end{eqnarray*}
If   $|y-x| >  \| \varphi-\varpsi \|_{C^\alpha}$, by \eqref{aqaz5} then
$$
\| w \circ \varpsi - w \circ \varphi \|_{C^0}
   \le   c_2 \|w\|_{C^{\alpha+\epsil}}  \| \varphi-\varpsi  
\|_{C^\alpha}^\epsil
  |y-x|^\alpha \ ,
$$
  which ends the proof of part (i) of this lemma.

Let us  prove part (ii) of this lemma.
By the mean value theorem, we obtain
\begin{equation}
\label{aqaz7}
\|w\circ  \varphi -w\circ \varpsi\|_{C^0} \le
\|w\|_{C^{\alpha}} \|\varphi-\varpsi\|_{C^0}^{\alpha}.
\end{equation}
Furthermore,
\begin{eqnarray*}
| w \circ \varphi (y) - w \circ \varphi (x) | & \le &
c_3  \|w\|_{C^{\alpha}}  \| \varphi'\|_{C^0}^{\alpha}  |y-x|^\alpha\ ,  \\
| w \circ \varpsi (y) - w \circ \varpsi (x) | & \le &
c_4 \|w\|_{C^{\alpha}}  \| \varpsi'\|_{C^0}^{\alpha}  |y-x|^\alpha \ ,
\end{eqnarray*}
and so
\begin{equation}
\label{aqaz9}
\|w\circ  \varphi -w\circ \varpsi\|_{C^\alpha} \le c_5  
\|w\|_{C^{\alpha}}
( \|\varphi'-\varpsi'\|_{C^0} + \| \varpsi'\|_{C^0} )^{\alpha}
\end{equation}
  which ends the proof of part (ii) of this lemma.
\Endproof

We shall need also some estimates on polynomial operators
coming from simple algebraic considerations.
For every polynomial $P$ of degree $d$ in $n$ variables
$x_1,x_2,\dots,x_n$ over $\mathbb{R}$, define
$\nu(P)$ as the sum of the absolute values of the
coefficients of $P$.
This is a well-known valuation in the ring  
$\mathbb{R}[x_1,x_2,\dots,x_n]$,
but all that really matters to us is that $\nu(P+Q) \le \nu(P)+\nu(Q)$
(sub-additivity), and that $\nu (\partial_{x_i} P ) \le d \nu (P)$.

\begin{lemma}\label{fr4}
Let $P \in \mathbb{R}[x_1,x_2,\dots,x_n]$ be a polynomial of degree
$d${\rm ,} and let $\phi_1, \phi_2,\dots,\phi_n \in C^s(I)$. Then{\rm ,} 
$$
\|P(\phi_1,\phi_2,\dots,\phi_n)\|_{C^s} \le \nu(P) 2^{sd}M^d \ ,
$$
where $M= \max\{1,  
\|\phi_1\|_{C^s},\|\phi_2\|_{C^s},\dots,\|\phi_n\|_{C^s}\}$.
Moreover{\rm ,}
if  $\psi_1,\psi_2,\dots,\psi_n \in C^s(I)$
also satisfy $\|\psi_i\|_{C^s} \le M$ for all $1 \le i\le n${\rm ,}
then
$$
\|P(\phi_1,\phi_2,\dots,\phi_n)-P(\psi_1,\psi_2,\dots,\psi_n)\|_{C^s}
\le d \nu(P) 2^{sd}M^{d-1} \sum_{i=1}^n \|\phi_i-\psi_i\|_{C^s} \ .
$$
\end{lemma}

\Proof
The first inequality is immediate from the definition of $\nu(P)$.
To prove the second, note that
$P:\left(C^s(I)\right)^n \to C^s(I)$ is a $C^1$ map
(norm of the sum in the domain of $P$).
Using the mean value inequality and the first inequality,
we see that
\begin{eqnarray*}
\lefteqn{  
\|P(\phi_1,\phi_2,\dots,\phi_n)-P(\psi_1,\psi_2,\dots,\psi_n)\|_{C^s}  
} \\
& \le &
  2^s \sup_{0 \le t \le 1} \max_i
\| \partial_{x_i} P(t\phi_1+(1-t)\psi_1, \dots, t\phi_n+(1-t)\psi_n)  
\|_{C^s}
\sum_{i=1}^n  \|\phi_i-\psi_i\|_{C^s}   \\
& \le &
2^s (d \nu(P) 2^{s(d-1)}M^{d-1}) \sum_{i=1}^n  \|\phi_i-\psi_i\|_{C^s}  
\ ,
\end{eqnarray*}
which is  the desired result.
\Endproof

We can now use the estimate given in the above lemmas to prove the
following general  proposition. Let $r, s\geq 1$ be real numbers and
for each $w\in C^r(I)$, let $$\Theta_w: C^s(I,I)\to C^s(I)$$ be the
operator given by $\Theta_w(\varphi) = w \circ \varphi$.

\begin{proposition}\label{rs1}
Let $r, s> 1$ be real numbers both noninteger{\rm ,} and let  $w\in C^r(I)${\rm ,}
$ \varphi ,  \varpsi\in C^s(I,I)$ with
$\| \varphi-\varpsi\|_{C^s} \le 1$.
\begin{enumerate}
\item If $r>s$ then there exists
  $K=K \left( \|\varphi\|_{C^s} \right)>0$ such that
$$
\|w \circ \varphi - w \circ \varpsi\|_{C^s}\leq  K \|w\|_{C^r}
\| \varphi-\varpsi\|_{C^s}^{\epsil}
$$
where $\epsil=\min\{1-\{s\},r-s\}$ {\rm (}$\{s\}$ denotes the fractional  
part of $s${\rm ).}
In particular{\rm ,} $\Theta_w:C^s(I,I)\to C^s(I)$ is $\epsil$-\older\  
continuous.
\item  If $r=s$  there exists   $c>0$ and $K=K \left(\|\varphi\|_{C^s}  
\right)>0$  such that
$$
\|w \circ \varphi - w \circ \varpsi\|_{C^r}  \leq
c  \|w\|_{C^r}  \|\varpsi' \|_{C^0}^r
   +
K \|w\|_{C^r} \|\varphi-\varpsi\|_{C^r}^\alpha   \ ,
$$
where  $\alpha=\{s\}$ is the fractional part of $s$.
\end{enumerate}
\end{proposition}

Part (ii) of the above proposition shows one of the main difficulties
in this theory which is the fact that  for   $w\in C^s(I)$
  the operator
$\Theta_w:C^s(I,I)\to C^s(I)$ is not even $C^0$.

\Proof
Let us write $s=k+\alpha$, with $k$ an integer  and
$0 < \alpha= \{s\} < 1$, and let
$$
A=w \circ \varphi - w \circ \varpsi \ .
$$
Since  $w, \varphi, \varpsi \in C^1$ and $\epsil \le 1-\alpha$,
applying Lemma \ref{conta0} we get
$$
\|A\|_{C^\alpha} \le  K_1 \|w\|_{C^1}  
\|\varphi-\varpsi\|_{C^1}^\epsil \ .
$$
By Faa-di-Bruno's Formula (see \cite {herman}, p.42), for all $1 \le l  
\le k$
we can write
$$
A^{(l)}=B_l \left(\varphi\right) -  B_l \left(\varpsi\right)
$$
where
$$
B_l(\phi) = \sum_{j=1}^l w^{(j)}  \circ \phi \cdot
P_{l,j}(\phi',\phi'', \dots, \phi^{(l-j)}) \ ,
$$
each $P_{l,j}$ being a (universal, homogeneous) polynomial of degree $j$
in $l-j$ variables (with integer coefficients explicitly computable  
from $l$ and $j$;
see \cite[p.~42]{herman}).
We only need the expression of $P_{l,j}$ for $j=l$;
it is easy to check that $P_{l,l}(\phi')= \left(\phi'\right)^l.$
Then, we can decompose $A^{(l)}=C_l+D_l$,
where
\begin{eqnarray*}
  C_l& = &  \sum_{j=1}^l
   w^{(j)} \circ \varphi \cdot \left(
P_{l,j} \left(\varphi',\varphi'', \dots, \varphi^{(l-j)} \right)
-
P_{l,j} \left(\varpsi',\varpsi'', \dots, \varpsi^{(l-j)} \right)
\right)\ , \\
D_l & = &  \sum_{j=1}^l
\left( w^{(j)}  \circ \varphi - w^{(j)} \circ \varpsi \right) \cdot
P_{l,j} \left(\varpsi',\varpsi'', \dots, \varpsi^{(l-j)} \right)  \ .
\end{eqnarray*}
By Lemma {\ref{fr4}} applied to each $P_{l,j}$, 
$$
\left \|
P_{l,j} \left(\varphi',\varphi'', \dots, \varphi^{(l-j)} \right)
-
P_{l,j} \left(\varpsi',\varpsi'', \dots, \varpsi^{(l-j)} \right)
\right \|_{C^\alpha}
  \le  K_2      \|\varphi-\varpsi\|_{C^s}   \ .
$$         
  Therefore, for all $1 \le l \le k$ we get
$$
\|C_l\|_{C^\alpha}  \le
  K_3   \|w\|_{C^s}
   \|\varphi-\varpsi\|_{C^s}
  \ .
$$
Let us now rewrite $D_l=E_l+F_l$ where,
\begin{eqnarray*}
\label{21}
E_l &= &
\sum_{j=1}^{l-1}
\left( w^{(j)} \circ \varphi - w^{(j)} \circ \varpsi \right) \cdot
P_{l,j} \left(\varpsi',\varpsi'', \dots, \varpsi^{(l-j)} \right)\ , \\
F_l &= &
\left( w^{(l)} \circ \varphi - w^{(l)} \circ \varpsi \right) \cdot
  \left(\varpsi' \right) ^{l}  \ .
\end{eqnarray*}
In bounding the first summation in $E_l$, we apply Lemma \ref{conta0}.
Since $w^{(j)}, \varphi,\varpsi$  is at least $C^{1}$
   we get
$$
\left \|
  w^{(j)} \circ \varphi - w^{(j)} \circ \varpsi
\right \|_{C^\alpha}
  \le
K_4 \left \|w^{(j)}\right  \|_{C^{1} }  \|\varphi-\varpsi  
\|_{C^1}^{1-\alpha}
$$
for all,  with $1 \le j \le l-1$.
 From this and  Lemma {\ref{fr4}}, for all $1 \le l \le k$
we obtain that
$$
\| E_l \|_{C^\alpha} \le K_5 \|w\|_{C^s} \|\varphi-\varpsi  
\|_{C^s}^{1-\alpha} \ .
$$
Our task has been reduced to bounding the $C^\alpha$ norm of $F_l$.
Here, we will do separately  the proof of part (i) and part (ii) of
this proposition.

Let us prove part (i) first.
Here,  for all $1 \le l \le k$ we have that $\varphi, \varpsi$ are at
least $C^1$ and that  $w^{(l)}$ is at least $C^{\alpha+\epsil}$, and
so by Lemma \ref{conta0} we get
$$
\left \|w^{(l)} \circ \varphi - w^{(l)} \circ \varpsi
\right \|_{C^{\alpha+\epsil}}
  \le
K_6 \|w\|_{C^{l+\alpha+\epsil}}  \|\varphi-\varpsi\|_{C^1}^\epsil \  
.
$$
Thus,  for all $1 \le l \le k$ we obtain
$$
\| F_l \|_{C^\alpha} \le
K_7 \|w\|_{C^r}  \|\varphi-\varpsi\|_{C^s}^\epsil \ ,
$$
which ends the proof of part (i).

Let us now prove part (ii). We know that $w^{(l)}, \varphi, \varpsi \in  
C^1$
for all $1 \le l \le k-1$, and so by Lemma \ref{conta0} we have
   $$
\left \|w^{(l)} \circ \varphi - w^{(l)} \circ \varpsi
\right \|_{C^{\alpha+\epsil}}
  \le
K_8 \|w\|_{C^{l+\alpha+\epsil}}  \|\varphi-\varpsi\|_{C^1}^\alpha \ .
$$
Thus, for all $1 \le l \le k-1$ we obtain
$$
\| F_l \|_{C^\alpha} \le
K_9 \|w\|_{C^r}  \|\varphi-\varpsi\|_{C^s}^\alpha \ .
$$
Therefore, we just have to  bound $\| F_k \|_{C^\alpha}$.
Here, $w^{(k)}$ is only $C^{\alpha}$. From the   inequalities  
\eqref{aqaz7} and \eqref{aqaz9} in
the proof of   Lemma \ref{conta0},
we get
   \begin{eqnarray*}
\left \|
  w^{(k)} \circ \varphi - w^{(k)} \circ \varpsi
\right \|_{C^0}
& \le   &
c_1 \|w\|_{C^r}  \|\varphi-\varpsi\|_{C^0}^\alpha  \ , \\
\left \|
  w^{(k)} \circ \varphi - w^{(k)} \circ \varpsi
\right \|_{C^\alpha}
& \le   &
c_2 \|w\|_{C^r} \left(\left \|\varpsi'\right\|_{C^0} + \left
\|\varphi'-\varpsi' \right \|_{C^0}
\right)^\alpha  \ ,
\end{eqnarray*}
and so
\begin{eqnarray*}
\| F_k \|_{C^\alpha} & \le   &
  c_3  \|\varpsi \|_{C^{1+\alpha}}^k
     \|w\|_{C^r} \|\varphi-\varpsi\|_{C^0}^\alpha \\
&& +
  c_4 \|\varpsi \|_{C^1}^k
  \|w\|_{C^r}  \left(\left  \| \varpsi'\right\|_{C^0}  +
  \left  \|\varphi'-\varpsi'\right\|_{C^0}
\right)^\alpha  \ ,
\end{eqnarray*}
which ends the proof of part (ii).
\hfq

\begin{lemma} \label{conta1}
Given $0 \le \alpha < 1$ and $0 < \epsil \le  1- \alpha${\rm ,}
let $f \in C^{1+\alpha+\epsil}(I)${\rm ,} $g \in C^1(I,I)$ and $v
\in C^1(I)$ with
$\|v\|_{C^1}\leq 1$ and $g+v\in C^1(I,I)$.
There exists $K=K \left( \|g\|_{C^1} \right)>0$  such
that
$$
\left\|f\circ (g+v)-f\circ g - f'\circ g \cdot
v\right\|_{C^\alpha} \leq K \|f\|_{C^{1+\alpha+\epsil}}
  \|v\|_{C^1} ^{1+\epsil}
\ .
$$
In particular{\rm ,} there exists   $K=K \left(\|g\|_{C^1} \right)>0$  such  
that
$$
\left\|f\circ (g+v)-f\circ g \right\|_{C^\alpha} \leq K  
\|f\|_{C^{1+\alpha+\epsil}}
  \|v\|_{C^1}
\ .
$$
\end{lemma}

\Proof Note that we have the following identity:
$$
(f\circ (g+v)-f\circ g - f'\circ g \cdot
v)(x)=v(x)\int_{0}^{1}\left[f'(g(x)+tv(x))-f'(g(x))\right]\,dt
\ .
$$
Applying Lemma \ref{conta0} with $w=f'$, $\varphi=g+tv$ and $\psi=g$,  
we can bound
the $C^\alpha$ norm of the integrand by
$Kt\|f'\|_{C^{\alpha+\epsil}}\|v\|_{C^1}^{\epsil}$. This proves the  
first stated
estimate in slightly stronger form. The second estimate is an immediate
consequence of the first.
\hfq

\begin{proposition}\label{rs}
Let $2 \leq s+1 < r$ be real numbers{\rm ,} and let
$f\in C^r(I)${\rm ,} $g\in C^s(I,I)$. There exists
$K=K \left( \|g\|_{C^s}  \right)>0$ such that{\rm ,} for all
$v\in C^s(I)$ with
$\|v\|_{C^s} \le 1$ and $g+v\in C^s(I,I)${\rm ,}   
\begin{equation}
\left\|f\circ (g+v)-f\circ g - f'\circ g \cdot v\right\|_{C^s}\leq
K \|f\|_{C^r}\|v\|_{C^s}^{1+\theta}
\ ,
\end{equation}
where $\theta=\min\{1-\{s\},r-s-1\}$.
In particular{\rm , (a)} the operator $\Theta_f: C^s(I,I)\to C^s(I)$ is
$C^1$ and its derivative is given by $D\Theta_f(g)v=f'\circ g\cdot v${\rm ,}
and {\rm (b)} there exists $K=K \left( \|g\|_{C^s} \right)>0$   such that
for all $v$ as above{\rm ,}
\begin{equation}
  \|f\circ (g+v)-f\circ g \|_{C^s}\leq
K \|f\|_{C^r}\|v\|_{C^s}
\ .
\end{equation}
\end{proposition}

\Proof  In this proof we use $K_1,K_2,\dots$ to denote constants
that depend only on $  \|g\|_{C^s}$. Consider the remainder term
$$
F= f\circ (g+v) - f\circ g -f'\circ g\cdot v
\ ,
$$
as well as its derivative $F'=A+B$, where
\begin{eqnarray*}
A & =& \left(f'\circ (g+v) - f'\circ g -f''\circ g\cdot
v\right)\cdot g' \ , \\
B & = &
\left(f'\circ (g+v) - f'\circ g \right)\cdot v'
\ .
\end{eqnarray*}
We want to show that
$$
\left\|F'\right\|_{C^s} \leq
K_1 \|f\|_{C^r}\|v\|_{C^s}^{1+\theta}\ .
$$
The proof will be by induction on the integral part of $s$.
Note however that the mean value theorem already gives us
$\|F\|_{C^0}\leq K_2\|f''\|_{C^0}\|v\|_{C^0}^2$ independently
of $s$.

First we deal with the base of induction, namely  when $1\leq s<2$, say,
$s=1+\alpha$.
By Lemma {\ref{conta1}}, we have
$$
\left\|A \right\|_{C^\alpha}
\leq
K_3 \|f'\|_{C^{1+\alpha+\theta}}
\left(\|v\|_{C^1}\right)^{1+\theta}
\ .
$$
The same Lemma {\ref{conta1}} yields
$$
\left\|B \right\|_{C^\alpha}
\leq
K_4 \|f'\|_{C^{1+\alpha+\theta}}
\left(\|v\|_{C^{1+\alpha}}\right)^2
\ .
$$
This establishes the base of induction.

Now suppose that our lemma holds for $s>1$. We will prove from this
that it holds for $s+1$. To do this, it suffices to show that
\begin{equation}\label{R1}
\left\|F'\right\|_{C^s} \leq
K_5 \|f\|_{C^r}\|v\|_{C^{s+1}}^{1+\theta}
\ .
\end{equation}
The proof is more of the same. By the induction hypothesis applied to
$f'$, we have
\begin{equation}\label{R2}
\left\|A \right\|_{C^s}
\leq
K_6 \|f'\|_{C^{r-1}}
\left(\|v\|_{C^s}\right)^{1+\theta}
\ .
\end{equation}
The same fact also gives
\begin{equation}\label{R3}
\left\|B\right\|_{C^s}
\leq
K_7 \|f'\|_{C^{r-1}}
\left(\|v\|_{C^{1+s}}\right)^2
\ .
\end{equation}
Putting \eqref{R2} and \eqref{R3} together we get \eqref{R1}, and so
the induction is complete.
\Endproof

\begin{proposition}\label{comp}
Let $2 \leq s+1 < r$ be real numbers. The
composition operator $\Theta: C^r(I)\times C^s(I,I)\to C^s(I)$ given
by $\Theta(f,g)=f\circ g$ is $C^{1+\theta}$ and its derivative is given  
by
$D\Theta(f,g)(u,v)=u\circ g +f'\circ g\cdot v$.
In particular{\rm ,} there exists
$K=K \left( \|f\|_{C^r}, \|g\|_{C^s}  \right)>0$
such that{\rm ,} for all $\|u\|_{C^r}\leq 1$
and $\|v\|_{C^s}\leq 1$ with $g+v\in C^s(I,I)${\rm ,}  
\begin{equation}
\label{1agas}
\|\Theta(f+u,g+v)-\Theta(f,g)-D\Theta(f,g)(u,v) \|_{C^s}
\le K  (\|u\|_{C^r}+\|v\|_{C^s})^{1+\theta} \ ,
\end{equation}
where $\theta=\min\{1-\{s\},r-s-1\}$.
\end{proposition}

\Proof In this proof, we denote by $K_1,K_2,\dots$ positive
constants depending only on $\|g\|_{C^s}$.
Let us take $u\in C^r(I)$ and $v\in
C^s(I)$ such that $\|u\|_{C^r}\leq 1$
and $\|v\|_{C^s}\leq 1$, respectively.
Now,
\begin{eqnarray*}
F & = & \Theta (f+u,g+v)-\Theta(f,g) -u\circ g-f'\circ g\cdot
v \\
& = &
f\circ (g+v)-f\circ g - f'\circ g\cdot v  +
u\circ (g+v)-u\circ g
\ .
\end{eqnarray*} \nopagebreak
Using Proposition {\ref{rs}}, we see that
$$
\left\| f\circ (g+v)-f\circ g - f'\circ g\cdot v  \right\|_{C^s}
\leq K_1 \|f\|_{C^r} (\|v\|_{C^s})^{1+\theta}
\ . \pagebreak
$$
The same Proposition {\ref{rs}} with $u$ replacing $f$ yields
\begin{eqnarray*}
\left\|u\circ (g+v)-u\circ g\right\|_{C^s}
& \leq &
\|u'\circ g\cdot v \|_{C^s} + K_2 \|u\|_{C^r} (\|v\|_{C^s})^{1+\theta}  
\\
& \leq &  K_3 \|u\|_{C^r} \|v\|_{C^s}
\ .
\end{eqnarray*}
Therefore we get
$$
\|F\|_{C^s} \leq
K_1 \|f\|_{C^r} (\|v\|_{C^s})^{1+\theta} +
K_3 \|u\|_{C^r} \|v\|_{C^s}
\ ,
$$
which proves that $\Theta$ is $C^1$ and that \eqref{1agas} is satisfied.
Now, we have that
$$
D \Theta (f+\phi,g+\psi) (u,v) -
D \Theta (f,g) (u,v)
= A + B + C
$$
where
\begin{eqnarray*}
A & = & u \circ (g+\psi) - u \circ g\ ,  \\
B  & = &  \left(f' \circ (g+\psi) - f' \circ g \right) \cdot v\ , \\
C  & = &  \phi' \circ (g+\psi) \cdot v \ .
\end{eqnarray*}
By Proposition \ref{rs1}, we
  obtain that
\begin{eqnarray*}
\|A\|_{C^s}& \leq & K_4 \|u\|_{C^{r-1}}  \|\psi\|_{C^s}^\theta \ ,\\
\|B\|_{C^s}& \leq & K_5 \|f\|_{C^r}  \|\psi\|_{C^s}^\theta \cdot  
\|v\|_{C^s}.
\end{eqnarray*}
Letting $k$ be the integral part of $s$ and $\varphi=g+\psi$,
and using Faa-di-Bruno's Formula,  we have
$$
\left(\phi' \circ \varphi \right)^{(k)}=
\sum_{j=1}^k  \phi^{(j+1)}  \circ \varphi \cdot
P_{k,j}(\varphi',\varphi'', \dots, \varphi^{(k-j)})
\ ,
$$
each $P_{k,j}$ being a (universal, homogeneous) polynomial of degree $j$
in $k-j$ variables. Hence, using Lemma \ref{fr4}, we get that
$\|C\|_{C^s}= K_6  \|\Delta f\|_{C^r}  \|v\|_{C^s}$.
Thus, $\Theta$ is a $C^{1+\theta}$ operator.
\hfq

\begin{corollary}\label{um}
Let $r,s > 0$ be real  numbers with  $r-1>s\ge 1$
and  for each positive integer $m${\rm ,} let
$Q_m:C^r(I,I) \to C^s(I,I)$ be  the operator  given by $Q_m(f)= f^m$.
\begin{enumerate}
\item
Let $0\leq t\leq r$  and let $U:C^t(I,I)\to C^s(I,I)$ be a  
$C^{1+\theta}$ operator
for some $0<\theta<1$.
Then  the operator $U_m:C^r(I,I) \to C^s(I,I)$  given by $U_m(f)= Q_m  
\circ U(f)$
is $C^{1+\theta'}$ for some  $0<\theta'=\theta'(\theta,r,s)<1$.
\item
In particular{\rm ,}  the operator $Q_m:C^r(I,I) \to C^s(I,I)$  is  
$C^{1+\theta''}$
  for some  $0<\theta''=\theta''(r,s)<1$
and
there exists $K=K \left(m, \|f\|_{C^r} \right)>0$
such that\end{enumerate}
\vglue-22pt
\begin{equation}
\label{3agas}
\|Q_m(f+u)-Q_m(f)-DQ_m(f)u \|_{C^s} \le K \|u\|_{C^r}^{1+\theta''} \ .
\end{equation}
\end{corollary}
 
\vskip8pt

\Proof
First note that $U_{m+1}(f)=\Theta(f, U_m(f))$. The operator
$U_1$ arises as the composition of the operator $C^r(I,I)\to
C^r(I,I)\times C^s(I,I)$ given by $f\mapsto (f,U(f))$,
which is $C^{1+\theta}$ because $U$ is $C^{1+\theta}$ (and $C^r(I,I)$  
embeds in $C^t$),
with the composition operator $\Theta: C^r(I,I)\times C^s(I,I)\to
C^s(I,I)$, which is $C^{1+\theta''}$
for some  $0<\theta''=\theta''(r,s)<1$
by Proposition {\ref{comp}}.
The desired result for part (i) then follows by induction.
Part (ii) is a  corollary of part (i),
  and by a computation  \eqref{3agas} follows from \eqref{1agas}.
\hfq

\begin{proposition}\label{Gees}
Let $r,s,t$ be real numbers with $2 \le s+1< r$  and $t\ge
0$. Let $U:C^t(I,I)\to C^s(I,I)$ be a $C^1$ operator. Then for each
$\phi\in C^r(I)$ and each $\psi\in C^t(I,I)$ there exists a function
$\sigma_\psi:\mathbb{R}^+ \to \mathbb{R}^+$ with $\sigma_\psi(h)/h\to  
0$ as
$h\to 0${\rm ,} varying continuously with $\psi${\rm ,} such that for all $v\in
C^t(I)$ with $\psi + v\in C^t(I,I)${\rm ,}
\begin{equation}
\left\|\phi\circ U(\psi + v)-\phi\circ U(\psi)- \phi'\circ
U(\psi)\cdot DU(\psi)v\right\|_{C^s}\leq
  \sigma_\psi (\|v\|_{C^t})
\ .
\end{equation}
\end{proposition}

\Proof As before, we denote by $K_1,K_2,\dots$ positive
constants that depend only on $\|\psi\|_{C^t}$.
We have that
$$
\phi\circ U(\psi + v)-\phi\circ U(\psi)- \phi'\circ
U(\psi)\cdot DU(\psi)v =  A+B
$$
where
\begin{eqnarray*}
  A & = &
  \phi\circ U(\psi + v)-\phi\circ U(\psi)- \phi'\circ
U(\psi)\cdot (U(\psi + v) - U(\psi)) \ , \\
  B & = & \phi'\circ U(\psi)\cdot\left( U(\psi + v) - U(\psi)  
-DU(\psi)v\right) \   .
\end{eqnarray*}
Since $U$ is $C^1$, there exists a continuous function
$\nu_\psi:\mathbb{R}^+ \to \mathbb{R}^+$ with $\nu_\psi(h)/h\break\to 0$ as
$h\to 0$, varying continuously with $\psi$, such that
$$
\left\|U(\psi + v) - U(\psi) -DU(\psi)v\right\|_{C^s}\leq
  \nu_\psi(\|v\|_{C^t})
\ .
$$
Hence, applying Proposition {\ref{rs}} with $f=\phi$ and $g=U(\psi)$
and $v$ replaced by $U(\psi + v) - U(\psi)$, we get
\begin{eqnarray*}
\|A\|_{C^s} & \le &
    K_2\|\phi\|_{C^r}
\left(\|U(\psi + v) -
U(\psi)\|_{C^s}\right)^{1+\theta}  \\
& \le &
K_3 \|\phi\|_{C^r}\left(\|DU(\psi)\|^{1+\theta}\|v\|_{C^t}^{1+\theta}  
\right) \ ,
\end{eqnarray*}
and
\begin{eqnarray*}
\|B\|_{C^s} & \le &
K_4\|\phi\|_{C^r}\nu_\psi(\|v\|_{C^t})  \ ,
  \end{eqnarray*}
where $K_3=K_3\left(\|U(\psi)\|_{C^s},\|DU(\psi)\|, \nu_\psi\right)$
and $K_4=K_2\left(\|U(\psi)\|_{C^s}\right)$.
Therefore,
$$\|A + B\|_{C^s}
\le
K_3\|\phi\|_{C^r}\|v\|_{C^t}^{1+\theta} +   
K_4\|\phi\|_{C^r}\nu_\psi(\|v\|_{C^t})   \ .
$$
This completes the proof.
\hfq

\begin{corollary}\label{vm}
Let $r,s,t$ be real numbers with  $r-1>s> 1$
and $0\leq t\leq r${\rm ,} and let $U:C^t(I,I)\to C^s(I,I)$ be a $C^1$
operator.
For each positive integer $n$ the operator
$V_n:C^r(I,I)\to C^s(I)$ given by $V_n(f)=(f^m)'\circ U(f)$
is differentiable at every $g\in C^{r+1}(I,I)\subseteq C^r(I,I)${\rm ,} and
as a map from $C^{r+1}(I,I)$ into $\mathcal{L}(C^r(I), C^s(I))${\rm ,}
the derivative operator $g\mapsto DV_n(g)$ is continuous.
\end{corollary}

\Proof
First note that by the chain rule,
$$
V_{n}(f)=\prod_{j=0}^{n-1}f'\circ \left(f^{(j)}\circ U(f)\right)
=\prod_{j=0}^{n-1}f'\circ U_j(f)
\ .
$$
This reduces the problem to the case $n=1$. We claim that the
linear operator
$$
L(v)=v'\circ U(g)+ g''\circ U(g)\cdot DU(g)v
$$
is the derivative of $V_1$ at $g\in C^{r+1}(I,I)$. Indeed, we have
$$
V_1(g+v)-V_1(g)-L(v) = A+B \ ,
$$
where
\begin{eqnarray*}
A &=& g'\circ U(g+v)-g'\circ U(g) -g''\circ U(g)\cdot
DU(g)v \ , \\
B &=&    v'  \circ U(g+v) -v'\circ U(g)
\ .
\end{eqnarray*}
By Proposition {\ref{Gees}} applied to $\phi=g'$ and $\psi=g$, there
exists $K_1=K_1(\|g\|_{C^{r+1}})$ such that
$$
\left\|A \right\|_{C^s}\leq
K_1\sigma_\psi (\|v\|_{C^r})
\ ,
$$
where $\sigma_\psi:\mathbb{R}^+\to\mathbb{R}^+$ is a continuous
function varying continuously with $\psi$ such that
$\sigma_\psi(h)/h\to 0$ as $h\to 0$.
On the other hand, by part (i) of Proposition {\ref{rs1}} and
since $U$ is $C^1$, we have
\begin{eqnarray*}
\left\|B \right\|_{C^s}&\leq&
K_2\|v\|_{C^r}\|U(g+v)-U(g)\|_{C^s}^{\epsil}\\
&\leq& K_3\left(\|v\|_{C^r}\right)^{1+\epsil}
\ ,
\end{eqnarray*}
where $0<\epsil=\min\{1-\{s\},r-s\}<1 $,
$K_2=K_2(\|U(g)\|_{C^s})$ and $K_3=K_3(\|U(g)\|_{C^s},$  
$\|DU(g)\|,\sigma_{\psi} (\|v\|_{C^r}))$.
Combining these inequalities, we deduce that $V_1$ is differentiable
at $g$ and $DV_1(g)=L$ as claimed. It is clear from the expression
defining it that $L$ varies continuously with $g\in C^{r+1}(I,I)$.
\hfq


\Subsec{Checking properties $\mathbf{B2}$ and $\mathbf{B3}$}
\label{sec:axiom5}
We now proceed to verify that the operator $T$ satisfies
  properties $\mathbf{B2}$  and $\mathbf{B3}$  of robustness.
They will follow respectively from Lemmas \ref{qn} and \ref{333}.
First it is necessary to analyze the
behavior of the linear scaling used in such operators.
Let us fix a positive integer $p$ and for each $f\in C^r(I,I)$ let
$\Lambda_f$ be the linear map $x\mapsto \lambda_f x$, where
$\lambda_f=f^p(0)$.

\begin{lemma} \label{lL}
For $r>2${\rm ,} the maps $\Lambda: C^r(I,I)\to  
\mathcal{L}(\mathbb{R},\mathbb{R})$
given by $\Lambda(f)=\Lambda_f$ and $\lambda: C^r(I,I)\to \mathbb{R}$
given by $\lambda(f)=\lambda_f$ are both $C^{1+\theta}$
for some $0<\theta=\theta(r,s)<1$. In particular{\rm ,} there is $K=K(p,  
\|f\|_{C^r})>0$
such that for all $v \in C^r(I)$ with $\|v\|_{C^r} \le 1$ and $f+v\in
C^r(I,I)${\rm ,}  
\begin{equation}\label{5agas}
\|\lambda(f+v)-\lambda(f)-D\lambda(f)v\|_{C^r} \le  K  
\|v\|_{C^r}^{1+\theta} \ .
\end{equation}
  The above  inequality also holds when $\lambda$ is replaced by  $\Lambda$.
\end{lemma}

\Proof Choosing $1< s <r-1$, we see that $\lambda= E\circ
Q_p$ where $Q_p: C^r(I,I)\to C^s(I,I)$ is the operator $Q_p(f)=f^p$,
which is $C^{1+\theta}$ for some $0<\theta=\theta(r,s)<1$,
  and $E: C^s(I,I)\to \mathbb{R}$ is the evaluation map
$E(g)=g(0)$, which is linear. Therefore, by Corollary \ref{um},   
$\lambda$ is
$C^{1+\theta}$ and
\eqref{5agas} follows from \eqref{3agas} and the linearity of $E$.
  The proof
for $\Lambda$ is entirely analogous.
\Endproof

We will also need to use the operators $U_n:C^r(I,I)\to C^s(I)$
given by $U_n(f)=f^n \circ \Lambda_f$ for all $n \ge 0$.

Property $\mathbf{B2}$ for the operator $T$ is a consequence of the
following lemma (the first assertion in $\mathbf{B2}$ is actually a
consequence of Lemma \ref{fr2} below).

\begin{lemma}\label{qn}
For $2< s +1 < r$  and for each $n\geq 0${\rm ,}
the operator  $U_n:C^r(I,I)\to C^s(I)$
is $C^{1+\theta}$
for some $0<\theta=\theta(r,s)<1$.
In particular{\rm ,}    $T^m:\mathbb{O}^r_m \to \mathbb{U}^s$  is also a
$C^{1+\theta}$ operator.
\end{lemma}

\Proof
This follows at once from Lemma {\ref{lL}} and Corollary {\ref{um}}
applied to $U=\Lambda$.
\Endproof

The following lemma is all we need to verify property $\mathbf{B3}$
for the operator~$T$. In this case $g$ is a map in the limit set
$\mathbb{K}$ of $T$, hence analytic, and $v=\mathbf{u}_g$ is a
tangent vector to the unstable manifold of $g$, which is analytic as
well.

\begin{lemma}
\label{333}
For $2< s +1< r${\rm ,} the map $\mathbb{O}^r_m \to \mathbb{U}^s$
given by $f\mapsto DT^m(f)v$ is differentiable at $f=g \in \mathbb{K}$.
Furthermore{\rm ,} for every $m \ge 1$ there exist  $C_m>1$ and $\nu_m > 0$
such that for each $g \in \mathbb{K}$ and $f \in \mathbb{O}^r_m$
with $\|f-g\|_{C^r} < \nu_m$   and all $v \in \mathbb{A}^r$
with $\|v\|_{C^r}=1${\rm ,}
\begin{equation}
\label{gtwewgg}
\|DT^m(f)v - DT^m(g)v  \|_{C^s} \leq
  C_m  \|f-g\|_{C^r}
\ .
\end{equation}
\end{lemma}

\vglue8pt
\Proof
Let $E: C^s(I,I)\to \mathbb{R}$ be the evaluation map
$E(g)=g(0)$, which is linear.
Recall that the derivative of $T^m$ is given by the expression
\begin{eqnarray}\label{deeT}
&&\qquad DT^m(f)v=
\frac {1}{\lambda_f}
\sum_{j=0}^{p-1}(f^j)'\circ U_{p-j}(f) \cdot v\circ U_{p-j-1}(f)  \\
&&\qquad
+\frac {1}{\lambda_f}
[\text{id} \cdot (T^mf)'-T^mf]
\sum_{j=0}^{p-1} E \left( (f^j)' \circ U_{p-j}(f) \right ) \cdot E  
\left( v \circ U_{p-j-1}(f) \right)
\ , \nonumber
\end{eqnarray}
where $\lambda_f=E \circ f^p$ and $\text{id}:\mathbb{R} \to \mathbb{R}$  
is the identity map.
  Each term of the first summation in \eqref{deeT} is
differentiable at $f=g$. To see this apply Lemma \ref{qn} and Corollary
{\ref{vm}} to each of the operators $f\mapsto (f^j)'\circ U_{p-j}(f)$
as well as Proposition {\ref{rs}} to each of the operators $f\mapsto
v\circ U_{p-j}(f)$. On the other hand, each term of the second
summation in \eqref{deeT} equals the corresponding term in the first
summation post-composed with the evaluation map $E$  (which is
linear), and is therefore differentiable at $f=g$. The analysis of the
expression in square brackets in \eqref{deeT} is similar.
By Lemma \ref{lL} and Corollary {\ref{vm}}, the
operator $f\mapsto T^m(f)'=(f^p)'\circ\Lambda_f$ is   differentiable
at $f=g$, and    the
operator $f\mapsto T^m(f)=\lambda_f \cdot  f^p \circ\Lambda_f$  is
also differentiable at $f=g$ by  Lemma \ref{lL} and Corollary \ref{um}.
 From this fact and compactness of $\mathbb K$ the inequality
\eqref{gtwewgg} follows.
\hfq


\Subsec{Checking property $\mathbf{B4}$}
\label{sec:gateaux}
The fifth property is verified in Lemma \ref{4444} below.
First we will need to prove the following two lemmas about  the  
operators
$U_i:C^{t+1+\varepsilon}(I,I) \to C^t(I)$  with $t \ge 1$.
Recall that $U_i(f)=f^i \circ \Lambda_f$.

\begin{lemma}\label{fr2}
For every $f\in C^{t+1+\epsil}(I,I)$
and all $v \in C^t(I)$
with small norm and such that $f+v \in C^{t}(I,I)${\rm ,}  
$$
\|U_i(f+v)-U_i(f)\|_{C^t} \le K \|v\|_{C^t}
$$
for all $0 \le i \le p$   where
$K=K\left(p,\|f\|_{C^{t+1+\epsil}}\right)$.
\end{lemma}

\Proof
Note that
$$
U_{i+1}(f+v)-U_{i+1}(f) =
f \circ U_i(f+v)- f \circ U_i(f) + v \circ U_i(f+v)\ .
$$
By Proposition \ref{rs}, there is  
$K_1=K_1\left(p,\|f\|_{C^{t+1+\epsil}}\right)$ such that
$$
\|U_{i+1}(f+v)-U_{i+1}(f)\|_{C^t} \le
K_1 \|U_i(f+v)- U_i(f)\|_{C^t} + \|v \circ U_i(f+v)\|_{C^t}.
$$
The required estimate now follows by induction,
because $U_0$ is $C^1$.
\hfq

\begin{lemma}\label{fr3}
For every $f\in C^{t+1+\epsil}(I,I)$
and all $v \in C^{t+1+\epsil}(I)$
with small norm and such that $f+v \in C^{t+1+\epsil}(I,I)${\rm ,} 
  $$
\|U_i(f+v)-U_i(f)-DU_i(f)v\|_{C^t} \le
K \|v\|_{C^t}^{1+\theta} \ ,
$$
   for all $0 \le i \le p${\rm ,} for some  $0 < \theta=\theta(t,\epsil) <1$  
  and
$K=K(p,\|f\|_{C^t})>0$.
\end{lemma}

\Proof
In this proof we denote by $K_1, K_2,\dots$ the positive constants  
depending only on $m$ and
$\|U_i(f)\|_{C^{t+1+\epsil}}$.
Again we use induction;
the case $i=0$ follows from  the differentiability of the scaling $f\to  
\Lambda_f$
and inequality \eqref{5agas}.
We have
$$
U_{i+1}(f+v)-U_{i+1}(f)-DU_{i+1}(f)v=A+B+C
$$
where
\begin{eqnarray*}
A &=& f \circ U_i(f+v)-f \circ U_i(f)-f' \circ U_i(f) \cdot  
\left(U_i(f+v)-U_i(f) \right) \ ,\\
B &=& f' \circ U_i(f) \cdot \left(U_i(f+v)-U_i(f)-DU_i(f)v \right)\ , \\
C &=&  v \circ U_i(f+v)- v \circ U_i(f) \ .
\end{eqnarray*}
By Proposition \ref{rs},
$$\|A\|_{C^t} \le K_1 \|U_i(f+v)-U_i(f)\|^{1+\epsil}_{C^t}.$$
Using Lemma {\ref{fr2}}, we get
$$
\|A\|_{C^{t-1}} \le K_2 \|v\|_{C^t}^{1+\epsil} \ .
$$
On the other hand, since $v$ is $C^{t+1+\epsil/2}$,
we know again from Proposition \ref{rs} that
$$
\|C\|_{C^t} \le K_3 \|v\|_{C^{t+1+\epsil/2}}
\|U_i(f+v)-U_i(f)\|_{C^t}
\le K_4 \|v\|_{C^{r+1+\epsil/2}}  \|v\|_{C^t}
  \ .
$$
Since $v$ has bounded $C^{t+1+\epsil}$ norm,
by an interpolation of norms,
we have $ \|v\|_{C^{t+1+\epsil/2}} \le K_5  \|v\|_{C^t}^{\theta_1}$
for some $\theta_1 > 0$. Therefore, taking $\theta=\min \{\epsil,  
\theta_1\}$
we get
$$
\|C\|_{C^t} \le K_6  \|v\|_{C^t}^{1+\theta} \ .
$$
This allows the induction as desired.
\Endproof
\vskip4pt

Property $\mathbf{B4}$ for the operator $T$ is a direct consequence
of  the following lemma.

\begin{lemma}
\label{4444}
For every $f\in \mathbb{O}^{s+1+\epsil}$
and all $v \in \mathbb{A}^{s+1+\epsil}$
with small norm such that $f+v \in \mathbb{O}^{s+1+\epsil}${\rm ,} 
$$
\|T(f+v)-T(f)-DT(f)(v)\|_{C^s }\leq
K  \|v\|^{1+\tau}_{C^s}
\ ,
$$
for some $0<\tau=\tau(s,\epsil) <1$ and  
$K(p,\|f\|_{C^{s+1+\epsil}})>0$.
\end{lemma}

\Proof
In this proof we denote by $K_1, K_2,\dots$ the positive constants  
depending only on $m$ and
$\|U_i(f)\|_{C^s}$.
Start, observing that since $T(f)=\lambda_f^{-1}\break \cdot U_p(f)$,
$$
T(f+v)-T(f)-DT(f)v=A+B+C
$$
where
\begin{eqnarray*}
A & = & \lambda_f^{-1} \cdot \left( U_p(f+v)-U_p(f)-DU_p(f)v\right)\ , \\
B &=& \left(\lambda_{f+v}^{-1}-\lambda_f^{-1}-D\lambda_f^{-1}  
(v)\right) \cdot U_p(f+v)\ , \\
C &=& D \lambda_f^{-1}(v)  \cdot \left( U_p(f+v)-U_p(f)\right) \ .
\end{eqnarray*}
Applying Lemma \ref{fr3} with $t=s$ we get
$$
\|A\|_{C^s}= K_1 \|v\|_{C^s}^{1+\theta_1} \ ,
$$
for some $0<\theta_1=\theta_1(s,\epsil) < 1$.
By Lemma \ref{lL} there is $0<\theta_2=\{s\}<1$ such that
$$
\|B\|_{C^s} \le  K_2 \|v\|_{C^s}^{1+\theta_2}
\ .
$$
By Lemma \ref{fr2}, we have $\|U_p(f+v)-U_p(f)\|_{C^s} \le K_3  
\|v\|_{C^s}$
and so
$$
\|C\|_{C^s} \le K_4 \|v\|_{C^s}^2
\ .
$$
Therefore, it is enough to take $\tau= \min \{\theta_1, \theta_2\}$.
\hfq

\Subsec{Checking properties $\mathbf{B5}$ and $\mathbf{B6}$}
\label{sec:frechet}
We now move on to the task of proving that the operator $T=R^N$ of
Theorem \ref{lyubhyp} satisfies properties $\mathbf{B5}$ and
  $\mathbf{B6}$ in the definition of robustness.
Unlike the previous ones, the verification of these
(last) two properties depends upon the geometry of
the post-critical sets of maps near in $\mathbb{A}$ to the limit set
$\mathbb{K}$ of $T$.
The estimates performed here are the most delicate,
and involve the results of \S \ref{sec:bdedgeom}.

Recall that $T^m$ is well-defined on an open set
$\mathbb{O}_m$ in the Banach space   $\mathbb{A}=\mathbb{A}_{\Omega_a}$
(see \S \ref{grtsfd}), which contains $\mathbb{K}$.
We shall denote the renormalization intervals
$\Delta_{0,mN},$ $\Delta_{1,mN},$ $\dots$, $\Delta_{p,mN}$ simply
by $\Delta_i=\Delta_{i,mN}$
(this shortened notation should cause no harm,
because $N$ is fixed since Theorem \ref{lyubhyp}, and $m$ will be fixed
in the particular estimates involving these intervals).

We can write the derivative of $T^m$ in the following form
$$ DT^m(f)v= A(f)  \sum_{j=0}^{p-1} B_j(f) \cdot C_j(f)
+A(f) \cdot D(f)  \sum_{j=0}^{p-1} E \circ  B_j(f) \cdot E \circ C_j(f)
$$
where $E$ is the evaluation map and
\begin{alignat*}{5}
A (f) & =   \left(\lambda_f\right)^{-1} \ ,&\qquad 
B_j (f)   &=    \left(f^j\right)'\circ U_{p-j}(f) \ , \\
C_j (f) & =  v\circ U_{p-j-1}(f) \ , & \qquad
D (f)  &=     \text{id} \cdot (f^p)'\circ U_0(f) -\lambda_f \cdot U_p
(f) \ .   
\end{alignat*}

To carry out our estimates for $T^m$, we shall use the operators
$U_i:f \mapsto f^i \circ \Lambda_f$ ($i\ge 0$).
Note that $U_0(f)=\Lambda_f$, hence $U_0$ is $C^1$ in whichever
space $C^r(I,I)$ we work in, because the scaling $f \mapsto \Lambda_f$
is $C^1$ by Lemma \ref{lL}.

First we need some estimates for $U_i$.
It is clear that $\|U_i(f)\|_{C^0} \le 1$ always, but more is true.

\begin{lemma}\label{fr1}
There exists $C > 0$ with the following property. For
every $m>0${\rm ,} there exists an open neigbourhood $\mathcal{O}_m \subset
\mathbb{O}_m$ of $\mathbb{K}$ such that for all  $f \in
\mathcal{O}_m${\rm ,}  
$$
\|B_j(f)\|_{C^0} \le
C \frac{|\Delta_{0}| }{ |\Delta_{p-j}|} \ ,
$$
for all $0  \le j \le p-1$.
Furthermore{\rm ,}
$\|U_i(f)' \|_{C^0} \le
C  |\Delta_{i}|${\rm ,} for all $0 \le i \le p$.
\end{lemma}

\Proof
Use bounded distortion and the real bounds (see \S \ref{sec:bdedgeom}).
\hfq

\begin{lemma}\label{fr23}
For all $f \in \mathcal{O}_m$ and all $v \in \mathbb{A}^r$
with small norm{\rm ,} 
$$
\|U_i(f+v)-U_i(f)\|_{C^r} \le K \|v\|_{C^r}
$$
for all $0 \le i \le p${\rm ,}  where $K=K(m)>0$.
\end{lemma}

\Proof This lemma follows from Lemma \ref{fr2}.
\Endproof

Next, we show an essential result to prove that the renormalization
operator satisfies properties $\mathbf{B5}$ and $\mathbf{B6}$.
  Here, we use again in a crucial way the geometric properties
of the postcritical set of $f \in \mathbb{O}_m$  proved in Section
\ref{hfgrtyresd}.

\begin{proposition}\label{fr33333}
{\rm (i)}
For every  $t > 2$ which is not an integer there exist $0 < \mu < 1$  
and $C>0$ with the
following property. For every  $g \in \mathbb{K}$ and for every $m>0${\rm ,}  
there is an $\eta>0$
such that for all $f \in \mathbb{O}_m$ with
$\|f-g\|_{\mathbb{A}} < \eta$  and for all $w \in \mathbb{A}^t$  with
$\|w\|_{C^t} < \eta${\rm ,}
\begin{equation}
\label{ooijkjkl}
    \left \|  A(f) \sum_{j=0}^{p-1} B_j(f)
\left (C_j(f+w)-C_j(f)   \right) \right \|_{C^t}
   \le  C  \mu^m \|v\|_{C^t} \ .
\end{equation}

{\rm (ii)} For every $  \mu > 1$ close to one{\rm ,} there is $s<2$ close to two
and $C>0$ with the following property\/{\rm :}\/
for every  $g \in \mathbb{K}$ and for every $m${\rm ,}  there is an $\eta>0$
such that for all $f \in \mathbb{O}_m$    with
$\|f-g\|_{\mathbb{A}} < \eta$ and for all $w \in \mathbb{A}^t$      with
$\|w\|_{C^t} < \eta${\rm ,} the inequality \eqref{ooijkjkl} above is also satisfied.
\end{proposition}

\Proof
Below, the positive constants $c_1,c_2,\dots$ depend only on $t$ (and
the real bounds), while the positive constants  $K_0,K_1, K_2,\dots$
may depend also on~$m$.

Let  $k$ and $0<\alpha <1$ be respectively the integer and the
fractional part of $t$ (when $t=k+\text{Lip}$, take $\alpha=1$).
We start observing  that for each $j$,
\begin{eqnarray}
\label{eedsfggr}
&&\\
\left \|      B_j(f)
\left (C_j(f+w) - C_j(f)   \right)  
\right \|_{C^t}
&\le&
   \|B_j(f)\|_{C^0} \|C_j (f+w)-C_j (f) \|_{C^t}   \nonumber\\
&& +   K_0
     \|B_j(f)\|_{C^t} \|C_j (f+w)-C_j (f) \|_{C^{k}}    \ .\nonumber
\end{eqnarray}
Note that on the right-hand side of \eqref{eedsfggr} only the second
term carries a constant $K_0$. By Lemma \ref{fr1}, there is $c_1 > 0$  
such that
for every integer $m$ there is an open neighbourhood $\mathcal{O}_m$ of
$\mathbb{K}$ with the property that for each $f\in \mathcal{O}_m$ we  
have
\begin{equation} \label{eeffeef}
   \| B_j(f)\|_{C^0}  \le  c_1 \frac { | \Delta_0 |} {|\Delta_{p-j}|} \ .
\end{equation}
In that neighborhood, we also have $\| B_j(f)\|_{C^t}\leq K_1$.
By Proposition \ref{rs1} and   Lemma \ref {fr23}, taking $0<\epsil <  
1$
such that $\alpha-\epsil > 0$, we obtain
\begin{eqnarray}\label{ikefp}\qquad \quad
\|C_j(f+w)-C_j(f) \|_{C^{k}} & \le &
  \|C_j(f+w)-C_j(f)  \|_{C^{t-\epsil}}   \\
& \le &
  K_2 \|v\|_{C^t} \|U_{p-j-1}(f+w)-U_{p-j-1}(f)\|_{C^t}^\epsil  
\nonumber  \\
& \le &  K_3 \|v\|_{C^t} \|w\|_{C^t}^\epsil   \ .\nonumber
\end{eqnarray}
On the other hand, putting together Proposition \ref{rs1} with
Lemma \ref{fr1} and with Lemma \ref {fr23}, we get
\begin{eqnarray} \label{okweoo} \qquad\quad
  \|C_j(f+w)-C_j(f)  \|_{C^t}
& \le & c_2  \|U_{p-j-1}(f)'\|_{C^0}^t  \|v\|_{C^t} \\
 & &   +\, K_4  \|  U_{p-j-1}(f+w)-U_{p-j-1}(f)\|_{C^t}^\alpha
\|v\|_{C^t}  \nonumber \\
 & \le & c_3 |\Delta_{p-j-1}|^t \|v\|_{C^t} + K_5  \|w\|_{C^t}^\alpha
\|v\|_{C^t} \ . \nonumber
\end{eqnarray}
The first term on the last line of \eqref{okweoo} looks a bit
dangerous. What saves us here is the geometric control on the
post-critical set of $f$ (hence on the intervals $\Delta_i$) that we
have had at our disposal since \S \ref{sec:bdedgeom}.
Substituting      \eqref{eeffeef}, \eqref{ikefp} and \eqref {okweoo}
in   \eqref{eedsfggr} and adding up the terms with $j=0,\dots,p-1$
we arrive~at
\begin{multline*} \label{h15}
  \left \|  A(f) \sum_{j=0}^{p-1} B_j(f)
\left (C_j(f+w)-C_j(f)   \right) \right \|_{C^t}   \\
  \le
c_4 \left\{\frac{1}{|\Delta_0|}   \sum_{j=0}^{p-1}
\frac { | \Delta_0 | \cdot |\Delta_{p-j-1}|^t}
{|\Delta_{p-j}|}
+ K_5 \|w\|_{C^t}^\epsil \right\}\|v\|_{C^t} \ . 
\end{multline*}
But as we have seen in \S \ref{sec:bdedgeom}:
\begin{enumerate}
\item By Proposition \ref{zaf} and Remark \ref{515151}, if $t > 2$  
there exist $0 <
\gamma < 1$ and $C>0$ with the following property. For every  $g \in
\mathbb{K}$ and every $m>0$, there exists an $\eta>0$
such that for all $f \in \mathbb{O}_m$ with
$\|f-g\|_{\mathbb{A}} < \eta$ we have
\end{enumerate}
\vglue-20pt
\begin{equation}
\label{2fgkewi4}
\sum_{j=0}^{p-1}
  \frac{|\Delta_j|^{t}}{|\Delta_{j+1}|}
\le C \gamma^{mN} \ .
\end{equation}
\begin{enumerate}\vglue-19pt\phantom{up}
\ritem{(ii)} By Propostion \ref{zef} and Remark \ref{515151}, for every $  
\gamma > 1$ close
to one, there exists $t<2$ close to two and $C>0$ with the following
property.  For every  $g \in \mathbb{K}$ and every $m>0$, there exists
$\eta>0$ such that for all $f \in \mathbb{O}_m$ with
$\|f-g\|_{\mathbb{A}} < \eta$ we have that the inequality
\eqref{2fgkewi4} above is also satisfied.
\end{enumerate}
These last  estimates  end the proof of this proposition, provided  $\mu=\gamma^N$ and $\eta<
\mu^{m/\epsil}$.
\Endproof

We arrive at last to the two main  results of this section.

\vglue-22pt
\phantom{up}
\begin{theorem}
\label{666666} {\rm (i)} If $t > 2$ is not an integer{\rm ,} there exist $0 < \mu
< 1$ and $C>0$ with the following property. For every  $g \in
\mathbb{K}$ and for every $m>0${\rm ,} there is an $\eta>0$
such that for all $f \in \mathbb{O}_m$ with
$\|f-g\|_{\mathbb{A}} < \eta$  and for all $w \in \mathbb{A}^t$       
with
$\|w\|_{C^t} < \eta${\rm ,}
\begin{equation}
\label{pdeps}
\|DT^m(f+w)v-DT^m(f)v\|_{C^t}
   \le C  \mu^m \|v\|_{C^t} \ .
\end{equation}

{\rm (ii)} For every $  \mu > 1$ close to one{\rm ,} there exist $t<2$ close to $2$
and $C>0$ with the following property.
For every  $g \in \mathbb{K}$ and every $m>0${\rm ,} there exists $\eta>0$
such that for all $f \in \mathbb{O}_m$    with
$\|f-g\|_{\mathbb{A}} < \eta$ and all $w \in \mathbb{A}^t$      with
$\|w\|_{C^t} < \eta${\rm ,} the inequality \eqref{pdeps} above is also  
satisfied.
\end{theorem}

Part (ii) of this theorem with $t=s$ implies property $\mathbf{B5}$
and part (i) is used later (for $t=r$) to prove   property
$\mathbf{B6}$.

\Proof
In this proof the positive constants $K_1, K_2,\dots$
  depend only on $r$ and  $\mathbb{O}_m$  and also on $m$.
Let $E: C^t(I,I)\to \mathbb{R}$ be the evaluation map
$E(f)=f(0)$, which is linear, and let  $U_n:C^s(I,I)\to C^s(I)$  be as
before. Let us write
$
  DT^m(f+w)v- DT^m(f)(v)= E_1+E_2+E_3+E_4+E_5+E_6+E_7
$,
   where
\begin{eqnarray*}
E_1 & = &  ( A(f+w) -A(f) ) \sum_{j=0}^{p-1} B_j(f+w) \cdot C_j(f+w)\ ,  \\
E_2 & = &  A(f)  \sum_{j=0}^{p-1} (B_j(f+w)-B_j(f) ) \cdot C_j(f+w) \ ,  \\
E_3 & = &  A(f)  \sum_{j=0}^{p-1} B_j(f)  \cdot (C_j(f+w) -C_j(f))\ ,  \\
E_4 & = &  (A(f+w)  - A(f)) \cdot D(f+w)  \sum_{j=0}^{p-1} E \circ   
B_j(f+w) \cdot E \circ C_j(f+w)\ ,  \\
E_5 & = &  A(f)  \cdot (D(f+w)-D(f))  \sum_{j=0}^{p-1} E \circ   
B_j(f+w) \cdot E \circ C_j(f+w)\ ,   \\
E_6 & = &  A(f)  \cdot D(f)  \sum_{j=0}^{p-1} (E \circ  B_j(f+w) -  (E  
\circ  B_j(f))
               \cdot E \circ C_j(f+w) \ ,  \\
E_7 & = &  A(f)  \cdot D(f)  \sum_{j=0}^{p-1} E \circ  B_j(f)  \cdot
              ( E \circ C_j(f+w) - E \circ C_j(f) ) \ .  
\end{eqnarray*}
By Lemma \ref{lL},  
\begin{equation}\label{h1}
\left | A(f+w)-A(f)\right | =
\left |{\lambda_{f+w}}^{-1}- {\lambda_f}^{-1} \right |_{C^t} \leq
  K_1 \|w\|_{C^t} \ .
\end{equation}
Hence,
$$
\|E_1\|_{C^t}   \le    K _2 \|w\|_{C^t} \|v\|_{C^t} ~~~~~{\rm and}~~~~~
\|E_4\|_{C^t}   \le   K _3 \|w\|_{C^t} \|v\|_{C^t} \ .
$$
By Proposition \ref{rs} and Lemma \ref{fr23}, we obtain
\begin{eqnarray}\label{heew2}
\|B_j(f+w)-B_j(f)\|_{C^t}
& \le &
  K_4 \|U_{p-j}(f+w)-U_{p-j}(f)\|_{C^t}  \\
& \le &
  K_5 \|w\|_{C^t}  \ . \nonumber
\end{eqnarray}
Since   $E$ is a bounded linear operator, from the last inequality,  
we obtain
$$
\|E_2\|_{C^t}  \le   K _6 \|w\|_{C^t} \|v\|_{C^t}  ~~~{\rm and}~~~
\|E_6\|_{C^t}  \le
   K _7 \|w\|_{C^t}  \|v\|_{C^t} \ .
$$
Taking $j=p$ in   \eqref{heew2}, we get
$$
\|B_p(f+w)-B_p(f)\|_{C^t}
  \le  K_8 \|w\|_{C^t} \ .
$$
By   Lemma \ref{fr23} and by \eqref{h1}, we get
$$
\|\lambda_{f+w} \cdot U_p (f+w)- \lambda_f \cdot U_p (f)\|_{C^t}
  \le  K_9 \|w\|_{C^t} \ .
$$
Combining the last two inequalities, we get
$
\|E_5\|_{C^t}  \le
  K _{10} \|w\|_{C^t}  \|v\|_{C^t} \ .
$

Let  $k$ and $0<\alpha <1$ be the integer and the fractional part of  
$t$,
and let  $0<\epsil < 1$ be
such that $\alpha-\epsil > 0$.
 From  inequality \eqref{ikefp}, and since $E$ is a bounded linear  
operator,
we get
$$
|E \left( C_j(f+w) \right)-
E \left( C_j(f) \right)|
  \le  K_{11}   \|w\|_{C^t}^\epsil \|v\|_{C^t} \ .
$$
Thus,
$
\|E_7\|_{C^t}  \le
  K _{12} \|w\|_{C^t}^\epsil  \|v\|_{C^t} \ .
$
The only thing left to do is to   bound $\|E_3\|_{C^t}$, and  this
follows at once from Proposition  \ref{fr33333}.
\hfq

\begin{theorem}
\label{55555}
If $r > 2$ is not an integer{\rm ,} there exist $0 < \mu < 1$ and $C>0$ with
the following property. For every  $g \in \mathbb{K}$ and for every
$m>0${\rm ,} there is an $\eta>0$ such that for all $f \in \mathbb{O}_m$ with
$\|f-g\|_{\mathbb{A}} < \eta$ and for all $v \in \mathbb{A}^r$ with
$\|v\|_{C^r} < \eta${\rm ,}
\begin{equation}
\label{qwwwe}
\|T^m(f+v)-T^m(f)-DT^m(f)v\|_{C^r}
   \le  C \mu^m \|v\|_{C^r}  \ .
\end{equation}
\end{theorem}

This theorem together with Theorem \ref{666666} (i) for $t=r$ imply
that the renormalization operator satisfies  property $\mathbf{B6}$.

\Proof
In this proof the constants $\theta,\theta_1,\theta_2,\dots$ are
greater than zero and smaller than one and just depend upon $r$.
The positive constants $c,c_1,c_2,\dots$ depend only on $r$ and
$\mathbb{O}_m$, and
the positive constants  $K,K_1, K_2,\dots$  depend also on $m$.
Start by observing that since $T^m(f)=\lambda_f^{-1} \cdot U_p(f)$,
we have\break
$ T^m(f+v)-T^m(f)-DT^m(f)v=A+B+C
$, where
\begin{eqnarray*}
A & = & \lambda_f^{-1} \cdot \left( U_p(f+v)-U_p(f)-DU_p(f)v\right) \ ,\\
B &=& \left(\lambda_{f+v}^{-1}-\lambda_f^{-1}-D\lambda_f^{-1}  
(v)\right) \cdot U_p(f+v) \ , \\
C &=& D \lambda_f^{-1}(v)  \cdot \left( U_p(f+v)-U_p(f)\right) \ .
\end{eqnarray*}
By Lemma \ref{lL}, we have that  $f \to  \lambda_f^{-1}$ is $C^1$ and  
that there is $\theta_1$
such that
$\|B\|_{C^r} \le K_1 \|v\|_{C^r}^{1+\theta_1}$.
Since $\|U_p(f+v)-U_p(f)\|_{C^r} \le K_2 \|v\|_{C^r}$,
we have also $\|C\|_{C^r} \le  K_3 \|v\|_{C^r}^2$.
Hence inequality \eqref{qwwwe}   will be established if we  prove the  
following claim.

\demo{\scshape Claim}
{\it If $r > 2$ there exist $0 < \mu < 1$ and $c_1>0$ with the following  
property\/{\rm :}\/
for every  $g \in \mathbb{K}$ and for every $m${\rm ,} there is an $\eta>0$
such that for all $f \in \mathbb{O}_m$ with
$\|f-g\|_{\mathbb{A}} < \eta$ and for all $v \in \mathbb{A}^r$ with
$\|v\|_{C^r} < \eta$ we have}
\begin{equation}
\label{03}
   \|U_p(f+v)-U_p(f)-DU_p(f)v\|_{C^r}
  \le c_1  \mu^m |\lambda_f|  \|v\|_{C^r} \ .
\end{equation}
\Enddemo

To prove this claim, we will proceed recursively.
Let us write for $i=0,\dots,p$,
$$ R_i= U_i(f+v)-U_i(f)-DU_i(f)v   \ .$$
Note that
$ R_{i+1}= E_i +F_i + f' \circ U_i(f) \cdot  R_i$,
where
\begin{eqnarray*}
E_i &=& f \circ U_i(f+v)-f \circ U_i(f)-f' \circ U_i(f) \cdot  
\left(U_i(f+v)-U_i(f) \right) \ , \\
F_i &=&  v \circ U_i(f+v)- v \circ U_i(f)   \ .
\end{eqnarray*}
Thus, working recursively from these expressions, we get
$$
R_p =  R_0  \cdot  G_{p} +  \sum_{i=0}^{p-1} ( E_i \cdot G_{p-i-1} +   
F_i \cdot G_{p-i-1} ) \ ,
$$
where
$G_{p-i-1}   =  (f^{p-i-1})' \circ  U_{i+1}(f)$ and
$R_0 =  \Lambda_{f+v}-\Lambda_f - D\Lambda_f (v)$.
Since $f \in \mathbb{O}_m$, by Proposition \ref{rs} and Lemma  
\ref{fr23}, we get
$$ \|E_i\|_{C^r}   \le    K_4 \|U_i(f+v)-U_i(f)\|_{C^r}^{1+\theta_2}    
\le  K_5 \|v\|_{C^r}^{1+\theta_2}
$$
for $\theta_2=1-\{r\}$.
Therefore,
$
\left \|
\sum_{i=0}^{p-1} E_i \cdot G_{p-i-1}
\right \|_{C^r} \le   K_6 \|v\|_{C^r}^{1+\theta_2}
$.
By Lemma \ref{lL}, there is $\theta_3$ such that
$  \|R_0 \|_{C^r} \le K_7 \|v\|_{C^r} ^{1+\theta_3}$.
Hence, $ \|R_0  \cdot  G_{p} \|_{C^r} \le K_8 \|v\|_{C^r}  
^{1+\theta_3}$.
Finally, by Proposition \ref{fr33333},  there exists $\theta_4>0$ such  
that
$$
\left \|
\sum_{i=0}^{p-1} F_i \cdot G_{p-i-1}
\right \|_{C^r} \le   K_9 \|v\|_{C^r}^{1+\theta_4}   +
  c_2 \mu^m |\lambda_f| \|v\|_{C^r} \ .
$$
This proves our original claim.
\hfq

\vskip-5pt\Subsec{Proof of Theorem {\rm \ref{rob}}}
\label{sec:rob} All the pieces of the puzzle may now be put
together. We want to check the robustness of $T$ relative to the spaces
${\mathcal A}=\mathbb{A}$, ${\mathcal B}=\mathbb{A}^r $, ${\mathcal
C}=\mathbb{A}^s$ and ${\mathcal D}=\mathbb{A}^0$. By Theorem
\ref{dfgdsfgg}, the pair $(\mathbb{A}^\gamma,\mathbb{A}^0)$ is
$\rho_\gamma$-compatible with $(T,\mathbb K)$ and $\rho_\gamma <
\lambda$ for $\gamma$ sufficiently close to $2$ and  is
$1$-compatible for $\gamma>2$. Hence  property $\mathbf{B1}$ is
satisfied because  $s>s_0$ with $s_0<2$ close to $2$ and $r>s+1>2$.
Since $r>s+1$, we know from Lemma \ref{qn}  that  $T$ satisfies
property $\mathbf{B2}$. It also satisfies property $\mathbf{B3}$  by
Lemma \ref{333}, and property $\mathbf{B4}$ by Lemma \ref{4444}.
Finally,  $T$ satisfies property $\mathbf{B5}$ by Theorem
\ref{666666}, and property $\mathbf{B6}$ by Theorem \ref{55555}.
Therefore the renormalization operator $T$ is indeed robust with
respect to $(\mathbb{A}^r,\mathbb{A}^s,\mathbb{A}^0)$.

\vskip-5pt\Subsec{Proof of the hyperbolic picture}
\label{klrlrlll} Having established that the renormalization
operator $T$ is robust, we are now ready to show that the hyperbolic
picture holds true for $T$ acting on each of the spaces
$\mathbb{U}^r$ and $\mathbb{V}^r$.

\vskip-5pt\Subsubsec{Proof of Theorem {\rm \ref{main}}}
We divide the proof of Theorem \ref{main} into two cases: (i) $r$ is  
not an
integer (including the Lipschitz case $r=k+\text{Lip}$), and (ii) $r=k$  
is an
integer.

\demo{Proof of case {\rm (i)}}  Putting together Theorem
\ref{main3} with Theorem \ref{rob} we deduce all the assertions of  
Theorem
\ref{main} except the fact that the holonomies are $C^{1+\beta}$ for  
some
$\beta>0$. This last fact follows if we combine Theorem \ref{trfef} with
Example \ref{fggjggh}.

\demo{Proof of case {\rm (ii)}}   To prove this case, let us consider the
Banach space $\mathbb{A}^{k-1+\text{Lip}}$.
Note that the natural inclusion
$i:\mathbb{A}^k \to \mathbb{A}^{k-1+\text{Lip}}$
is an isometric embedding.
Indeed, for all $v \in \mathbb{A}^k$
we have $\|v\|_{C^k}=\|v\|_\mathcal{B}$,
by the mean-value theorem.
Applying case (i)  to $\mathbb{A}^{k-1+\text{Lip}}$,
we see that
for every $g \in \mathbb K$ the
local stable set $W_\varepsilon^{s,{k-1+\text{Lip}}}(g)$ is a  
codimension-one
$C^1$ Banach submanifold of $\mathbb{A}^{k-1+\text{Lip}}$.
In fact, there exists a $C^1$ function
$\Phi:\mathcal{O}_0 \to \mathbb{R}$, where
$\mathcal{O}_0 \subseteq \mathbb{A}^{k-1+\text{Lip}}$ is an open set  
containing $g$,
such that $0 \in \mathbb{R}$ is a regular value for $\Phi$,
with
$$
\Phi^{-1}(0)= \mathcal{O}_0 \cap W_\varepsilon^{s,{k-1+\text{Lip}}}(g)
$$
and such that $D\Phi(g) \mathbf{u}_g \ne 0$.
Let $\mathcal{O}_1=i^{-1}(\mathcal{O}_0) \subseteq \mathbb{A}^k$.
Then $\mathcal{O}_1$ is open and $\Phi \circ i: \mathcal{O}_1 \to  
\mathbb{R}$
is $C^1$.
Since $\mathbf{u}_g \in \mathbb{A}^k$ and
$D(\Phi \circ i)(g) \mathbf{u}_g = D \Phi(g) \mathbf{u}_g \ne 0$,
it follows that $0 \in \mathbb{R}$ is a regular value for
$\Phi \circ i$ at $g$.
Hence, by the implicit function theorem,
$$
\mathcal{O}_1 \cap W_\varepsilon^{s,k}(g)
= \mathcal{O}_1 \cap W_\varepsilon^{s,{k-1+\text{Lip}}}(g)
= \mathcal{O}_1 \cap (\Phi \circ i)^{-1}(0) \ ,
$$
is a $C^1$, codimension-one Banach submanifold of $\mathbb{A}^k$.
Since by case (i) the local stable manifolds in
$\mathbb{A}^{k-1+\text{Lip}}$ form a continuous lamination, we
deduce that the same is true for the local stable manifolds in
$\mathbb{A}^k$, because $i$ is an isometric embedding. Finally, if
$F$ is a $C^2$ ordered transversal (in the sense of \S
\ref{gholfff}) to the stable lamination in $\mathbb{U}^k$, then $i
\circ F$ is a $C^2$ ordered transversal to the stable lamination in
$\mathbb{A}^{k-1+\text{Lip}}$, and therefore by case (i)  its
holonomy in $\mathbb{U}^{k-1+\text{Lip}}$ is $C^{1+\theta}$ for some
$\theta >0$. But then it follows that the holonomy of the
transversal $F$ in $\mathbb{U}^k$ is $C^{1+\theta}$ also.

\Subsubsec{Proof of Corollary {\rm \ref{mainmain}}}
A similar argument to the one used in the proof of Theorem \ref{main}  
can be used
here. The map $i$ is replaced throughout by the inclusion  
$j:\mathbb{B}^r\to
\mathbb{A}^r$, which is a bounded linear operator (see \S  
\ref{quadunimodal}).
Hence the pre-images by $j$ of the local stable leaves in  
$\mathbb{U}^r$ are
$C^1$ manifolds and form a $C^0$ lamination in $\mathbb{V}^r$. Using
\cite{MMM} (see Remark \ref{referee} below), we see that the leaves of
such a lamination contain the local stable sets of each $g\in \mathbb{K}$  
in $\mathbb{V}^r$.

\section{Global stable manifolds and one-parameter families}
\label{sec:mainglobal}

In this section, we prove Theorem \ref{corol}. The first part will  
follow from
Theorem \ref{dfgdhkh} and the second part will follow from Theorem  
\ref{gfgfrtere1}.

\Subsec{The global stable manifolds of renormalization}
\label{sec:global}
In this section we construct the global stable manifolds
of the renormalization operator $T$ in $\mathbb{V}^r$,
for all $r$ sufficiently large.

Let $g$ be an element of the (bounded-type) invariant set $\mathbb{K}$  
of
$T$. Recall that the global stable set $W^{s,r}(g)$ of $g \in  
\mathbb{V}^r$
is given by
$$
W^{s,r}(g)=
\left \{ f \in \mathbb{V}^r: \|T^n(f)-T^n(g)\|_{C^r} \to 0 \
\, \text{ when}\ \,  n \to \infty \right \} \ .
$$
 From Corollary \ref{mainmain}, we know that the convergence is
exponential, and the exponential rate of convergence is independent of  
$f$ and
$g$, provided $r \ge 2+\alpha$ with $0 < \alpha < 1$ close to one.

\begin{theorem}
\label{dfgdhkh}
For every $r \ge 3+\alpha$ with $\alpha <1$ sufficiently close to $1${\rm ,}
and every $g \in \mathbb{K}${\rm ,}
the global stable set $W^{s,r}(g)$
is an immersed{\rm ,} codimension one $C^1$ Banach submanifold of  
$\mathbb{V}^r$.
\end{theorem}

\begin{remark}\label{referee}
By \cite{MMM}, if the invariant set $\mathbb{K}$
of the renormalization  operator is of bounded type then for every
$r \ge 3$ and every $g \in \mathbb{K}$ we have that
  $W^{s,r} (g)$
coincides with the set of all
maps $f \in \mathbb{V}^r$ with the same combinatorial type of~$g$.
\end{remark}

\Proof
We already know that the local stable sets are $C^1$ submanifolds.
The idea is to pull-back such  a manifold structure by $T$
using the implicit function theorem.
More precisely, by Corollary \ref{mainmain} there exist
$\varepsilon,\beta > 0$ so small  that
  $W^{s,r-1-\beta}_\varepsilon(g)$ is a
  codimension one $C^1$ Banach submanifold of
$\mathbb{V}^{r-1-\beta}$, for all $g \in \mathbb{K}$.
We may assume that $\varepsilon > 0$ is so small that
the vector $\mathbf{u}_g$ is transversal to the local
stable set $W^{s,r-1-\beta}_\varepsilon(g)$ at each one of its points.

Now fix $g \in \mathbb{K}$ and let $f \in  W^{s,r}(g)$.
There exists $N=N(f)>0$ so large that
$$
T^N(f)  \in W^{s,r}_\varepsilon(T^N(g) )
\subset W^{s,r-1-\beta}_\varepsilon(T^N(g) ) \ .
$$
Since $v=\mathbf{u}_{T^N(g)}$ is transversal at $T^N(f)$
to $W^{s,r-1-\beta}_\varepsilon(T^N(g) )$,
There exist a small open set $\mathbb{O}_0 \subset  
\mathbb{V}^{r-1-\beta}$
containing $T^N(f)$  and   a $C^1$ function
$\Phi:\mathbb{O}_0 \to  \mathbb{R}$ such that
$\Phi^{-1}(0)=W^{s,r-1-\beta}_\varepsilon(T^N(g)) \subset \mathbb{O}_0$
for which $0 \in \mathbb{R}$ is a regular value
and   $D\Phi (T^N(f)) v \ne 0$.
The operator $T^N$ is $C^1$ as a map from
$\mathbb{V}^r$ into $\mathbb{V}^{r-1-\beta}$.
Let $\mathbb{O}_1 \subset \mathbb{V}^r$
be an open set containing $f$ such that $T^N(\mathbb{O}_1) \subset  
\mathbb{O}_0$.
We want to show that $0 \in \mathbb{R}$ is a regular value for
$\Phi \circ T^N: \mathbb{O}_1 \to  \mathbb{R}$.
Defining $F_t=T^N(f)+tv$ (for $|t|$ small),
we get a   $C^1$ family $\{F_t\}$ of   maps in $\mathbb{V}^r$
which is transversal to $W^{s,r-1-\beta}_\varepsilon(T^N(g))$ at
$F_0=T^N(f)$. Now, we have the following claim.

\demo{Claim}
{\it   There exists a $C^1$ family $\{f_t\}$
  with $f_t \in \mathbb{V}^r$ such that for all small $t$
we have $T^N (f_t) =F_t$.}
\Enddemo

Let us assume this claim for a moment.
Setting
$$
w= \left. \frac{d}{dt} \right|_{t=0} f_t \ ,
$$
  we obtain that
$$
D(\Phi \circ T^N)(f) w
  = D \Phi  (F_0) v \ne 0 \ .
$$
Therefore, $\Phi \circ T^N$ is a $C^1$ local submersion at $f$.
By the implicit function theorem
$(\Phi \circ T^N)^{-1}(0)$ is a codimension one,
$C^1$ Banach submanifold of $\mathbb{O}_1$ (or $\mathbb{V}^r$).
Furthermore, if $h\in (\Phi \circ T^N)^{-1}(0)$ then
$T^N(h) \in W^{s,r-1-\beta}_\varepsilon(T^N(g) )$, and so $h$ belongs  
to the
global stable set $W^{s,r-1-\beta}(g)$. Using \cite{MMM} (see Remark
\ref{referee}), we deduce that $h$ belongs in fact to $W^{s,r}(g)$.  
This proves that
$W^{s,r}(g)$ is an immersed $C^1$ manifold as asserted.

It remains to prove the claim. We first note that $F_t=h_t \circ
F_0$ where each $h_t \in C^r(I,I)$ is a $C^r$ diffeomorphism of
$I=[-1,1]$. Since $T^N(f)=F_0$, there exist   $p>0$  and  closed,
pairwise disjoint intervals $0 \in \Delta_0, \Delta_1,\dots,
\Delta_{p-1} \subseteq I$ with  $f(\Delta_i) \subseteq \Delta_{i+1}$
for $0 \le i < p-1$ and $f (\Delta_{p-1}) \subseteq \Delta_0$, such
that
$$
F_0=T^N(f)=\Lambda_{f}^{-1} \circ f^p \circ \Lambda_{f} \ ,
$$
where $\Lambda_{f}:I \to \Delta_0$ is the map $x \mapsto f^p(0) x$.
Let $\overline{h}_t: \Delta_0 \to \Delta_0$ be the $C^r$
diffeomorphism given by $\overline{h}_t= \Lambda_{f} \circ h_t \circ
\Lambda^{-1}_{f}$. Consider a $C^r$ extension of $\overline{h}_t$ to
a diffeomorphism $H_t:I \to I$ with the property  that $H_t |
\Delta_i$ is the identity for all $i \ne 0$. Then let $f_t \in
\mathbb{V}^r$ be the map $f_t = H_t \circ f$. Note that $f_t^i(0) =
f^i(0)$ for all $0 \le i \le p$, that $f_t$ is $N$-times
renormalizable (under $T$) and that
\begin{eqnarray*}
T^N(f_t)  & = & \Lambda^{-1}_{f} \circ f_t^p \circ \Lambda_{f} \\
& = & \Lambda^{-1}_{f} \circ
\left(  H_t \circ  f   \right) |_{\Delta_{p-1}}  \circ
\left(  H_t \circ  f   \right) |_{\Delta_{p-2}}  \circ
\dots  \circ
\left(  H_t \circ  f   \right) |_{\Delta_0}
  \circ \Lambda_{f}  \\
& = &
  \Lambda^{-1}_{f} \circ  \overline{h}_t \circ f^p \circ \Lambda_{f}   =  
  \Lambda^{-1}_{f} \circ  \Lambda_{f}  \circ
  h_t \circ \left( \Lambda_{f}^{-1} \circ f^p \circ \Lambda_{f} \right)  
\\
& = &  h_t \circ F_0  
  =   F_t  \ ,
\end{eqnarray*}
which proves the claim.
\hfq

\vskip-5pt\Subsec{One-parameter families}
\label{jhgyttre}
A one-parameter family of maps is a map
$\psi:[0,1] \times I \to I$
(where $I=[-1,1]$ is the phase space) such that
$\psi_t = \psi(t, \cdot)$ belongs to $\mathbb{V}^r$ for all $t \in  
[0,1]$.
If $\psi$ is a $C^k$ map, then we say that $\psi$ is a $C^k$ family
(of $C^r$ unimodal maps).
We often identify the family $\psi$ with the curve
$\{\psi_t\}_{0 \le t \le 1}$ of unimodal maps in $\mathbb{V}^r$.
We shall denote by $\mathcal{UF}^k$ the space of all $C^k$ families  
with the
$C^k$ topology ($\mathcal{UF}^k$ is a subset of $C^k([0,1] \times I)$).

We say that two families are $C^{1+}$ equivalent if there exists a
diffeomorphism from one into the other which sends each
   infinitely renormalizable map (with a fixed bounded combinatorial  
type)
to a map with the same combinatorics.
We  are now in a position to state the result we have in mind.

\vglue-20pt \phantom{up}
\begin{theorem}
\label{gfgfrtere1}
Let $r \ge 3 + \alpha$ with $\alpha > 0$ close to $1${\rm ,}
and let $2 \le k \le r$. There exists an open and dense subset
$\mathcal{O} \subseteq \mathcal{UF}^k$
of one-parameter $C^k$ families of $C^r$
unimodal maps having the following properties\/{\rm :}\/

\vskip-25pt \phantom{up}
\begin{enumerate}
\item
Every family $\psi \in \mathcal{O}$ intersects the global stable  
lamination
$\mathcal{L}^s$ of renormalization transversally. 
\end{enumerate}
\vglue-36pt
\phantom{up}
\begin{enumerate}
 
\ritem{(ii)}
For every $\psi \in \mathcal{O}${\rm ,} there exist
$0=t_0 < t_1 < \dots < t_n = 1$ such that for each
$i=0,1,\dots,n-1$ the sub-arc
$\{\psi_t:t_i \le t \le t_{i+1} \}$ is $C^{1+\beta}$ diffeomorphic{\rm ,}
via a holonomy-preserving diffeomorphism{\rm ,}
to a corresponding sub-arc in the quadratic family.
Here $\beta>0$ is as  given in Corollary {\rm \ref{mainmain}.}
\end{enumerate}
\end{theorem}

The proof will require   a few lemmas.
The first lemma says that every $C^k$ family can be approximated
(in the $C^k$ sense) by a real analytic family.

\begin{lemma}
\label{fghhrttsd}
If $\psi \in \mathcal{UF}^k${\rm ,} then for each $\varepsilon > 0$
there exists a real analytic family $f \in \mathcal{UF}^\omega$
such that $\|\psi-f\|_{C^k([0,1] \times I)} < \varepsilon$.
\end{lemma}

\Proof
Write each $\psi_t \in \mathbb{V}^r$ as $h_t \circ q$, where
$q(x)=x^2$ and $h_t$ is a diffeomorphism, and consider the $C^k$ map
$h:[0,1] \times I \to I$ given by $h(t,x)=h_t(x)$, a $C^k$ family of
$C^r$ diffeomorphisms.
To approximate $h$ by a real analytic family of diffeomorphisms,
consider the convolution of $h$ with the heat kernel
$k(t,x,\varepsilon)=e^{-(t^2+x^2)/4 \varepsilon}$
for $\varepsilon > 0$ sufficiently small (see \cite{tangerman}).
\Endproof

Given this ``denseness'' result, the idea will be to show that  
arbitrarily
close to an $f$ as in Lemma \ref{fghhrttsd} we can find a $C^k$ family
which is also transversal to the global stable lamination  
$\mathcal{L}^s$
of renormalization, by some kind of perturbation argument,
to eliminate possible tangencies between $\{f_t\}$ and $\mathcal{L}^s$.

We will reduce our problem to the following general result about  
laminations
with  complex analytic leaves, whose
elegant proof is due to Douady.

\begin{lemma}
\label{fgtgtrhhh}
Let $\mathcal{L} \subseteq \mathbb{C}^2$ be a $C^0$ lamination
whose leaves are complex one-manifolds{\rm ,}
and let $F:\mathbb{D} \to \mathbb{C}$ be a holomorphic function whose  
graph
is tangent of finite order at $(0,F(0))$ to a leaf $L_0 \in  
\mathcal{L}$.
Then the tangency is isolated\/{\rm :}\/ there exists a neighborhood of $(0,F(0))$
in $\mathbb{C}^2$ on which every other
intersection of the graph of $F$ with the leaves of $\mathcal{L}$ is a  
transversal
intersection.
\end{lemma}

\Proof
Using a suitable chart, we may assume that the leaf $L_0$
is the horizontal plane $w=0$ in $\mathbb{C}^2$,
and that the other leaves of  $\mathcal{L}$ in that chart are the graphs
of holomorphic functions $\varphi_\mu:\mathbb{D} \to \mathbb{C}$
(with $\varphi_\mu(0)=\mu \in D$,
where $D \subseteq \mathbb{C}$ is some open disk around zero, and  
$\varphi_0 \equiv 0$).

Since $\mathcal{L}$ is a $C^0$ lamination,
$\varphi_\mu$ converges to $0$ uniformly in $\mathbb{D}$
as $\mu$ tends to $0$. Hence, for $|\mu|$ small enough,
we have $\varphi_\mu (\mathbb{D}) \subset \mathbb{D}$.
Moreover, $\varphi_\mu(z) \ne 0$ for all $z \in \mathbb{D}$
(leaves cannot intersect), so in fact
$\varphi_\mu (\mathbb{D}) \subset \mathbb{D}^*$.

Now, we have $F(0)=F'(0)=\dots=F^{(k-1)}(0) =0 \ne F^{(k)}(0)$,
for some $k \ge 2$.
Composing the chart with a bi-holomorphic map if necessary,
we may therefore assume that $F(z)=z^k$.

Let us fix $\mu \in D\setminus \{0\}$ and suppose that $z_0 \in  
\mathbb{D}$ is such that
$\varphi_\mu(z_0)=F(z_0)$.
We assume that $|z_0|< 1/2$ (taking $|\mu|$ small enough).
To show that this intersection between $\varphi_\mu$ and $F$
is transversal, it suffices to show that
$\varphi_\mu'(z_0) \ne F'(z_0)$.
But, by Schwarz's Lemma, the derivative $\varphi_\mu'(z_0)$
measured with respect to the Poincar\'e metrics of domain $\mathbb{D}$
and range $\mathbb{D}^*$ must be less than or equal to $1$, that is to  
say
$$
\left \|\varphi_\mu'(z_0)\right\|_P=
\frac
{\left|\varphi_\mu'(z_0)\right| \left(1-\left |z_0 \right |^2 \right) }
{\left|\varphi_\mu(z_0)\right|
\log \left( \left|\varphi_\mu(z_0)\right| ^{-1}   \right) }
\le 1 \ .
$$
Thus, 
$$
\left |\varphi_\mu'(z_0)\right| \le
\frac{4}{3} k |z_0|^k \log \left( \left|z_0 \right| ^{-1}   \right)
\ .
$$
On the other hand,
$$
\left|F'(z_0) \right|=k|z_0|^{k-1} \ .
$$
This shows that
$\left|\varphi_\mu'(z_0)\right| /
\left|F'(z_0) \right|$
converges to $0$ as $\mu$ tends to $0$,
whence $\varphi_\mu'(z_0) \ne F'(z_0)$ for all sufficiently
small $|\mu|$. Therefore $(0, F(0))$ is an isolated tangency as claimed.
\Endproof

We may now state and prove the result on  laminations with real  
analytic leaves
which is needed
for the proof of Theorem \ref{gfgfrtere1}.

\begin{lemma}
\label{hjkyhtgsef}
Let $\mathcal{F} \subseteq [a,b] \times \mathbb{R}$
be a $C^0$ foliation whose leaves are the graphs of real analytic  
functions
$\varphi_\mu:[a,b] \to \mathbb{R}$ with{\rm ,} say{\rm ,}
$\varphi_\mu(a)=\mu \in [0,1]$.
Let $\mathcal{L} \subseteq \mathcal{F}$ be a sub-lamination which is  
transversally
totally disconnected
  \/{\rm (}\/i.e.\ $K_0=\{\mu \in [0,1]: gr (\varphi_\mu)
  \subseteq \mathcal{L} \}$
is a totally disconnected set\/{\rm ).}\/
If $F:[a,b] \to \mathbb{R}$ is a real analytic function  whose graph is not a leaf{\rm ,} then
\begin{enumerate}
\item
${\rm gr}(F)$ is tangent to $\mathcal{F}$ at only finitely many points\/{\rm ;}\/
\item
for all $\varepsilon > 0$ and all $k \ge 0${\rm ,}
there exists a real analytic $G:[a,b] \to \mathbb{R}$
such that $\|F-G\|_{C^k} < \varepsilon$ and all tangencies of
${\rm gr}(G)$ with $\mathcal{F}$ belong to $\mathcal{F} \setminus  
\mathcal{L}${\rm ;}
in particular{\rm ,} ${\rm gr}(G)$ is transversal to $\mathcal{L}$.
\end{enumerate}
\end{lemma}

\Proof
(i) Complexifying $\mathcal{F}$ (i.e.\ the leaves $\varphi_\mu$)
as well as $F$, we put ourselves in the situation of Lemma  
\ref{fgtgtrhhh}.
All tangencies are therefore isolated, and since $[a,b]$ is compact,
there are only finitely many such, say at $x_i \in [a,b]$,
$i=1,2,\dots,n$.

\smallbreak
(ii) Let $d_i$ be the order of tangency of $F$ with $\mathcal{F}$
at $(x_i,F(x_i))$. Then for every real analytic $G$ sufficiently close
to $F$ in the $C^k$ topology with $k$ large
($k \ge \sum_{i=1}^n d_i$ will do), the number $n(G)$ of tangencies of
${\rm gr}(G)$ with $\mathcal{F}$ -- not counting multiplicies -- is bounded by
$\sum_{i=1}^n d_i$.
Hence we can find $G_0:[a,b] \to \mathbb{R}$ real analytic with
$\|F-G_0\|_{C^k} < \varepsilon / 2$ such that $n(G_0)$ is maximal.
All tangencies of $G_0$ with $\mathcal{F}$ must be first-order  
tangencies
($d_i=1$).
Indeed, if, say, $d_1 > 1$, then adding a suitable polynomial with small
$C^k$ norm to $G_0$, vanishing of very high order at
$x_2, x_3, \dots,x_{n(G_0)}$, we could unfold the tangency at $x_1$
to produce a new real analytic $G$ with $n(G) > n(G_0)$.
Now we may consider $G_t:[a,b] \to \mathbb{R}$ given by
$G_t(x)=G_0(x)+t$ for $|t| < \varepsilon/2$.
Since first-order tangencies are persistent, each tangency
$(x_i, G_0(x_i))$ of $G_0$ with $\mathcal{F}$ generates a continuous,
nonconstant path $(x_i(t), G_t(x_i(t))) \in {\rm gr}(\varphi_{\mu_i(t)})$
of (first-order) tangencies of $G_t$ with $\mathcal{F}$.
Each function $t \mapsto \mu_i(t)$, $i=1,2,\dots,n(G_0)$,
is continuous and nonconstant.
Since $K_0$ is totally disconnected,
there exists $t$ (with $|t| < \varepsilon/2$) such that
$\mu_i(t) \in [0,1] \setminus K_0$ for all $i$.
Therefore, all tangencies of $G_t$ with $\mathcal{F}$ fall in
$\mathcal{F} \setminus \mathcal{L}$,
whence $G_t$ is transversal to $\mathcal{L}$.
\hfq

\demo{Proof of Theorem {\rm \ref{gfgfrtere1}}} 
Both properties (i) and (ii) are easily seen to be open,
hence we concentrate on proving that they are dense.
Let $\varepsilon>0$.

Take any family $\psi \in \mathcal{UF}^k$.
By Lemma \ref{fghhrttsd}, there exists a real analytic family
$f \in \mathcal{UF}^\omega$ whose $C^k$ distance from $\psi$ is
less than $\varepsilon/2$.
The corresponding curve $\{f_t\}$ in $\mathbb{V}^r$ may fail to be  
transversal
to the global stable lamination $\mathcal{L}^s$, so let us show how
to perturb it locally to get a transversal family.
Let $t_0 \in [0,1]$ be such that $f_{t_0} \in \mathcal{L}^s$
(and $\{f_t\}$ is tangent to $\mathcal{L}^s$ at $f_{t_0}$).
Since $f_{t_0}$ is infinitely renormalizable and real analytic,
there exists $N>0$ such that $R^N(f_{t_0}) \in \mathbb{A}_{\Omega_a}$
(where $a>0$ is the constant in Theorem \ref{lyubhyp}).
Let $J \subseteq  [0,1]$ be an interval containing $t_0$ such that
$R^N(f_t)$ is well-defined and belongs to $\mathbb{A}_{\Omega_a}$
for all $t \in J$. We restrict our attention to the
sub-family $\{f_t\}_{t \in J}$ from now on.

First we embed $\{f_t\}_{t \in J}$ in a two-parameter family in the  
following
way. Note that each $f_t$ belongs to $\mathbb{A}_{\Omega_\alpha}$
for some (fixed) $\alpha > 0$. As a map from (an open subset of)
$\mathbb{A}_{\Omega_\alpha}$ into $\mathbb{A}_{\Omega_a}$,
$R^N$ is a real analytic operator.

\demo{\scshape Claim}
{\it There exist analytic vectors $v \in \mathbb{A}_{\Omega_\alpha}$
and $w \in \mathbb{A}_{\Omega_a}$ with the property that
$DR^N(f_{t_0}) v = w$ and $w$ is transversal to
$\mathcal{L}^s_a = \mathcal{L}^s \cap \mathbb{A}_{\Omega_a}$ at
$R^N(f_{t_0}) \in \mathcal{L}^s_a$.}
\Enddemo

To see this, take any $w_0 \in \mathbb{A}_{\Omega_a}$ transversal
to the (co-dimension one) lamination $\mathcal{L}^s_a$ at  
$R^N(f_{t_0})$.
The same construction used in the proof of Theorem \ref{dfgdhkh}
yields a $C^\infty$ vector $v_0$ at $f_{t_0}$ such that
$DR^N(f_{t_0})v_0=w_0$. Now approximate $v_0$ by an analytic vector
$v \in \mathbb{A}_{\Omega_\alpha}$
(in the $C^m$ sense for $m \ge r$).
Then $w=DR^N(f_0) v$ will still be transversal to $\mathcal{L}_a^s$.
Shrinking $J$ if necessary, we may in fact assume that
$DR^N(f_t)v$ is transversal to $\mathcal{L}_a^s$ for all $t \in J$.
Hence, we consider the two-parameter family of maps
$f_{t,s} \in  \mathbb{A}_{\Omega_\alpha}$ given by
$f_{t,s}=f_t+s v$ with $t \in J$ and $|s| \le \delta$ with $\delta$  
small.
Now,
$$
W=\left\{f_{t,s}:\ t \in J, \ s \in [-\delta,\delta] \right\} \cong
J \times [-\delta,\delta] \subseteq \mathbb{R}^2 \ ,
$$
and $R^N|W:W \to A_{\Omega_a}$ is an injective, real analytic map.
Recall now that in $\mathbb{A}_{\Omega_a}$ we have a  $C^0$ foliation
$\mathcal{F}$ with real analytic leaves
(coming from hybrid classes, {\it cf.} \S \ref{grtsfd})
and that $\mathcal{L}_a^s \subseteq \mathcal{F}$
is the sub-lamination corresponding to the stable leaves of  
renormalization,
which is transversally totally disconnected.
Taking $\mathcal{F}_W=R^{-N}(\mathcal{F}) \subseteq W$
and $\mathcal{L}^s_W=R^{-N}(\mathcal{L}_a^s) \subseteq W$
and noting that
$DR^N(f_{t,s}) v = w$ is transversal to $\mathcal{L}_a^s$ for all
$t \in J$, $s \in [-\delta,\delta]$
(making $\delta$ smaller if necessary) we deduce that
$\mathcal{F}_W$ is a $C^0$ foliation (in $W$) by real analytic curves,
and $\mathcal{L}_W^s \subseteq \mathcal{F}_W$ is a sub-lamination.
Therefore we can apply Lemma~\ref{hjkyhtgsef} to this situation
(with $\mathcal{F}= \mathcal{F}_W$ and $\mathcal{L}= \mathcal{L}_W^s$),
obtaining a new analytic curve
$\{g_t\}_{t \in J}$ with $\|f_t-g_t\|_{C^k} < \varepsilon/2$,
transversal to $\mathcal{L}_W^s$ in $W$,
and such that $\{R^N(g_t)\}$ is transversal to
$\mathcal{L}_a^s$ at $R^N(g_{t_0})$. Since by Corollary \ref{mainmain}  
the
holonomy of $\mathcal{L}_a^s$ is $C^{1+\beta}$ for some $\beta>0$
(and the quadratic family is transversal to $\mathcal{L}_a^s$)
we deduce that $\{g_t\}$ satisfies properties (a) and (b) of the  
statement.
This completes the proof.
\hfill\qed

\section{A short list of symbols}

For the reader's convenience, we present below a short list of symbols  
used in
this paper.

\begin{table}[ht]
\renewcommand\arraystretch{1.5}
\noindent\[
\begin{array}{||c|c||}
\hline \hline
{p}& {\text{Period of renormalization}}\\
\hline
{\lambda_f}& {\text{Scaling factor}\  \lambda_f=f^p(0)}\\
\hline
{\Lambda_f} & {\text{Linear scaling}\ \Lambda_f:x \to f^p(0)
\cdot
x}\\
\hline
{R} & {\text{Renormalization operator}\  R^N f=
\Lambda_f^{-1} \circ
f^p \circ \Lambda_f} \\
\hline
{\mathbb{K}}& {\text{Bounded type limit set of}\ R}\\
\hline
{p_k}&{\text{Number of renormalization intervals at level}\ k}\\
\hline
{\Delta_{j,k}(f)} &{\text{ Renormalization intervals at
level}\ k\
(0\leq j\leq p_k-1)}\\
\hline
{\mathcal{I}_f}&{\text{Post-critical set of}\ f}\\
\hline
 {\mathbb{A}_V}&\begin{array}{c}{\text{Real Banach space of continuous maps}\
\overline{V} \to \mathbb{C},}\\[-6pt]
 {\text{holomorphic in}\ V,
\ \text{symmetric about real axis}} \end{array}\\
\hline 
{T=R^N:\mathbb{O} \to \mathbb{A}_{\Omega_a}}&
\begin{array}{c}{\text{Real
analytic operator for which}}\\[-6pt]
 { \mathbb{K}\subset  \mathbb{O}\ \text{is  a
hyperbolic basic set}}\end{array}\\
\hline {u_g(t)}&{\text{Parametrization of local unstable manifold}\
W^u_\varepsilon(g)}\\
\hline
{\mathbf{u}_g}&{\text{Unit vector tangent to}\
W^u_\varepsilon(g)\
\text{at}\ g}\\
\hline
{\delta_g}&{\text{Unique real number such that}\
DT(g)\mathbf{u}_g=\delta_g \mathbf{u}_{T(g)}}\\
\hline {\delta_g^{(n)}}&{\text{The product}\  \delta_g \delta_{T(g)}
\dots \delta_{T^{n-1}(g)}}\\
\hline
{L_f=DT(f)}&{\text{Derivative of}\  T\ \text{at}\ f}\\
\hline
{\mathbb{V}^r}&{C^r\  \text{unimodal maps with quadratic
critical
point at}\ 0}\\
\hline {\mathbb{A}^r}&{\text{Tangent space to unimodal maps
contained in}\ \mathbb{V}^r}\\ \hline\hline
\end{array}
\]
\end{table}

\newpage


\end{document}